\documentclass[a4paper,12pt]{article}
\usepackage[english]{babel}
\usepackage[utf8]{inputenc}
\usepackage{amsmath}
\usepackage{amsthm}
\usepackage{amssymb}
\usepackage{ytableau}
\usepackage{youngtab}
\usepackage{color}
\usepackage{xcolor}
\usepackage{graphicx}
\usepackage{tikz}
\usepackage{mathrsfs}
\usepackage{setspace}
\usepackage{enumitem}
\usepackage{aliascnt}
\usepackage{etoolbox}
\usepackage{rotating}
\usepackage[margin=20mm]{geometry}
\usepackage[colorlinks=true]{hyperref}

\setlist[enumerate]{label=(\alph*)}

\def\NewTheorem#1{
	\newaliascnt{#1}{equation}
	\newtheorem{#1}[#1]{#1}
	\aliascntresetthe{#1}
	\expandafter\def\csname #1autorefname\endcsname{#1}
}
\def\equationautorefname~#1\null{(#1)\null}
\def\itemautorefname~#1\null{(#1)\null}

\numberwithin{equation}{section}
\swapnumbers

\newtheorem*{main theorem*}{Main Theorem}
\newtheorem*{theorem*}{Theorem}
\newtheorem*{proposition*}{Proposition}

\NewTheorem{Proposition}
\NewTheorem{Theorem}
\NewTheorem{Main Theorem}
\NewTheorem{Corollary}
\NewTheorem{Lemma}
\NewTheorem{Definition}
\theoremstyle{definition}
\NewTheorem{Example}
\AtEndEnvironment{Example}{\null\hfill$\Diamond$}%
\NewTheorem{Examples}
\AtEndEnvironment{Examples}{\null\hfill$\Diamond$}%
\theoremstyle{remark}
\NewTheorem{Remark}

\newcommand{\residue}{\mathop{\rm res}}

\newcommand{\Std}{\mathop{\rm Std}}
\newcommand{\SStd}{\mathop{\rm SStd}}
\newcommand{\radic}{\mathop{\rm Rad}}
\newcommand{\shape}{\mathop{\rm shape}}
\newcommand{\toppo}{\mathop{\rm top}}
\newcommand{\bottom}{\mathop{\rm bot}}

\usetikzlibrary{patterns,decorations.pathreplacing}

\tikzset{
	ultra thin/.style= {line width=0.05pt},
	very thin/.style=  {line width=0.2pt},
	thin/.style=       {line width=0.1pt},
	semithick/.style=  {line width=0.6pt},
	thick/.style=      {line width=0.8pt},
	very thick/.style= {line width=1.2pt},
	ultra thick/.style={line width=1.6pt}
}

\tikzset{wei/.style={red,double=red,double
		distance=1pt}}

\tikzset{wei2/.style={red,double=red,double
		distance=1pt}}

\tikzset{ghost/.style={gray,densely dashed}}

\tikzset{
	pics/solidstring/.style args={#1,#2}{
		code = {
			\foreach \pt [count=\c, remember=\pt as \Pt] in {#2} {
				\ifnum\c>1
				\draw \Pt--\pt;
				\draw[ghost, xshift=-#1]\Pt--\pt;
				\fi
			}
		}
	}
}

\let\emph\relax
\DeclareTextFontCommand{\emph}{\bfseries\em}

\tikzset{
	centered/.style = {
		baseline = {([yshift=1.0ex]current bounding box.center)}
	},
	dots/.style={line width=1pt, line cap=round,
		dash pattern=on 0pt off 2\pgflinewidth},
	frame/.style={fill=FloralWhite, draw=BurlyWood,very thin},
	circled/.style = {fill=Tan},
	redstring/.style = {line width=1.0mm, red},
	semisolid/.style = {draw=DarkBlue, thin},
	pics/dot/.style = {
		code={\node[fill=DarkBlue,circle,inner sep=1.8pt,outer sep=0]at(0,0){};}
	},
	wall/.style={black!80, thick},
	innerwall/.style={gray, thin, fill=white},
	pics/russiantableau/.style 2 args={
		code = { \begin{scope}[#1]
				\def\lastC{0}
				\foreach \row [count=\r, remember=\r as \rr] in {#2} {
					\foreach \ent [count=\c, remember=\c as \cc] in \row {
						\draw[innerwall] (\r-\c,\c+\r-2) --++(-1,1)--++(1,1)--++(1,-1)--cycle;
						\node at (\r-\c,\c+\r-1) {\ent};
					}
					\ifnum\r=1\draw[wall](0,0)--(-\cc,\cc);\fi
					\ifnum\lastC>\cc\draw[wall](\r-\cc-1,\cc+\r-1)--(\r-\lastC-1,\lastC+\r-1);\fi
					\draw[wall](\r-1,\r-1)--++(1,1);
					\draw[wall](\r-\cc-1,\cc+\r-1)--++(1,1);
					\xdef\lastC{\cc}
				}
				\draw[wall](\rr,\rr)--++(-\lastC,\lastC);
			\end{scope}
		}
	}
}

\title{On simple modules of cyclotomic quiver Hecke algebras of type A}
\author{Alexander Ferdinand Kerschl}
\date{}

\begin{document}

\ytableausetup{centertableaux}

\maketitle

\abstractname: We give a classification of the graded simple modules of cyclotomic quiver Hecke algebras of type A using the diagram calculus of the diagrammatic Cherednik algebra. We also obtain a non-trivial lower bound for the dimension of the simple modules which is independent of the field characteristic.

\section{Introduction}

Khovanov and Lauda \cite{KhovLaud:diagI} and, independently, Rouquier \cite{Rouq:2KM} introduced a family of diagrammatic algebras that depend on a quiver. These algebras are now known as quiver Hecke or KLR algebras. They are widely studied because they categorify the negative half of the quantum group $U_q^{-}(\mathfrak{g})$, where $\mathfrak{g}$ is the Kac-Moody Lie algebra associated with the quiver. Recently Webster \cite{Webster:RouquierConjecture} constructed diagrammatic algebras depending on a parameter $\theta$, which are generalisations of the KLR algebras. Hence, he called them weighted KLR algebras but they are also known as diagrammatic Cherednik algebras.

For quivers of type A, each quiver Hecke algebra has a cyclotomic quotient $\mathscr{R}_n^\Lambda$ that is isomorphic to the cyclotomic Ariki-Koike Hecke algebra $\mathscr{H}_n^\Lambda(q)$, by Brundan and Kleshchev \cite{BK:GradedKL}. Ariki \cite{Ariki:class} showed that the simple modules of $\mathscr{H}_n^\Lambda(q)$ are indexed by a certain set $\mathscr{K}^\ell$ of multipartitions and Brundan and Kleshchev \cite{BK:GradedDecomp} used the isomorphism to prove that the same holds for the simple modules of $\mathscr{R}_n^\Lambda$. Webster \cite{Webster:RouquierConjecture}, on the other hand, showed that the diagrammatic Cherednik algebra contains a subalgebra that is isomorphic to $\mathscr{R}_n^\Lambda$ independently of the choice of $\theta$. He \cite{Webster:RouquierConjecture} and Bowman \cite{Bowman:ManyCellular} used this isomorphism to construct a family of different graded cellular bases, each depending on $\theta$, for $\mathscr{R}_n^\Lambda$. This gives $\mathscr{R}_n^\Lambda$ the structure of a graded cellular algebra, as in \cite{HuMathas:GradedCellular}. The corresponding cell modules $S_\theta^\lambda$ of $\mathscr{R}_n^\Lambda$ are indexed by multipartitions $\lambda\in\mathscr{P}^\ell(\theta)$ and they come equipped with an inner product. As in \cite{GL} or \cite{HuMathas:GradedCellular}, one can then find all simple modules of $\mathscr{R}_n^\Lambda$ by looking at the list of non-zero quotient modules $D_\theta^\lambda = S_\theta^\lambda/\radic(S_\theta^\lambda)$, where $\radic(S_\theta^\lambda)$ is the radical of the inner product on $S_\theta^\lambda$. It is a non-trivial problem to determine when $D_\theta^\lambda$ is non-zero.

For each $\theta$ there is a crystal graph that has a recursively defined set $\mathscr{U}^\ell(\theta)$ of multipartitions as vertices. Jacon \cite{Jacon:Staggered} gives an alternative description of the set $\mathscr{U}^\ell(\theta)$ and allows us to construct special paths $\textbf{d}$ in a directed graph for each multipartition in $\mathscr{U}^\ell(\theta)$. We can now state our main result:
\begin{main theorem*}
	Let $\lambda\in\mathscr{P}_n^\ell(\theta)$. Then the following are equivalent:
	\begin{enumerate}
		\item $D_\theta^\lambda\neq 0$
		\item $\lambda\in\mathscr{U}_n^\ell(\theta)$
		\item $|\langle c_{\mathfrak{t}_{\textbf{d}}}^\lambda, c_{\mathfrak{t}_{\textbf{d}}^\diamond}^\lambda\rangle_\lambda| = 1$, for any decomposition path $\textbf{d}$ for $\lambda$
	\end{enumerate}
\end{main theorem*}
In the special case when $\theta$ is 'well-separated' (see after \autoref{D: Dominance order}), $\mathscr{K}^\ell = \mathscr{U}^\ell(\theta)$ so this recovers the classification of the simple $\mathscr{R}_n^\Lambda$-modules given by Brundan and Kleshchev \cite{BK:GradedDecomp} and Ariki \cite{Ariki:class}. In this case we also obtain a result that implies a generalisation of the Dipper-James-Murphy conjecture \cite[Conjecture 8.13]{DJM} (see also \cite[Question 5.10 (i)]{GL}):
\begin{theorem*}
	Let $\lambda\in\mathscr{P}_n^\ell(\theta)$ and $\theta$ well-separated. Then $D_\theta^\lambda \neq 0$ if and only if there exists a DJM-tableau of shape $\lambda$.
\end{theorem*}
Ariki and Jacon \cite{ArikiJacon:DJM} proved the original conjecture in type B. The generalised conjecture was recently proved by Jacon \cite{Jacon:Staggered}. However, our result gives us an alternative proof and more importantly, this result and the application of our constructions also give us a non-trivial lower bound for the graded dimension of the graded simple module $D_\theta^\lambda$ which is independent of the characteristic of the field. And in the special case when $\theta$ is well-separated we actually get the dimensions of certain pieces of the simple module:
\begin{theorem*}
	Let $\lambda\in\mathscr{U}_n^\ell(\theta)$ and $\mathbf{d} = \mathbf{s}_1\ldots \mathbf{s}_k$ a decomposition path for $\lambda$ with $|\mathbf{s}_j| = m_j$. Then
	\begin{equation*}
	\dim_t(D_\theta^\lambda \cdot e_\omega(\residue(\mathfrak{t}_{\mathbf{d}}^\diamond))) \geq \prod_{j=1}^{k} [m_j]_t^!
	\end{equation*}
	Moreover, if $\lambda\in\mathscr{K}_n^\ell$ then we have equality.
\end{theorem*}

The paper is structured into 4 sections. In section 2 we start by recalling the theory of graded cellular algebras following \cite{GL, HuMathas:GradedCellular} and by defining the diagrammatic Cherednik algebra with the necessary combinatorics and diagrams, as in \cite{Webster:RouquierConjecture, Bowman:ManyCellular}. For section 3 we define the set $\mathscr{U}^\ell(\theta)$ and recall Jacon's \cite{Jacon:Staggered} recent alternative description using staggered sequences. These allow us to construct the decomposition paths $\textbf{d}$ and a family of standard tableaux associated to it. In the key theorem we use the basis elements in $S_\theta^\lambda$ associated to these tableaux to compute the inner product, which implies that $D_\theta^\lambda\neq 0$ if $\lambda\in\mathscr{U}^\ell(\theta)$. This and consequences of this theorem are part of section 4. Section 5 is mostly dedicated to prove the key theorem using diagram calculations.

\vfill
\noindent\hrulefill\\
This paper is written based on Version~4 of Bowman's paper \cite{Bowman:ManyCellular}, as available from arXiv.org. Just prior to submitting this paper the author learned that Bowman has recently updated his paper on the ArXiv. The most recent version of his paper uses substantially different notation from the previous versions. This paper uses the notation from version 4 of~\cite{Bowman:ManyCellular}. In addition, the latest version of Bowman's paper gives a classification of the simple $\mathscr{R}_n^\Lambda(\theta)$-modules, which gives an independent proof of parts of \autoref{T: Main Theorem}.

\section{Diagrammatic Cherednik Algebra}

This section introduces the diagrammatic Cherednik algebras, which are the main objects of study in this paper. The diagrammatic Cherednik Algebras contain a cyclotomic quiver Hecke or KLR Algebra of type A as a subalgebra. Both of these algebras come with bases that give them the structure of a graded cellular algebra. The representation theory of graded cellular algebras was established by Hu and Mathas \cite{HuMathas:GradedCellular}, extending Graham and Lehrer's work in the ungraded case \cite{GL}. As we recall, cellular algebras come equipped with a family of cell modules and these can be used to construct all simple modules.

\subsection{Graded cellular algebras}

In this paper graded will always mean $\mathbb{Z}$-graded. Throughout let $K$ be a field of any characteristic. A \emph{graded algebra} $A$ is a unital associative $K$-algebra such that
\begin{equation*}
A = \bigoplus_{j\in\mathbb{Z}} A_j \quad\text{ and }\quad A_x A_y \subseteq A_{x+y},\text{ for all } x,y\in\mathbb{Z}.
\end{equation*}
Elements $a\in A_x$, for some $x\in\mathbb{Z}$, are \emph{homogeneous} of degree $x$. Set $\deg(a) = x$. Let $\underline{A}$ be the ungraded algebra obtained by forgetting the grading.

A graded (right) $A$-module $M$, similarly, has a direct sum decomposition $M = \bigoplus_{j\in\mathbb{Z}} M_j$ and must satisfy $\underline{M}$ being an $\underline{A}$-module and $M_x A_y \subseteq M_{x+y}$. Given a graded $A$-module $M$ let $M\langle x\rangle$ be the graded $A$-module that is obtained by shifting the grading on $M$ up by $x$, i.e. $(M\langle x\rangle)_y = M_{y-x}$. The graded version of submodules, left modules, and so on are defined accordingly.

Let $t$ be an indeterminate over $\mathbb{N}_0$. If $M$ is a finite dimensional graded $A$-module such that $M_j$ is finite dimensional over $K$, then its \emph{graded dimension} is the Laurent polynomial $\dim_t(M) = \sum_{j\in\mathbb{Z}} \dim(M_j) t^j\in\mathbb{N}_0[t,t^{-1}]$.

\begin{Definition}[{\cite[Definition 1.1]{GL}},{\cite[Definition 2.1]{HuMathas:GradedCellular}}]\label{D: Graded cellular}
	Let A be a finite dimensional graded $K$-algebra. A \emph{graded cell datum} for $A$ is an ordered quadruple $(\mathscr{P}, T, c, \deg)$, where $(\mathscr{P}, \triangleright)$ is a poset, $T(\lambda)$ a finite set for $\lambda\in\mathscr{P}$, and
	\begin{equation*}
	c: \bigcup\limits_{\lambda\in\mathscr{P}} T(\lambda) \times T(\lambda) \longrightarrow A: (s,t)\mapsto c_{st}^\lambda, \quad\text{ and }\quad \deg:\bigcup\limits_{\lambda\in\mathscr{P}} T(\lambda) \longrightarrow \mathbb{Z}
	\end{equation*}
	are two functions such that $c$ is injective and the following hold:
	\begin{enumerate}
		\item The set $\{c_{st}^\lambda|\lambda\in\mathscr{P},s,t\in T(\lambda)\}$ is a $K$-basis of $A$.
		\item\label{D: Graded cellular b} Let $\lambda\in\mathscr{P}$, for any $s,t\in T(\lambda)$ and $a\in A$ there exist scalars $r_{tv}(a)\in K$ not depending on $s$ such that
		\begin{equation*}
		c_{st}^\lambda a \equiv \sum_{v\in T(\lambda)} r_{tv}(a) c_{sv}^\lambda \pmod {A^{\triangleright\lambda}},
		\end{equation*}
		where $A^{\triangleright\lambda}$ is the $K$-submodule of $A$ spanned by $\{c_{ab}^\mu | \mu\triangleright\lambda \text{ and } a,b\in T(\mu)\}$.
		\item There is a $K$-linear anti-isomorphism $\ast: A \longrightarrow A$ determined by $(c_{st}^\lambda)^\ast = c_{ts}^\lambda$, for any $\lambda\in\mathscr{P}$ and $s,t\in T(\lambda)$.
		\item Each element $c_{st}^\lambda$ is homogeneous of degree $\deg(c_{st}^\lambda) = \deg(s) + \deg(t)$, for $\lambda\in\mathscr{P}$ and $s,t\in T(\lambda)$.
	\end{enumerate}
	A \emph{graded cellular algebra} is a graded algebra with a graded cell datum. The basis $\{c_{st}^\lambda\}$ is a \emph{graded cellular basis} of $A$.
\end{Definition}

Given a graded cellular basis for an algebra $A$ we can canonically define special modules using property \ref{D: Graded cellular b} of the definition.

\begin{Definition}
	Let $A$ be a graded cellular algebra and let $\lambda\in\mathscr{P}$. Then the \emph{graded cell module} $C^\lambda$ is the graded right $A$-module
	\begin{equation*}
	C^\lambda = \bigoplus_{j\in\mathbb{Z}} C_j^\lambda,
	\end{equation*}
	where $C_j^\lambda$ is the $K$-span of $\{c_t^\lambda | t\in T(\lambda) \text{ and } \deg(t) = j\}$ and where the $A$-action is given by
	\begin{equation*}
	c_{t}^\lambda a = \sum_{v\in T(\lambda)} r_{tv}(a) c_{v}^\lambda,
	\end{equation*}
	and the scalars $r_{tv}(a)$ are the same as in \autoref{D: Graded cellular}\ref{D: Graded cellular b}.
\end{Definition}

As in \cite[After Definition 1.1]{GL}, applying $\ast$ to \autoref{D: Graded cellular}\ref{D: Graded cellular b}:
\begin{enumerate}[label = (\alph*')]
	\setcounter{enumi}{1}
	\item\label{D: Graded cellular b'} Let $\lambda\in\mathscr{P}$, for any $s,t\in T(\lambda)$ and $a\in A$ there exist scalars $r_{vs}(a)\in K$ not depending on the choice of $t$ such that
	\begin{equation*}
	a c_{st}^\lambda \equiv \sum_{v\in T(\lambda)} r_{vs}(a) c_{vt}^\lambda \pmod {A^{\triangleright\lambda}},
	\end{equation*}
	where $A^{\triangleright\lambda}$ is the $K$-submodule of $A$ spanned by $\{c_{ab}^\mu | \mu\triangleright\lambda \text{ and } a,b\in T(\mu)\}$.
\end{enumerate}

By \ref{D: Graded cellular b} and \ref{D: Graded cellular b'}, there is a well-defined inner product on $C^\lambda$ that is uniquely determined by
\begin{equation}\label{eq: def inner product}
c_{st}^\lambda c_{uv}^\lambda \equiv \langle c_t^\lambda , c_v^\lambda \rangle_\lambda c_{sv}^\lambda \pmod {A^{\triangleright\lambda}}.
\end{equation}
Equivalently $c_s^\lambda c_{tu}^\lambda = \langle c_s^\lambda, c_t^\lambda \rangle_\lambda c_u^\lambda$.

\begin{Proposition}[{\cite[Lemmata 2.6, 2.7]{HuMathas:GradedCellular}}]\label{P: Inner product properties}
	Let $\lambda\in\mathscr{P}$	and $a\in A$, $x,y\in C^\lambda$. Then $\langle\ , \ \rangle_\lambda$ is symmetric, homogeneous of degree $0$, and
	\begin{equation*}
	\langle xa, y \rangle_\lambda = \langle x, ya^\ast \rangle_\lambda.
	\end{equation*}
\end{Proposition}

It follows that $\radic(C^\lambda) = \{x\in C^\lambda | \langle x, y \rangle_\lambda = 0, \text{ for all } y\in C^\lambda\}$ is a graded submodule of $C^\lambda$.

\begin{Definition}\label{D: Simple}
	Suppose $\lambda\in\mathscr{P}$. Let $D^\lambda := C^\lambda / \radic(C^\lambda)$.
\end{Definition}

Define the subset $\mathscr{P}_0 = \{\lambda\in\mathscr{P} | D^\lambda \neq 0\} \subseteq \mathscr{P}$.

\begin{Theorem}[{\cite[Theorem 3.4]{GL}, \cite[Theorem 2.10]{HuMathas:GradedCellular}}]\label{T: Classification Cell Alg}
	Let $\lambda,\mu\in\mathscr{P}_0$ and $K$ be any field. Then
	\begin{enumerate}
		\item $D^\lambda$ is an absolutely irreducible graded $A$-module.
		\item $D^\lambda \cong D^\mu\langle x\rangle$, for some $x\in\mathbb{Z}$, if and only if $\lambda = \mu$ and $x=0$.
		\item $\{D^\lambda\langle x\rangle | \lambda\in\mathscr{P}_0 \text{ and } x\in\mathbb{Z}\}$ is a complete set of pairwise non-isomorphic graded simple $A$-modules.
	\end{enumerate}
\end{Theorem}

\subsection{Multipartitions}\label{S: Cherednik section}

This section sets up some of the combinatorics that we will need.

Fix an integer $e$ such that $e\geq 2$ or $e=0$ and define an indexing set $I = \mathbb{Z}/e\mathbb{Z}$. Let $\Gamma_e$ be the oriented quiver with vertices $i\in I$ and arrows $i \rightarrow i+1$, for all $i\in I$. So $\Gamma_e$ is either the linear quiver of type $A_\infty$, if $e=0$, or the cyclic quiver of type $A_e^{(1)}$, if $e \neq 0$. To this quiver we associate its symmetric Cartan matrix $(c_{ij})_{i,j\in I}$, where $c_{ij} = 2\delta_{ij} - \delta_{i(j+1)} - \delta_{i(j-1)}$.

Following Kac \cite{Kac}, let $\widehat{\mathfrak{sl}}_e$, for $e>0$, and $\mathfrak{sl}_\infty$, for $e=0$, be the \emph{Kac-Moody algebra} of $\Gamma_e$ with \emph{simple roots} $\{\alpha_i|i\in I\}$ and \emph{fundamental weights} $\{\Lambda_i|i\in I\}$. Let $(\ ,\ )$ be the bilinear form determined by
\begin{equation*}
(\alpha_i,\alpha_j) = c_{ij} \quad\text{ and }\quad (\Lambda_i,\alpha_j) = \delta_{ij}, \text{ for } i,j\in I.
\end{equation*}
Let $P_+ = \bigoplus_{i\in I} \mathbb{N}_0\Lambda_i$ be the \emph{dominant weight lattice} and $Q_+ = \bigoplus_{i\in I} \mathbb{N}_0\alpha_i$ the \emph{positive root lattice}. Fix a dominant weight $\Lambda\in P_+$ and set $\ell = \sum_{i\in I} (\Lambda,\alpha_i)$.

Fix an $n\in\mathbb{N}_0$. An $\ell$-multipartition of $n$ is an $\ell$-tuple $\lambda = (\lambda^{(1)}, \ldots, \lambda^{(\ell)})$ of partitions $\lambda^{(k)}$, such that
\begin{equation*}
|\lambda| = \sum_{k=1}^\ell \sum_{j=1}^\infty \lambda_j^{(k)} = |\lambda^{(1)}| + \ldots |\lambda^{(\ell)}| = n.
\end{equation*}
We identify a multipartition with its \emph{Young diagram} of boxes with coordinates $(r,c,l)$ of row, column, and component. We visualize these diagrams in the \emph{Russian style}, where every box has diagonal length $2\ell$.

\begin{Example}
	Let $\lambda = ((2,1),(3,1))$.
	\[
	\begin{tikzpicture}
	\draw (0,0) pic{russiantableau={scale=0.5}{{\phantom{1},\phantom{1}},{\phantom{1}}}};
	\draw (4,0) pic{russiantableau={scale=0.5}{{\phantom{1},\phantom{1},\phantom{1}},{\phantom{1}}}};
	\end{tikzpicture}
	\]
\end{Example}

Let $\mathscr{N}_n^\ell = \{(r,c,l) | r,c,l\in\mathbb{N}, r+c\leq 2n, \text{ and } 1\leq l\leq \ell\}$ be a set of boxes. Fix an $e$-multicharge $\kappa = (\kappa_1,\ldots,\kappa_\ell)\in I^\ell$ such that $(\Lambda,\alpha_i) = |\{1\leq k\leq \ell | \kappa_k = i\}|$. With the multicharge we define the \emph{residue} of a box $\alpha\in\mathscr{N}_n^\ell$ to be $\residue(\alpha) = c - r + \kappa_l$. We define an equivalence relation on boxes by $(r,c,l) \sim (r',c',l')$ if $l = l'$ and $c-r = c'-r'$. We refer to its equivalence classes as \emph{diagonals}. By definition, all boxes in a diagonal have the same residue. A diagonal is an $i$-diagonal if the residue of its boxes is $i$, for $i\in I$.

A \emph{loading} is an $\ell$-tuple $\theta = (\theta_1,\ldots,\theta_\ell) \in\mathbb{Z}^\ell$ such that $\theta_i - \theta_j \notin \ell\mathbb{Z}$, for $1\leq i < j \leq \ell$. Fix a loading $\theta$ and a quantity $\varepsilon < \frac{1}{n}$. We define a function $\theta_\varepsilon:\mathscr{N}_n^\ell \longrightarrow \mathbb{R}$, by $\theta_\varepsilon(r,c,l) = \theta_l + (r-c)\ell + (r+c)\varepsilon$. For a box $\alpha\in\mathscr{N}_n^\ell$ we then assign to it a real valued coordinate $x_\alpha = \theta_\varepsilon(\alpha)\in\mathbb{R}$.

\begin{Definition}[{\cite[Definition 1.2]{Bowman:ManyCellular}}]\label{D: Total Order}
	Let $\alpha,\beta\in\mathscr{N}_n^\ell$. Then $\alpha \prec_\theta \beta$ if $x_\alpha < x_\beta$. We say that $\alpha$ is to the \emph{left} of $\beta$ or equivalently $\beta$ is to the \emph{right} of $\alpha$.
\end{Definition}

\begin{Lemma}[{\cite[Remark 7]{Bowman:ManyCellular}}]
	$(\mathscr{N}_n^\ell,\prec_\theta)$ is a totally ordered set.
\end{Lemma}

We visualize this order by arranging our Russian diagrams such that for any box $\alpha$ its top corner aligns vertically with the coordinate $x_\alpha$.

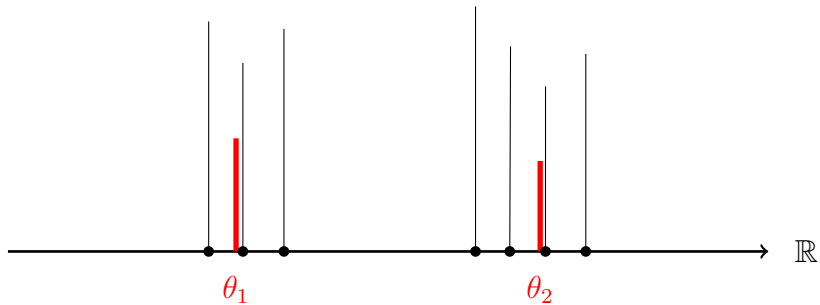
\begin{Example}
	\[
	\begin{tikzpicture}
	\begin{turn}{-5}
	\draw (0,0) pic{russiantableau={scale=0.5}{{\phantom{1},\phantom{1}},{\phantom{1}}}};
	\end{turn}
	\begin{turn}{-5}
	\draw (4,0) pic{russiantableau={scale=0.5}{{\phantom{1},\phantom{1},\phantom{1}},{\phantom{1}}}};
	\end{turn}

	\draw[line width = 1,->] (-3,-1.5) -- (7,-1.5);
	\node at (7.5,-1.5) {$\mathbb{R}$};

	\draw[red,line width = 2] (0,0) -- (0,-1.5);
	\node[red] at (0,-2) {$\theta_1$};

	\draw[red,line width = 2] (4,-0.3) -- (4,-1.5);
	\node[red] at (4,-2) {$\theta_2$};

	\draw (-0.36,1.55) -- (-0.36,-1.5);
	\draw (0.09,1) -- (0.09,-1.5);
	\draw (0.63,1.45) -- (0.63,-1.5);

	\fill (-0.36,-1.5) circle (2pt);
	\fill (0.09,-1.5) circle (2pt);
	\fill (0.63,-1.5) circle (2pt);

	\draw (3.15,1.75) -- (3.15,-1.5);
	\draw (3.61,1.22) -- (3.6,-1.5);
	\draw (4.07,0.69) -- (4.07,-1.5);
	\draw (4.6,1.12) -- (4.6,-1.5);

	\fill (3.15,-1.5) circle (2pt);
	\fill (3.6,-1.5) circle (2pt);
	\fill (4.07,-1.5) circle (2pt);
	\fill (4.6,-1.5) circle (2pt);
	\end{tikzpicture}
	\]
\end{Example}

A \emph{box configuration} is a subset of $\mathscr{N}_n^\ell$. Let $\mathscr{C}_n^\ell(\theta)$ be the set of all box configurations with $n$ boxes. The ordering of boxes above defines a partial order on $\mathscr{C}_n^\ell(\theta)$.

\begin{Definition}[{\cite[Definition 1.2]{Bowman:ManyCellular}}]\label{D: Dominance order}
	Let $\lambda,\mu\in\mathscr{C}_n^\ell(\theta)$. Then $\lambda$ dominates $\mu$ or $\mu \trianglelefteq_\theta \lambda$ if for all $\alpha\in\mathscr{N}_n^\ell$
	\begin{equation*}
	|\{\beta\in\mu | \residue(\beta) = \residue(\alpha), \beta \prec_\theta \alpha\}| \leq |\{\beta\in\lambda | \residue(\beta) = \residue(\alpha), \beta \prec_\theta \alpha\}|
	\end{equation*}
\end{Definition}

\noindent Note that by this definition, $\mu \trianglelefteq_\theta \lambda$ if and only if $\mu\setminus\lambda \trianglelefteq_\theta \lambda\setminus\mu$, for $\lambda,\mu\in\mathscr{C}_n^\ell(\theta)$.	Set $(\lambda\setminus\mu)_i = \{\alpha\in\lambda\setminus\mu | \residue(\alpha) = i\}$, for $\lambda,\mu\in\mathscr{C}_n^\ell(\theta)$ and $i\in I$.

Let $\mathscr{P}_n^\ell(\theta) \subset \mathscr{C}_n^\ell(\theta)$ be the poset of $\ell$-multipartitions of $n$ ordered by the dominance order $\trianglelefteq_\theta$. Also we set $\mathscr{P}^\ell(\theta) = \bigcup_{n\in\mathbb{N}_0} \mathscr{P}_n^\ell(\theta)$.

As Bowman \cite[Example 1.3]{Bowman:ManyCellular} remarks, $\trianglelefteq_\theta$ coincides with the normal dominance order on multipartitions, if $\theta_{i+1} - \theta_i > n\ell$, for all $i= 1,\ldots,\ell-1$. We call $\theta$ \emph{well-separated} in this case.

A box $\alpha$ is \emph{addable} for $\lambda\in\mathscr{P}_n^\ell(\theta)$ if $\lambda \cup \alpha\in \mathscr{P}_{n+1}^\ell(\theta)$ and it is \emph{removable} if $\lambda \setminus \alpha \in \mathscr{P}_{n-1}^\ell(\theta)$. Let $\lambda\in\mathscr{P}_n^\ell(\theta)$. A \emph{standard tableau} $\mathfrak{t}$ for $\lambda$ is a bijection $\mathfrak{t}: \lambda \longrightarrow \{1,\ldots,n\}$ such that
\begin{equation*}
\mathfrak{t}(r+1,c,l) > \mathfrak{t}(r,c,l) \text{ and } \mathfrak{t}(r,c+1,l) > \mathfrak{t}(r,c,l),
\end{equation*}
for all $(r,c,l),(r+1,c,l),(r,c+1,l)\in\lambda$. The \emph{shape} of a standard tableau $\mathfrak{t}$ for $\lambda$ is $\shape(\mathfrak{t}) = \lambda$. Let $\Std(\lambda)$ be the set of all standard tableaux of shape $\lambda$. Let $\mathfrak{t}_{\downarrow k}$ be the subtableau of $\mathfrak{t}$ that consists only of the boxes containing the numbers up to $k$, for $1\leq k\leq n$. For a standard tableau $\mathfrak{t}$, the box containing the number $k$ is $\mathfrak{t}^{-1}(k)$. Then the \emph{residue sequence} of $\mathfrak{t}$, $\residue(\mathfrak{t})$ is the $n$-tuple $(\residue(\mathfrak{t}^{-1}(1)),\ldots, \residue(\mathfrak{t}^{-1}(n)))\in I^n$. More general the \emph{box configuration residue sequence} $\residue(\lambda)$, for $\lambda\in\mathscr{C}_n^\ell(\theta)$, is the sequence of residues of the boxes of $\lambda$ ordered from left to right.

Additionally, for a tableau $\mathfrak{t}$ of shape $\lambda\in\mathscr{P}_n^\ell(\theta)$ and for some $1\leq k\leq n$ let $\mathscr{A}_\mathfrak{t}^\theta(k)$ and $\mathscr{R}_\mathfrak{t}^\theta(k)$ be the sets of addable and respectively removable boxes of the multipartition $\shape(\mathfrak{t}_{\downarrow k})$ of residue $\residue(\mathfrak{t}^{-1}(k))$ and are to the right of $\mathfrak{t}^{-1}(k)$ with respect to $\prec_\theta$.

\begin{Definition}[{\cite[Definition 2.15]{Webster:RouquierConjecture}},{\cite[Definition 1.7]{Bowman:ManyCellular}}]\label{D: Degree function}
	For $\lambda\in\mathscr{P}_n^\ell(\theta)$ and $\mathfrak{t}\in\Std(\lambda)$, define the degree of $\mathfrak{t}$ to be
	$\deg_\theta(\mathfrak{t}) = \sum_{k=1}^n \left(\left|\mathscr{A}_\mathfrak{t}^\theta(k)\right| - \left|\mathscr{R}_\mathfrak{t}^\theta(k)\right|\right)$.
\end{Definition}

\subsection{Diagrammatic Cherednik Algebra}\label{S: Cherednik algebra section}

In this section, following Webster \cite{Webster:RouquierConjecture}, we introduce the diagrammatic Cherednik algebra, which is a graded cellular algebra.

To define the Diagrammatic Cherednik algebra, we need to define the diagrams that span this diagrammatic algebra. Let $x,y\in\mathbb{R}$ then a \emph{string of type} $(x,y)$ is a coloured path in the Euclidean plane from $(x,1)$ to $(y,0)$ such that it is diffeomorphic to the unit interval from $(0,1)$ to $(0,0)$. We assign colours to strings to help distinguish them. A \emph{solid string} of type $(x,y)$ is a black string of type $(x,y)$. A \emph{ghost string} of a solid string is a dotted grey string and a translation of the solid string by $\ell$ units to the left. A \emph{red string} is a red string of type $(\theta_l,\theta_l)$, for a $1\leq l\leq \ell$.

A solid string and its ghost string carry a residue $i\in I$, we call such strings $i$-\emph{strings}. Red strings of type $(\theta_l,\theta_l)$ carry the residue $\kappa_l$. Solid strings are allowed to have a finite number of dots on them.

Given $\lambda\in\mathscr{C}_n^\ell(\theta)$, let $X_\lambda = \{x_\alpha | \alpha\in\lambda\} \subset \mathbb{R}$ be the discrete subset of points on the real line corresponding to all boxes $\alpha\in\lambda$. Let $\hom(\lambda,\mu)$ be the set of all bijections from $X_\lambda$ into $X_\mu$.

Fix $\lambda,\mu\in\mathscr{C}_n^\ell(\theta)$. A $\theta$-\emph{diagram of type} $(\lambda,\mu)$ is a diagram that
\begin{itemize}
	\item has $\ell$ red strings of type $(\theta_l,\theta_l)$, for $1\leq l\leq \ell$,
	\item has $n$ solid strings of type $(x_\alpha,f(x_\alpha))$, for $x_\alpha\in X_\lambda$ and $f\in\hom(\lambda,\mu)$,
	\item has $n$ ghost strings,
	\item has any finite number of dots on its solid strings,
	\item has a residue $i\in I$ assigned to each solid string,
	\item has at most two strings intersect at any point,
	\item has no dots on any intersection of strings, and
	\item has no tangencies of strings.
\end{itemize}
Two $\theta$-diagrams are equivalent if they are the same diagram up to isotopy.

Let $D$ be a diagram of type $(\lambda,\mu)$. Define the pairs $\toppo(D) = (\lambda,\mathbf{i})$ and $\bottom(D) = (\mu,\mathbf{j})$, where $\mathbf{i},\mathbf{j}\in I^n$ are the sequences of residues read from left to right at the top and bottom of the diagram $D$ respectively.

\begin{Example}
	Let $\theta=(0,8)$, $\lambda = ((2,1),(3,1))$, and $\mu =(\emptyset,(1^7))$. Then a possible $\theta$-diagram of type $(\lambda,\mu)$ is:
	
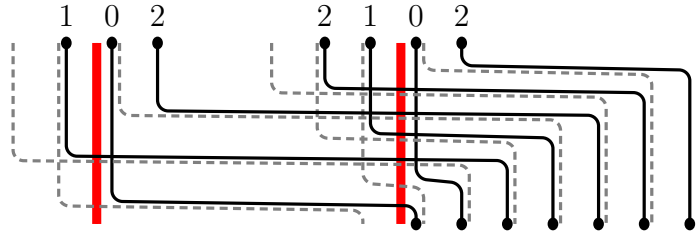
\[
	\begin{tikzpicture}[very thick,baseline,rounded corners,yscale = 1.2]
	\draw[wei](0,0) -- (0,2);
	\draw[wei](4,0) -- (4,2);

	\fill (-0.4,2) circle (2pt);
	\fill (0.2,2) circle (2pt);
	\fill (0.8,2) circle (2pt);

	\fill (4.2,2) circle (2pt);
	\fill (3.6,2) circle (2pt);
	\fill (3,2) circle (2pt);
	\fill (4.8,2) circle (2pt);

	\node at (-0.4,2.3) {$1$};
	\node at (0.2,2.3) {$0$};
	\node at (0.8,2.3) {$2$};

	\node at (4.2,2.3) {$0$};
	\node at (3.6,2.3) {$1$};
	\node at (3,2.3) {$2$};
	\node at (4.8,2.3) {$2$};

	\fill (4.2,0) circle (2pt);
	\fill (4.8,0) circle (2pt);
	\fill (5.4,0) circle (2pt);
	\fill (6,0) circle (2pt);
	\fill (6.6,0) circle (2pt);
	\fill (7.2,0) circle (2pt);
	\fill (7.8,0) circle (2pt);

	\draw[ghost] (-0.5,2) +(0,0) -- +(0,-1.8) -- +(4,-1.85) -- +(4,-2);
	\draw[solid] (0.2,2) +(0,0) -- +(0,-1.75) -- +(4,-1.8) -- +(4,-2);

	\draw[ghost] (3.5,2) +(0,0) -- +(0,-1.55) -- +(0.8,-1.6) -- +(0.8,-2);
	\draw[solid] (4.2,2) +(0,0) -- +(0,-1.5) -- +(0.6,-1.55) -- +(0.6,-2);

	\draw[ghost] (-1.1,2) +(0,0) -- +(0,-1.3) -- +(6,-1.35) -- +(6,-2);
	\draw[solid] (-0.4,2) +(0,0) -- +(0,-1.25) -- +(5.8,-1.3) -- +(5.8,-2);

	\draw[ghost] (2.9,2) +(0,0) -- +(0,-1.05) -- +(2.6,-1.1) -- +(2.6,-2);
	\draw[solid] (3.6,2) +(0,0) -- +(0,-1) -- +(2.4,-1.05) -- +(2.4,-2);

	\draw[ghost] (0.3,2) +(0,0) -- +(0,-0.8) -- +(5.8,-0.85) -- +(5.8,-2);
	\draw[solid] (0.8,2) +(0,0) -- +(0,-0.75) -- +(5.8,-0.8) -- +(5.8,-2);

	\draw[ghost] (2.3,2) +(0,0) -- +(0,-0.55) -- +(4.4,-0.6) -- +(4.4,-2);
	\draw[solid] (3,2) +(0,0) -- +(0,-0.5) -- +(4.2,-0.55) -- +(4.2,-2);

	\draw[ghost] (4.3,2) +(0,0) -- +(0,-0.3) -- +(3,-0.35) -- +(3,-2);
	\draw[solid] (4.8,2) +(0,0) -- +(0,-0.25) -- +(3,-0.3) -- +(3,-2);
	\end{tikzpicture}
	\]
	Here $\toppo(D) = (\lambda,(1,0,2,2,1,0,2))$ and $\bottom(D) = (\mu,(0,0,1,1,2,2,2))$.
\end{Example}

\emph{Note} that throughout this paper we draw our diagrams with solid strings that have segments in which they are drawn almost horizontally. It will often look like the solid and its ghost string overlap in these segments. However, this is just due to difficulties in plotting these diagrams. In fact these strings are not overlapping but are just very close to each other in these segments.

Let $D,E$ be two $\theta$-diagrams of types $(\lambda,\mu)$ and $(\nu,\xi)$. Then we construct a diagram $D\times E$ by translating $D$ upwards by $1$ unit on top of $E$ and rescaling the stacked diagram vertically by a factor of $\frac{1}{2}$. We then define a multiplication on $\theta$-diagrams by
\begin{equation*}
D\circ E = \begin{cases}
D\times E, &\text{if } \bottom(D) = \toppo(E)\\
0, &\text{otherwise}.
\end{cases}
\end{equation*}
If $D\circ E \neq 0$ then it is a $\theta$-diagram of type $(\lambda,\xi)$.

There is an anti-isomorphism $\ast$ on $\theta$-diagrams that is given by reflection through the line $y=\frac{1}{2}$, which turns a $\theta$-diagram of type $(\lambda,\mu)$ into a $\theta$-diagram of type $(\mu,\lambda)$.

A $\theta$-diagram is \emph{unsteady} if there is a solid string $n\ell$ units to the left of the leftmost red string.

\begin{Definition}[{\cite[Definition 4.2]{Webster:RouquierConjecture}}]\label{D: Cherednik Algebra}
	The diagrammatic Cherednik algebra $\mathbb{A}_n^\kappa(\theta)$ is the unital associative algebra spanned by all $\theta$-diagrams with $n$ solid strings and with multiplication $\circ$ and subject to the following relations:
	\begin{enumerate}
		\item(Dots and crossings)\label{D: Cherednik Dots and crossings} A dot on an $i$-string can pass through any crossing except:
		\[\begin{tikzpicture}[very thick,baseline,rounded corners,scale = 0.5]
		\draw[solid] (0,0) +(0,0) -- +(2,2);
		\draw[solid] (2,0) +(0,0) -- +(-2,2);
		\fill (0.5,1.5) circle (5pt);

		\node at (3,1) {$-$};

		\draw[solid] (4,0) +(0,0) -- +(2,2);
		\draw[solid] (6,0) +(0,0) -- +(-2,2);
		\fill (5.5,0.5) circle (5pt);

		\node at (7,1) {$=$};

		\draw[solid] (8,0) +(0,0) -- +(0,2);
		\draw[solid] (10,0) +(0,0) -- +(0,2);

		\node at (11,1) {$=$};

		\draw[solid] (12,0) +(0,0) -- +(2,2);
		\draw[solid] (14,0) +(0,0) -- +(-2,2);
		\fill (12.5,0.5) circle (5pt);

		\node at (15,1) {$-$};

		\draw[solid] (16,0) +(0,0) -- +(2,2);
		\draw[solid] (18,0) +(0,0) -- +(-2,2);
		\fill (17.5,1.5) circle (5pt);

		\node at (0,-0.5) {$i$};
		\node at (2,-0.5) {$i$};
		\node at (4,-0.5) {$i$};
		\node at (6,-0.5) {$i$};
		\node at (8,-0.5) {$i$};
		\node at (10,-0.5) {$i$};
		\node at (12,-0.5) {$i$};
		\node at (14,-0.5) {$i$};
		\node at (16,-0.5) {$i$};
		\node at (18,-0.5) {$i$};
		\end{tikzpicture}\]
		\item(Double crossings)\label{D: Cherednik double crossings} A double crossing between any two strings can be
		pulled apart except in the following cases, where $j=i+1\in I$:
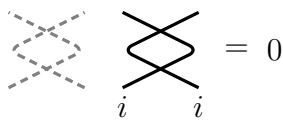
		\begin{enumerate}[label=(\Roman*)]
			\item\label{D: Cherednik double i crossing} \[\begin{tikzpicture}[very thick,baseline,rounded corners,scale = 0.5]
			\draw[solid] (0,0) +(0,0) -- +(2,1) -- +(0,2);
			\draw[solid] (2,0) +(0,0) -- +(-2,1) -- +(0,2);

			\draw[ghost] (-3,0) +(0,0) -- +(2,1) -- +(0,2);
			\draw[ghost] (-1,0) +(0,0) -- +(-2,1) -- +(0,2);

			\node at (3,1) {$=$};

			\node at (4,1) {$0$};

			\node at (0,-0.5) {$i$};
			\node at (2,-0.5) {$i$};
			\end{tikzpicture}\]
			\item
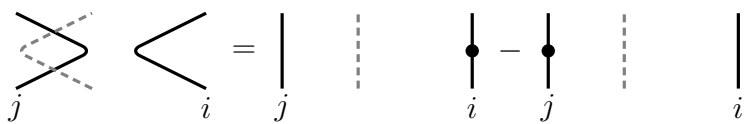
\label{D: Cherednik solid ghost double} \[\begin{tikzpicture}[very thick,baseline,rounded corners,scale = 0.5]
			\draw[solid] (-3,0) +(0,0) -- +(2,1) -- +(0,2);
			\draw[ghost] (-1,0) +(0,0) -- +(-2,1) -- +(0,2);
			\draw[solid] (2,0) +(0,0) -- +(-2,1) -- +(0,2);

			\node at (3,1) {$=$};

			\draw[solid] (4,0) +(0,0) -- +(0,2);
			\draw[ghost] (6,0) +(0,0) -- +(0,2);
			\draw[solid] (9,0) +(0,0) -- +(0,2);
			\fill (9,1) circle (5pt);

			\node at (10,1) {$-$};

			\draw[solid] (11,0) +(0,0) -- +(0,2);
			\draw[ghost] (13,0) +(0,0) -- +(0,2);
			\draw[solid] (16,0) +(0,0) -- +(0,2);
			\fill (11,1) circle (5pt);

			\node at (-3,-0.5) {$j$};
			\node at (2,-0.5) {$i$};
			\node at (4,-0.5) {$j$};
			\node at (9,-0.5) {$i$};
			\node at (11,-0.5) {$j$};
			\node at (16,-0.5) {$i$};
			\end{tikzpicture}\]
			\item
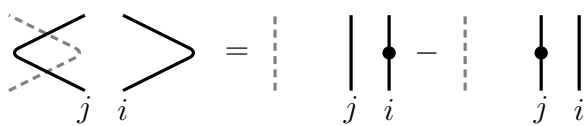
\label{D: Cherednik ghost solid double} \[\begin{tikzpicture}[very thick,baseline,rounded corners,scale = 0.5]
			\draw[ghost] (-3,0) +(0,0) -- +(2,1) -- +(0,2);
			\draw[solid] (-1,0) +(0,0) -- +(-2,1) -- +(0,2);
			\draw[solid] (0,0) +(0,0) -- +(2,1) -- +(0,2);

			\node at (3,1) {$=$};

			\draw[ghost] (4,0) +(0,0) -- +(0,2);
			\draw[solid] (6,0) +(0,0) -- +(0,2);
			\draw[solid] (7,0) +(0,0) -- +(0,2);
			\fill (7,1) circle (5pt);

			\node at (8,1) {$-$};

			\draw[ghost] (9,0) +(0,0) -- +(0,2);
			\draw[solid] (11,0) +(0,0) -- +(0,2);
			\draw[solid] (12,0) +(0,0) -- +(0,2);
			\fill (11,1) circle (5pt);

			\node at (-1,-0.5) {$j$};
			\node at (0,-0.5) {$i$};
			\node at (6,-0.5) {$j$};
			\node at (7,-0.5) {$i$};
			\node at (11,-0.5) {$j$};
			\node at (12,-0.5) {$i$};
			\end{tikzpicture}\]
			\item
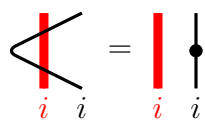
\label{D: Cherednik red solid double} \[\begin{tikzpicture}[very thick,baseline,rounded corners,scale = 0.5]
			\draw[wei] (0,0) +(0,0) -- +(0,2);
			\draw[solid] (1,0) +(0,0) -- +(-2,1) -- +(0,2);

			\node at (2,1) {$=$};

			\draw[wei] (3,0) +(0,0) -- +(0,2);
			\draw[solid] (4,0) +(0,0) -- +(0,2);
			\fill (4,1) circle (5pt);

			\node[red] at (0,-0.5) {$i$};
			\node at (1,-0.5) {$i$};
			\node[red] at (3,-0.5) {$i$};
			\node at (4,-0.5) {$i$};
			\end{tikzpicture}\]
			\item
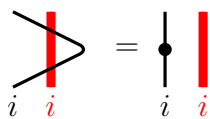
\label{D: Cherednik solid red double} \[\begin{tikzpicture}[very thick,baseline,rounded corners,scale = 0.5]
			\draw[wei] (0,0) +(0,0) -- +(0,2);
			\draw[solid] (-1,0) +(0,0) -- +(2,1) -- +(0,2);

			\node at (2,1) {$=$};

			\draw[wei] (4,0) +(0,0) -- +(0,2);
			\draw[solid] (3,0) +(0,0) -- +(0,2);
			\fill (3,1) circle (5pt);

			\node[red] at (0,-0.5) {$i$};
			\node at (-1,-0.5) {$i$};
			\node[red] at (4,-0.5) {$i$};
			\node at (3,-0.5) {$i$};
			\end{tikzpicture}\]
		\end{enumerate}
		\item(Triple crossings)\label{D: Cherednik triple crossings} A string can be pulled through a triple crossing except
		in the following cases, where $j=i+1\in I$:
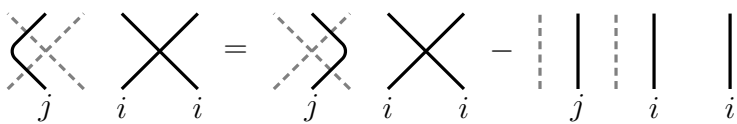
		\begin{enumerate}[label=(\Roman*)]
			\item\label{D: Cherednik 2-ghost triple} \[\begin{tikzpicture}[very thick,baseline,rounded corners,scale = 0.5]
			\draw[ghost] (0,0) +(0,0) -- +(2,2);
			\draw[ghost] (2,0) +(0,0) -- +(-2,2);
			\draw[solid] (1,0) +(0,0) -- +(-1,1) -- +(0,2);

			\draw[solid] (3,0) +(0,0) -- +(2,2);
			\draw[solid] (5,0) +(0,0) -- +(-2,2);

			\node at (6,1) {$=$};

			\draw[ghost] (7,0) +(0,0) -- +(2,2);
			\draw[ghost] (9,0) +(0,0) -- +(-2,2);
			\draw[solid] (8,0) +(0,0) -- +(1,1) -- +(0,2);

			\draw[solid] (10,0) +(0,0) -- +(2,2);
			\draw[solid] (12,0) +(0,0) -- +(-2,2);

			\node at (13,1) {$-$};

			\draw[ghost] (14,0) +(0,0) -- +(0,2);
			\draw[ghost] (16,0) +(0,0) -- +(0,2);
			\draw[solid] (15,0) +(0,0) -- +(0,2);

			\draw[solid] (17,0) +(0,0) -- +(0,2);
			\draw[solid] (19,0) +(0,0) -- +(0,2);

			\node at (1,-0.5) {$j$};
			\node at (3,-0.5) {$i$};
			\node at (5,-0.5) {$i$};

			\node at (8,-0.5) {$j$};
			\node at (10,-0.5) {$i$};
			\node at (12,-0.5) {$i$};

			\node at (15,-0.5) {$j$};
			\node at (17,-0.5) {$i$};
			\node at (19,-0.5) {$i$};
			\end{tikzpicture}\]
			\item
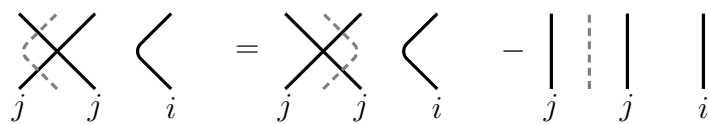
\label{D: Cherednik 2-solid triple} \[\begin{tikzpicture}[very thick,baseline,rounded corners,scale = 0.5]
			\draw[solid] (0,0) +(0,0) -- +(2,2);
			\draw[solid] (2,0) +(0,0) -- +(-2,2);
			\draw[ghost] (1,0) +(0,0) -- +(-1,1) -- +(0,2);

			\draw[solid] (4,0) +(0,0) -- +(-1,1) -- +(0,2);

			\node at (6,1) {$=$};

			\draw[solid] (7,0) +(0,0) -- +(2,2);
			\draw[solid] (9,0) +(0,0) -- +(-2,2);
			\draw[ghost] (8,0) +(0,0) -- +(1,1) -- +(0,2);

			\draw[solid] (11,0) +(0,0) -- +(-1,1) -- +(0,2);

			\node at (13,1) {$-$};

			\draw[solid] (14,0) +(0,0) -- +(0,2);
			\draw[solid] (16,0) +(0,0) -- +(0,2);
			\draw[ghost] (15,0) +(0,0) -- +(0,2);

			\draw[solid] (18,0) +(0,0) -- +(0,2);

			\node at (0,-0.5) {$j$};
			\node at (2,-0.5) {$j$};
			\node at (4,-0.5) {$i$};

			\node at (7,-0.5) {$j$};
			\node at (9,-0.5) {$j$};
			\node at (11,-0.5) {$i$};

			\node at (14,-0.5) {$j$};
			\node at (16,-0.5) {$j$};
			\node at (18,-0.5) {$i$};
			\end{tikzpicture}\]
			\item
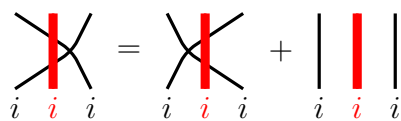
\label{D: Cherednik red triple} \[\begin{tikzpicture}[very thick,baseline,rounded corners,scale = 0.5]
			\draw[solid] (0,0) +(0,0) -- +(1.5,1) -- +(2,2);
			\draw[solid] (2,0) +(0,0) -- +(-0.5,1) -- +(-2,2);
			\draw[wei] (1,0) +(0,0) -- +(0,2);

			\node at (3,1) {$=$};

			\draw[solid] (4,0) +(0,0) -- +(0.5,1) -- +(2,2);
			\draw[solid] (6,0) +(0,0) -- +(-1.5,1) -- +(-2,2);
			\draw[wei] (5,0) +(0,0) -- +(0,2);

			\node at (7,1) {$+$};

			\draw[solid] (8,0) +(0,0) -- +(0,2);
			\draw[solid] (10,0) +(0,0) -- +(0,2);
			\draw[wei] (9,0) +(0,0) -- +(0,2);

			\node at (0,-0.5) {$i$};
			\node[red] at (1,-0.5) {$i$};
			\node at (2,-0.5) {$i$};

			\node at (4,-0.5) {$i$};
			\node[red] at (5,-0.5) {$i$};
			\node at (6,-0.5) {$i$};

			\node at (8,-0.5) {$i$};
			\node[red] at (9,-0.5) {$i$};
			\node at (10,-0.5) {$i$};
			\end{tikzpicture}\]
		\end{enumerate}
		\item(Unsteady)\label{D: Cherednik unsteady} Any unsteady diagram is zero.
	\end{enumerate}
\end{Definition}

The Cherednik algebra comes with a natural grading, which is defined locally.
\begin{align*}
\deg(\ \begin{tikzpicture}[very thick,baseline,rounded corners,scale = 0.25]
\draw[solid] (0,0) +(0,0) -- +(0,2);
\fill (0,1) circle (10pt);
\node at (0,-1) {i};
\end{tikzpicture}\ ) &= 2 \\
\deg(\ \begin{tikzpicture}[very thick,baseline,rounded corners,scale = 0.25]
\draw[solid] (0,0) +(0,0) -- +(2,2);
\draw[solid] (2,0) +(0,0) -- +(-2,2);
\node at (0,-1) {i};
\node at (2,-1) {j};
\end{tikzpicture}\ ) &= -2\delta_{ij}\\
\deg(\ \begin{tikzpicture}[very thick,baseline,rounded corners,scale = 0.25]
\draw[ghost] (0,0) +(0,0) -- +(2,2);
\draw[solid] (2,0) +(0,0) -- +(-2,2);
\node at (0,-1) {i};
\node at (2,-1) {j};
\end{tikzpicture}\ ) &= \delta_{i+1,j} =
\deg(\ \begin{tikzpicture}[very thick,baseline,rounded corners,scale = 0.25]
\draw[solid] (0,0) +(0,0) -- +(2,2);
\draw[ghost] (2,0) +(0,0) -- +(-2,2);
\node at (0,-1) {j};
\node at (2,-1) {i};
\end{tikzpicture}\ ) \\
\deg(\ \begin{tikzpicture}[very thick,baseline,rounded corners,scale = 0.25]
\draw[wei] (0,0) +(0,0) -- +(0,2);
\draw[solid] (1,0) +(0,0) -- +(-2,2);
\node[red] at (0,-1) {i};
\node at (1,-1) {j};
\end{tikzpicture}\ ) &= \delta_{ij} =
\deg(\ \begin{tikzpicture}[very thick,baseline,rounded corners,scale = 0.25]
\draw[wei] (0,0) +(0,0) -- +(0,2);
\draw[solid] (-1,0) +(0,0) -- +(2,2);
\node[red] at (0,-1) {i};
\node at (-1,-1) {j};
\end{tikzpicture}\ )
\end{align*}

\begin{Definition}[{\cite[Definition 2.1]{Bowman:ManyCellular}}]
	Given $\lambda\in\mathscr{P}_n^\ell(\theta)$ and $\mu\in\mathscr{C}_n^\ell(\theta)$ a \emph{semistandard tableau} $\mathfrak{T}$ of shape $\lambda$ and weight $\mu$ is a bijection $\mathfrak{T}:\lambda \longrightarrow X_\mu$ such that
	\begin{enumerate}
		\item $\mathfrak{T}(1,1,l) > \theta_l$,
		\item $\mathfrak{T}(r,c,l) > \mathfrak{T}(r-1,c,l) + \ell$,
		\item $\mathfrak{T}(r,c,l) > \mathfrak{T}(r,c-1,l) - \ell$.
	\end{enumerate}
	The set of all semistandard tableau of shape $\lambda$ and weight $\mu$ is denoted by $\SStd_\theta(\lambda,\mu)$
\end{Definition}

Given a semistandard tableau $\mathfrak{T}$ of shape $\lambda$ and weight $\mu$ we can associate to it a $\theta$-diagram $C_\mathfrak{T}$ of type $(\lambda,\mu)$, where the solid strings are strings of type $(x_\alpha,\mathfrak{T}(\alpha))$ and residue $\residue(\alpha)$, for all $\alpha\in\lambda$, and such that there is a minimal number of intersections.

This choice of diagram is not unique, so we choose one for every semistandard tableau and fix this choice until further notice. In fact, later we will fix special diagrams for some tableaux. By Bowman \cite[4.1]{Bowman:ManyCellular}, the results in this section are independent of the choice of diagrams. Taking a pair of  semistandard tableaux $(\mathfrak{S},\mathfrak{T})\in\SStd_\theta(\lambda,\mu)\times\SStd_\theta(\lambda,\nu)$, we define the $\theta$-diagram $C_{\mathfrak{ST}} = C_\mathfrak{S}^\ast \circ C_\mathfrak{T}$, where $\ast$ is the anti-isomorphism on $\mathbb{A}_n^\kappa(\theta)$ that flips the diagram about the line $y=\frac{1}{2}$.

\begin{Theorem}[{\cite[Theorem 4.10]{Bowman:ManyCellular}}]
	The algebra $\mathbb{A}_n^\kappa(\theta)$ is a graded cellular algebra with poset $(\mathscr{P}_n^\ell(\theta),\vartriangleleft_\theta)$ and graded cellular basis
	\begin{equation*}
	\{C_{\mathfrak{ST}}^\lambda | (\mathfrak{S},\mathfrak{T})\in\SStd\phantom{}_\theta(\lambda,\mu)\times\SStd\phantom{}_\theta(\lambda,\nu), \lambda\in\mathscr{P}_n^\ell(\theta), \text{ and } \mu,\nu\in\mathscr{C}_n^\ell(\theta)\},
	\end{equation*}
	with $\deg_\theta(C_{\mathfrak{ST}}^\lambda) = \deg_\theta(\mathfrak{S}) + \deg_\theta(\mathfrak{T})$.
\end{Theorem}

Fix the weight $\omega := (\emptyset,\ldots,\emptyset,(1^n)) \in\mathscr{P}_n^\ell(\theta)$. Then there is a correspondence $\Psi:\Std_\theta(\lambda) \rightarrow \SStd\phantom{}_\theta(\lambda,\omega)$ \cite[Proposition 4.4]{Bowman:ManyCellular} between standard tableaux of shape $\lambda$ and semistandard tableaux of shape $\lambda$ and weight $\omega$, for any $\lambda\in\mathscr{P}_n^\ell(\theta)$.

For any $\mu\in\mathscr{C}_n^\ell(\theta)$ let $e_\mu(\mathbf{i})$ be the $\theta$-diagram of type $(\mu,\mu)$ such that all strings are straight vertical and the $k$-th solid string from the left is an $i_k$-string, for $1\leq k\leq n$. We call $e_\mu(\mathbf{i})$ the \emph{identity diagram of type} $\mu$ \emph{and residue} $\mathbf{i}$. Define the idempotent $E_\omega = \sum_{\mathbf{i}\in I^n} e_\omega(\mathbf{i})$ and consider the graded cellular subalgebra $E_\omega \mathbb{A}_n^\kappa(\theta)E_\omega$. Define $c_\mathfrak{st}^\lambda = C_{\mathfrak{ST}}^\lambda$, where $\mathfrak{S} = \Psi(\mathfrak{s})$ and $\mathfrak{T} = \Psi(\mathfrak{t})$, for $\lambda\in\mathscr{P}_n^\ell(\theta)$ and $\mathfrak{s},\mathfrak{t}\in\Std_\theta(\lambda)$. Then $E_\omega \mathbb{A}_n^\kappa(\theta)E_\omega$ has graded cellular basis
\begin{equation*}
\{C_{\mathfrak{ST}}^\lambda | \lambda\in\mathscr{P}_n^\ell(\theta),\ \mathfrak{S}, \mathfrak{T}\in\SStd\phantom{}_\theta(\lambda,\omega)\} = \{c_\mathfrak{st}^\lambda | \lambda\in\mathscr{P}_n^\ell(\theta),\ \mathfrak{s},\mathfrak{t}\in \Std(\lambda)\}.
\end{equation*}

Let $\mathscr{R}_n^\Lambda$ be the \emph{cyclotomic quiver Hecke algebra}, defined by Brundan and Kleshchev \cite[2.3]{BK:GradedDecomp}. As we do not need the actual definition of this algebra for our computations to avoid unnecessary technicalities we refrain from giving it here.

\begin{Theorem}[{\cite[Proposition 2.19]{Webster:RouquierConjecture}}, {\cite[Theorem 4.5]{Bowman:ManyCellular}}]\label{T: KLR isom}
	Let $\theta\in\mathbb{Z}^\ell$ be any loading, then $\mathscr{R}_n^\Lambda \cong E_\omega \mathbb{A}_n^\kappa(\theta)E_\omega$ as graded cellular algebras. In particular, $E_\omega \mathbb{A}_n^\kappa(\theta)E_\omega$ is independent of $\theta$.
\end{Theorem}

\begin{Corollary}[{\cite[Theorem 2.26]{Webster:RouquierConjecture}}, {\cite[Theorem 5.1]{Bowman:ManyCellular}}]\label{C: R is Cellular Alg}
	The cyclotomic KLR algebra $\mathscr{R}_n^\Lambda$ is a graded cellular algebra with graded cellular basis $\{c_\mathfrak{st}^\lambda\}$, where $\deg_\theta(c_\mathfrak{st}^\lambda) = \deg_\theta(\mathfrak{s}) + \deg_\theta(\mathfrak{t})$.
\end{Corollary}

Henceforth, we set $\mathscr{R}_n^\Lambda(\theta) = E_\omega \mathbb{A}_n^\kappa(\theta)E_\omega$. Define $S_\theta^\lambda$ to be the graded cell module of $\mathscr{R}_n^\Lambda(\theta)$ determined by $\lambda$ and the graded cellular basis $\{c_\mathfrak{st}^\lambda\}$ and define $D_\theta^\lambda = S_\theta^\lambda / \radic(S_\theta^\lambda)$. Note that even though $\mathscr{R}_n^\Lambda(\theta)$ does not depend on $\theta$, the basis $\{c_\mathfrak{st}^\lambda\}$ and, hence, $S_\theta^\lambda$ and $D_\theta^\lambda$ all do. For later reference, we want to relate these results to the (ungraded) Ariki-Koike Hecke algebras.

\begin{Definition}[see {\cite[Definition 3.1]{AK}}]\label{D: ArikiKoike}
	Let $1\neq q\in K^x$ such that $q^e = 1$. The \emph{cyclotomic Ariki-Koike Hecke algebra} $\mathcal{H}_n^\Lambda(q)$ is the unital associative $K$-algebra with generators $\{T_0,\ldots,T_{n-1}\}$ and relations:
	\begin{align*}
	(T_0 - q^{\kappa_1}) \ldots (T_0 - q^{\kappa_\ell}) = 0,&\\
	T_0 T_1 T_0 T_1 = T_1 T_0 T_1 T_0,&\\
	(T_i + 1)(T_i - q) = 0,& \text{ for } 1\leq i\leq n-1,\\
	T_{i+1} T_i T_{i+1} = T_i T_{i+1} T_i,& \text{ for } 1\leq i\leq n-2,\\
	T_i T_j = T_j T_i,& \text{ for } 0\leq i < j-1 \leq n-2.
	\end{align*}
\end{Definition}

\begin{Theorem}[{\cite[Main Theorem]{BK:GradedKL}}]\label{T: BK isom}
	There is a graded isomorphism $\mathcal{H}_n^\Lambda(q) \cong \mathscr{R}_n^\Lambda(\theta)$.
\end{Theorem}

\section{Decomposition paths}\label{S: Decomposition paths}

Recently, Jacon \cite{Jacon:Staggered} introduced staggered sequences in order to better understand the set of Uglov multipartitions. These multipartitions are important because they index the vertices in crystal graphs of the integral highest weight modules for $U_q(\widehat{\mathfrak{sl}}_e)$. We extend Jacon's ideas to define decomposition paths. We use these paths to compute the inner products in the graded cell modules in \autoref{S: Inner products}.

\subsection{Uglov multipartitions and decomposition paths}

In this section we define the set of Uglov multipartitions that is originally defined recursively in the literature using good boxes; see for example \cite{MisraMiwa}, \cite[Introduction]{AM:simples}, and \cite[Section 2.2]{Uglov}. We recall Jacon's \cite{Jacon:Staggered} recent result that gives an alternative recursive description of this set using staggered sequences. We use these sequences to define decomposition paths that help us further understand the elements in this set.

Let $\lambda\in\mathscr{P}_n^\ell(\theta)$. Define $\mathcal{A}_\theta(\lambda)$ the set of addable boxes of $\lambda$ and $\mathcal{R}_\theta(\lambda)$ the set of removable boxes of $\lambda$ and let $\mathcal{B}_\theta(\lambda) = \mathcal{A}_\theta(\lambda) \cup \mathcal{R}_\theta(\lambda)$. For $i\in I$, $\mathcal{A}_\theta^i(\lambda) = \{\alpha\in\mathcal{A}_\theta(\lambda) | \residue(\alpha) = i\}$ and $\mathcal{R}_\theta^i(\lambda) = \{\rho\in\mathcal{R}_\theta(\lambda) | \residue(\rho) = i\}$. Similarly, $\mathcal{B}_\theta^i(\lambda) = \mathcal{A}_\theta^i(\lambda) \cup \mathcal{R}_\theta^i(\lambda)$. $\alpha\in\mathcal{A}_\theta^i$ is an \emph{addable normal box} if for all $\rho\in\mathcal{R}_\theta^i(\lambda)$ with $\rho \prec_\theta \alpha$
\begin{equation*}
|\{\tilde{\alpha}\in\mathcal{A}_\theta^i(\lambda) | \rho \prec_\theta \tilde{\alpha} \prec_\theta \alpha\}| \geq |\{\tilde{\rho}\in\mathcal{R}_\theta^i(\lambda) | \rho \preceq_\theta \tilde{\rho} \prec_\theta \alpha\}|.
\end{equation*}
$\rho\in\mathcal{R}_\theta^i$ is a \emph{removable normal box} if for all $\alpha\in\mathcal{A}_\theta^i(\lambda)$ with $\rho \prec_\theta \alpha$
\begin{equation*}
|\{\tilde{\rho}\in\mathcal{R}_\theta^i(\lambda) | \rho \prec_\theta \tilde{\rho} \prec_\theta \alpha\}| \geq |\{\tilde{\alpha}\in\mathcal{A}_\theta^i(\lambda) | \rho \prec_\theta \tilde{\alpha} \preceq_\theta \alpha\}|.
\end{equation*}
Let
\begin{align*}
\mathcal{NA}_\theta^i(\lambda) &= \{\alpha\in\mathcal{A}_\theta^i(\lambda) | \alpha \text{ is normal}\} \quad\text{ and}\\
\mathcal{NR}_\theta^i(\lambda) &= \{\rho\in\mathcal{R}_\theta^i(\lambda) | \rho \text{ is normal}\}.
\end{align*}
$\alpha\in\mathcal{NA}_\theta^i(\lambda)$ is an \emph{addable good box} if for all $\alpha'\in\mathcal{NA}_\theta^i(\lambda)$: $\alpha' \preceq_\theta \alpha$. $\rho\in\mathcal{NR}_\theta^i(\lambda)$ is a \emph{removable good box} if for all $\rho'\in\mathcal{NR}_\theta^i(\lambda)$: $\rho \preceq_\theta \rho'$. 

\begin{Definition}
	The set of $\theta$-\emph{Uglov multipartitions} $\mathscr{U}^\ell(\theta)$ is the subset of $\mathscr{P}^\ell(\theta)$ determined by $\lambda\in\mathscr{U}^\ell(\theta)$ if either $\lambda = \emptyset$ or $\lambda$ has a removable good box $\rho$ and $\lambda\setminus\rho\in\mathscr{U}^\ell(\theta)$. Define $\mathscr{U}_n^\ell(\theta) = \mathscr{U}^\ell(\theta) \cap \mathscr{P}_n^\ell(\theta)$ to be the set of \emph{Uglov $\ell$-multipartitions of $n$}. If $\theta$ is well-separated then $\mathscr{K}_n^\ell = \mathscr{U}_n^\ell(\theta)$ is the set of \emph{Kleshchev multipartitions}.
\end{Definition}

\noindent Jacon \cite{Jacon:Staggered} gave an alternative recursive description of the set $\mathscr{U}^\ell(\theta)$.

\begin{Definition}
	Let $i\in I$ and $\alpha,\beta\in\mathcal{B}_\theta^i(\lambda)$ with $\alpha \prec_\theta \beta$. Then $\alpha$ and $\beta$ are \emph{adjacent} if when $\alpha \preceq_\theta \gamma \preceq_\theta \beta$, for $\gamma\in\mathcal{B}_\theta^i(\lambda)$, then $\gamma\in\{\alpha,\beta\}$. If two boxes are adjacent to each other we write $\alpha \prec_\theta^a \beta$.
\end{Definition}

\begin{Definition}[{\cite[Definition 3.1.1]{Jacon:Staggered}}]\label{D: Staggered sequence}
	Let $i\in I$. A $(\theta,i)$-\emph{staggered sequence} of $\lambda$ is a non-empty sequence $\mathbf{s} = (s_1,\ldots,s_m)$ of removable boxes in $\mathcal{B}_\theta^i(\lambda)$ such that
	\begin{enumerate}
		\item $s_i \prec_\theta^a s_{i+1}$, for $1\leq i < m$,
		\item $\alpha \preceq_\theta s_m$, for all $\alpha\in\mathcal{B}_\theta^i(\lambda)$,
		\item if there is an $\alpha\in\mathcal{B}_\theta^i(\lambda)$ with $\alpha\prec_\theta^a s_1$ then $\alpha\in\mathcal{A}_\theta^i(\lambda)$.
	\end{enumerate}
	We define $|\mathbf{s}| := m$ and $\residue(\mathbf{s}) := i$.
\end{Definition}

Staggered sequences depend on $\theta$ but we omit this from our notation. If $\mathbf{s} = (s_1,\ldots,s_m)$ is a staggered sequence of $\lambda$ then define $\lambda\setminus \mathbf{s} := \lambda\setminus\{s_1,\ldots,s_m\}$.  

Note that an arbitrary multipartition $\lambda$ does not need to have a staggered sequence for any $i\in I$. However, we now recursively define a subset of $\mathscr{P}_n^\ell(\theta)$ consisting of multipartitions that always have staggered sequences.

\begin{Definition}[{\cite[Definition 3.1.2]{Jacon:Staggered}}]\label{D: Uglov/Kleshchev}
	The set of $\theta$-\emph{Staggered multipartitions} $\mathscr{S}^\ell(\theta)$ is the subset of $\mathscr{P}^\ell(\theta)$ determined by $\lambda\in\mathscr{S}^\ell(\theta)$ if
	\begin{enumerate}
		\item $\lambda = \emptyset$ or
		\item $\lambda$ has a staggered sequence $\mathbf{s} = (s_1,\ldots,s_m)$ and $\lambda\setminus\mathbf{s}\in\mathscr{S}^\ell(\theta)$.
	\end{enumerate}
	Define $\mathscr{S}_n^\ell(\theta) = \mathscr{S}^\ell(\theta) \cap \mathscr{P}_n^\ell(\theta)$ the set of \emph{Staggered $\ell$-multipartitions of $n$}.
\end{Definition}

\begin{Theorem}[{\cite[Theorem 5.1.1]{Jacon:Staggered}}]\label{T: Jacon Uglov}
	As sets $\mathscr{S}^\ell(\theta) = \mathscr{U}^\ell(\theta)$.
\end{Theorem}

\begin{Theorem}\label{T: Number of simples}
	The number of simple $\underline{\mathscr{R}}_n^\Lambda(\theta)$-modules is $|\mathscr{S}_n^\ell(\theta)| = |\mathscr{U}_n^\ell(\theta)|$.
\end{Theorem}

\begin{proof}
	By \autoref{T: Jacon Uglov}, $\mathscr{S}_n^\ell(\theta) = \mathscr{U}_n^\ell(\theta)$ is the set of Uglov multipartitions as in \cite{Uglov}. In particular, $\mathscr{U}_n^\ell(\theta)$ and $\mathscr{K}_n^\ell$ both label the crystal graph of the highest weight module, hence, $|\mathscr{U}_n^\ell(\theta)| = |\mathscr{K}_n^\ell|$. Hence, in view of \autoref{T: BK isom} this proves the result because Ariki and Mathas \cite{AM:simples} showed that $\underline{\mathcal{H}}_n^\Lambda(q)$ has $|\mathscr{K}_n^\ell|$ simple modules.
\end{proof}

\noindent So for the remainder of the paper we use \autoref{D: Uglov/Kleshchev} as the definition of $\mathscr{U}^\ell(\theta)$.

\begin{Definition}
	Let $c\geq 2$ be an integer and $\lambda = (\lambda^{(1)}, \ldots, \lambda^{(\ell)})\in\mathscr{P}^\ell(\theta)$. Then $\lambda$ is \emph{$c$-restricted} if $\lambda_i^{(l)} - \lambda_{i+1}^{(l)} < c$, for all $i\geq 1$ and $1\leq l\leq\ell$.
\end{Definition}

The following is a well known result for Kleshchev multipartitions but we have not found it in the literature for Uglov multipartitions. However, the proof is short and is analogous to the Kleshchev case.

\begin{Lemma}\label{L: e-restricted}
	Let $\lambda\in\mathscr{U}_n^\ell(\theta)$. Then $\lambda$ is $e$-restricted.
\end{Lemma}

\begin{proof}
	We argue by induction on $|\lambda|=n$. If $n=1$, then there is nothing to show. So assume $n>1$ and the statement being true for $n-1$. Let $\mathbf{s} = (s_1,s_2,\ldots,s_m)$ be a staggered sequence of $\lambda$, for some $m\in\mathbb{N}$, then $\lambda' = \lambda\setminus s_1 \in\mathscr{U}_{n-1}^\ell(\theta)$ because $\mathbf{s}' = (s_2,\ldots,s_m)$ is a staggered sequence of $\lambda'$. By induction, $\lambda'$ is $e$-restricted. Assume $\lambda$ is not $e$-restricted and let $s_1 = (r,c,l)$. Then $\lambda_r^{(l)} - \lambda_{r+1}^{(l)} = e$. However, this implies $\alpha = (r+1,c-e,l)$ is an addable box of $\lambda$ with $\residue(\alpha) = \residue(s_1)$ and $s_1 \prec_\theta \alpha$. This is a contradiction to $\mathbf{s}$ being a staggered sequence.
\end{proof}

For the special case $\ell = 1$ the reverse implication is also true.

\begin{Theorem}[\cite{MisraMiwa}]\label{T: MisraMiwa}
	Let $\lambda\in\mathscr{P}_n^1(\theta)$. Then $\lambda\in\mathscr{U}_n^1(\theta)$ if and only if $\lambda$ is $e$-restricted.
\end{Theorem}

\begin{Definition}\label{D: Decomposition Path}
	Let $\lambda\in\mathscr{U}_n^\ell(\theta)$. The \emph{decomposition graph} $\Gamma(\lambda)$ of $\lambda$ is the edge labelled directed graph with vertex set $\{\mu\in\mathscr{U}^\ell(\theta) | \mu\subset\lambda\}$ and edges $\mu\xrightarrow{\textbf{s}}\nu$ if $\textbf{s}$ is a staggered sequence of $\mu$ and $\nu = \mu\setminus\textbf{s}$. A \emph{decomposition path} $\mathbf{d} = \mathbf{s}_1\mathbf{s}_2\ldots \mathbf{s}_k$ for $\lambda$ is a finite path
	\begin{equation*}
	\lambda = \lambda_0 \xrightarrow{\mathbf{s}_1} \lambda_1 \xrightarrow{\mathbf{s}_2} \ldots \xrightarrow{\mathbf{s}_k} \lambda_k = \emptyset
	\end{equation*}
	in $\Gamma(\lambda)$. The \emph{length} of a decomposition path $\mathbf{d} = \mathbf{s}_1\mathbf{s}_2\ldots \mathbf{s}_k$ is $|\mathbf{d}| = k$. 
\end{Definition}

\begin{Example}\label{E: Braid Move}
	Let $\lambda = (3,1,1)$, $e=3$. Then $\lambda$ has two decomposition paths that are encoded in $\Gamma(\lambda)$:
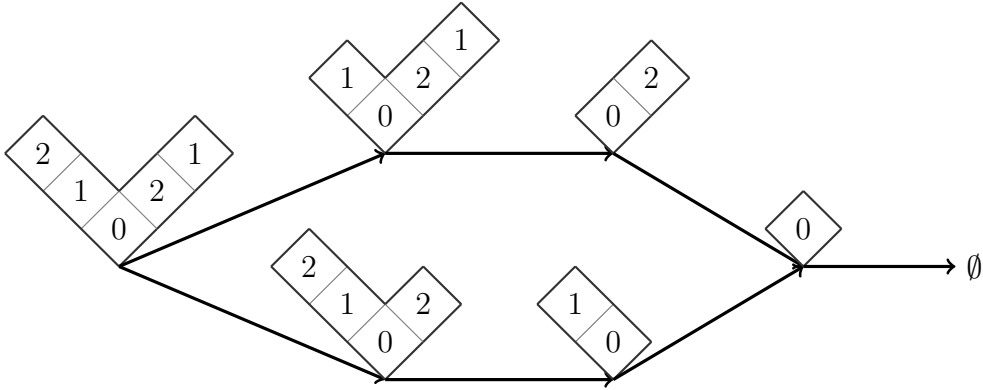
	\[
	\begin{tikzpicture}[very thick,baseline,scale=0.5]
	\draw (0,0) pic{russiantableau={scale=0.5}{{0,1,2},{2},{1}}};
	\draw[->] (0,0) -- (7,3);
	\draw[->] (0,0) -- (7,-3);

	\draw (7,3) pic{russiantableau={scale=0.5}{{0,1},{2},{1}}};
	\draw[->] (7,3) -- (13,3);

	\draw (7,-3) pic{russiantableau={scale=0.5}{{0,1,2},{2}}};
	\draw[->] (7,-3) -- (13,-3);

	\draw (13,3) pic{russiantableau={scale=0.5}{{0},{2}}};
	\draw[->] (13,3) -- (18,0);

	\draw (13,-3) pic{russiantableau={scale=0.5}{{0,1}}};
	\draw[->] (13,-3) -- (18,0);

	\draw (18,0) pic{russiantableau={scale=0.5}{{0}}};
	\draw[->] (18,0) -- (22,0);

	\node at (22.5,0) {$\emptyset$};
	\end{tikzpicture}
	\]
\end{Example}

\subsection{Decomposition tableaux}\label{S: Decomposition Tableaux}

This section, constructs special standard tableaux for Uglov multipartitions, that will be used for our computations of the inner product on the graded cell modules $S_\theta^\lambda$ in \autoref{S: Inner products}.

\begin{Definition}\label{D: Decomposition Tableaux}
	Let $\lambda$ be an Uglov multipartition and let $\mathbf{d} = \mathbf{s}_1\mathbf{s}_2\ldots \mathbf{s}_k$ be a decomposition path for $\lambda$. The \emph{positive decomposition tableau} $\mathfrak{t}_{\mathbf{d}}$ is the unique standard tableau of shape $\lambda$ such that
	\begin{enumerate}
		\item $\mathfrak{t}_{\mathbf{d}}(\alpha) < \mathfrak{t}_{\mathbf{d}}(\beta)$ if $\alpha,\beta\in \mathbf{s}_j$, for some $1\leq j\leq k$, and $\alpha \prec_\theta \beta$,
		\item $\mathfrak{t}_{\mathbf{d}}(\alpha) < \mathfrak{t}_{\mathbf{d}}(\beta)$ if $\alpha\in \mathbf{s}_j$ and $\beta\in \mathbf{s}_{j'}$ such that $1\leq j' < j \leq k$.
	\end{enumerate}
	Similarly, the \emph{negative decomposition tableau} $\mathfrak{t}_{\mathbf{d}}^\diamond$ is the unique standard tableau of shape $\lambda$ such that
	\begin{enumerate}
		\item $\mathfrak{t}_{\mathbf{d}}^\diamond(\alpha) > \mathfrak{t}_{\mathbf{d}}^\diamond(\beta)$ if $\alpha,\beta\in \mathbf{s}_j$, for some $1\leq j\leq k$, and $\alpha \prec_\theta \beta$,
		\item $\mathfrak{t}_{\mathbf{d}}^\diamond(\alpha) < \mathfrak{t}_{\mathbf{d}}^\diamond(\beta)$ if $\alpha\in \mathbf{s}_j$ and $\beta\in \mathbf{s}_{j'}$ such that $1\leq j' < j \leq k$.
	\end{enumerate}
\end{Definition}

\begin{Example}\label{E: decomposition tableaux}
	Let $\lambda = ((2,1),(3,1))$, $e=3$, $\theta=(0,8)$, and $\kappa=(0,0)$. Then there is a unique decomposition path $\mathbf{d} = \color{red}{\mathbf{s}_1}\color{blue}{\mathbf{s}_2}\color{green}{\mathbf{s}_3}$, where $\color{red}{\mathbf{s}_1}$ is of residue $2$, $\color{blue}{\mathbf{s}_2}$ of residue $1$, and $\color{green}{\mathbf{s}_3}$ of residue $0$.
	
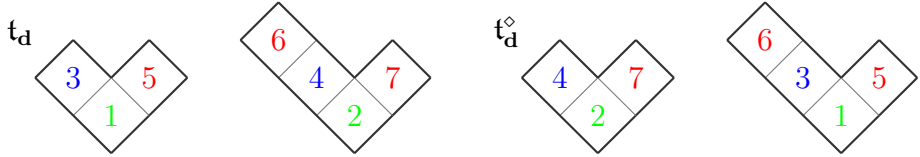
\[
	\begin{tikzpicture}[scale = 0.8]
	\node at (-1.5,2) {$\mathfrak{t}_{\mathbf{d}}$};

	\draw (0,0) pic{russiantableau={scale=0.5}{{\color{green}{1},\color{blue}{3}},{\color{red}{5}}}};
	\draw (4,0) pic{russiantableau={scale=0.5}{{\color{green}{2},\color{blue}{4},\color{red}{6}},{\color{red}{7}}}};

	\node at (6.5,2) {$\mathfrak{t}_{\mathbf{d}}^\diamond$};

	\draw (8,0) pic{russiantableau={scale=0.5}{{\color{green}{2},\color{blue}{4}},{\color{red}{7}}}};
	\draw (12,0) pic{russiantableau={scale=0.5}{{\color{green}{1},\color{blue}{3},\color{red}{6}},{\color{red}{5}}}};
	\end{tikzpicture}
	\]
\end{Example}

First we give some immediate consequences of this definition. Let $\mathbf{d}$ be a decomposition path of length $k$, for $k\in\mathbb{N}_0$, and $m_j = |\mathbf{s}_j|$, for $1\leq j\leq k$. The Young subgroup $\mathfrak{S}_\mathbf{d} := \prod_{j=1}^{k} \mathfrak{S}_{m_j}$ acts on the decomposition tableaux of $\mathbf{d}$ by
\begin{equation*}
(\mathfrak{t}_{\mathbf{d}}\cdot\sigma)(\alpha) := \sigma_j(\mathfrak{t}_{\mathbf{d}}(\alpha)) \quad\text{ and }\quad (\mathfrak{t}_{\mathbf{d}}^\diamond\cdot\sigma)(\alpha) := \sigma_j(\mathfrak{t}_{\mathbf{d}}^\diamond(\alpha)),
\end{equation*}
where $\sigma = (\sigma_1,\ldots,\sigma_k)\in\mathfrak{S}_\mathbf{d}$ and $\alpha\in\mathbf{s}_j$, for $1\leq j\leq k$. By definition, staggered sequences are removable boxes for their corresponding multipartition, so $\mathfrak{t}_{\mathbf{d}}\cdot\sigma$ and $\mathfrak{t}_{\mathbf{d}}^\diamond\cdot\sigma$ are standard tableaux, for $\sigma\in\mathfrak{S}_\mathbf{d}$. Let $\omega_0^{\mathbf{d}}\in \mathfrak{S}_\mathbf{d}$ be the \emph{longest element}. By definition, $\mathfrak{t}_{\mathbf{d}}\cdot\omega_0^{\mathbf{d}} = \mathfrak{t}_\mathbf{d}^\diamond$. The tableaux in the set $\{\mathfrak{t}_{\mathbf{d}}\cdot\sigma | \sigma\in \mathfrak{S}_\mathbf{d}\}$ are pairwise distinct.

\begin{Lemma}\label{L: Residue sequence}
	Let $\mathbf{d}$ and $\mathbf{e}$ be decomposition paths of $\lambda\in\mathscr{U}^\ell(\theta)$. Then for all $\sigma\in \mathfrak{S}_\mathbf{d}$, $\residue(\mathfrak{t}_\mathbf{d}) = \residue(\mathfrak{t}_\mathbf{d}\cdot\sigma)$ and $\residue(\mathfrak{t}_{\mathbf{d}}) = \residue(\mathfrak{t}_{\mathbf{e}}) \Longleftrightarrow \mathbf{d} = \mathbf{e}$.
\end{Lemma}

\begin{proof}
	Let $\mathbf{d} = \mathbf{s}_1\ldots \mathbf{s}_k$ and $\sigma = (\sigma_1,\ldots,\sigma_k)\in\mathfrak{S}_d$. $(\mathfrak{t}_{\mathbf{d}}\cdot\sigma)(\alpha) = \sigma_j(\mathfrak{t}_{\mathbf{d}}(\alpha))$, for all $\alpha\in\mathbf{s}_j$. But by definition, all boxes of $\mathbf{s}_j$ are boxes of the same residue. So the residue sequence of $\mathfrak{t}_{\mathbf{d}}$ does not change under the action of $\mathfrak{S}_d$.

	For the second part, if $\mathbf{d} = \mathbf{e}$ then $\residue(\mathfrak{t}_{\mathbf{d}}) = \residue(\mathfrak{t}_{\mathbf{e}})$ trivially. So assume $\mathbf{d} \neq \mathbf{e}$ and let $\mathbf{d} = \mathbf{s}_1\ldots \mathbf{s}_k$ and $\mathbf{e} = \mathbf{t}_1\ldots \mathbf{t}_l$, for some $k,l\in\mathbb{N}$. Then there exists $j\in\mathbb{N}$ minimal such that $\mathbf{s}_j\neq \mathbf{t}_j$ and $\mathbf{s}_i = \mathbf{t}_i$ for all $1\leq i<j$. But this implies the residue of the boxes of $\mathbf{s}_j$ are different from those in $\mathbf{t}_j$, as there is only one staggered sequence for a given residue. Hence, $\residue(\mathfrak{t}_{\mathbf{d}}) \neq \residue(\mathfrak{t}_{\mathbf{e}})$.
\end{proof}

\begin{Corollary}
	Let $\mathbf{d}$ and $\mathbf{e}$ be decomposition paths of $\lambda\in\mathscr{U}^\ell(\theta)$ and $\sigma\in\mathfrak{S}_{\mathbf{d}}$ and $\tau\in\mathfrak{S}_{\mathbf{e}}$. Then $\mathfrak{t}_{\mathbf{d}}\cdot\sigma = \mathfrak{t}_{\mathbf{e}}\cdot\tau$ if and only if $\mathbf{d} = \mathbf{e}$ and $\sigma = \tau$.
\end{Corollary}

Next we determine the degree of the tableaux of the form $\mathfrak{t}_{\mathbf{d}}\cdot\sigma$ for some decomposition path $\mathbf{d}$. Let $\ell$ be the usual \emph{length function} on $\mathfrak{S}_{\mathbf{d}}$. Recall \autoref{D: Degree function}, $\deg_\theta(\mathfrak{t}) = \sum_{j=1}^n \left(\left|\mathscr{A}_\mathfrak{t}^\theta(j)\right| - \left|\mathscr{R}_\mathfrak{t}^\theta(j)\right|\right)$, where
\begin{align*}
\mathscr{A}_\mathfrak{t}^\theta(j) = \{\alpha\in\mathcal{A}_\theta^{i} (\shape(\mathfrak{t}_{\downarrow j})) | \mathfrak{t}^{-1}(j) \prec_\theta \alpha\},\\
\mathscr{R}_\mathfrak{t}^\theta(j) = \{\rho\in\mathcal{R}_\theta^{i} (\shape(\mathfrak{t}_{\downarrow j})) | \mathfrak{t}^{-1}(j) \prec_\theta \rho\},
\end{align*}
and $i=\residue(\mathfrak{t}^{-1}(j))$, for some $i\in I$.

\begin{Lemma}\label{L: degree of tableaux}
	Let $\mathbf{d} = \mathbf{s}_1\ldots \mathbf{s}_k$ a decomposition path for $\lambda$. Then
	\begin{equation*}
	\deg(\mathfrak{t}_{\mathbf{d}}\cdot\sigma) = \ell(\omega_0^{\mathbf{d}}) - 2\ell(\sigma),
	\end{equation*}
	for all $\sigma\in\mathfrak{S}_{\mathbf{d}}$. In particular,
	\begin{equation*}
	\deg(\mathfrak{t}_{\mathbf{d}}\cdot\sigma) + \deg(\mathfrak{t}_{\mathbf{d}}^\diamond\cdot\sigma^{-1}) = 0.
	\end{equation*}
\end{Lemma}

\begin{proof}
	First we claim $\deg(\mathfrak{t}_{\mathbf{d}}) = \ell(\omega_0^{\mathbf{d}})$. Use induction on $|\mathbf{d}| = k$. Assume $k=1$ then $\mathbf{d} = \mathbf{s}$, where $\mathbf{s} = (s_1,\ldots,s_m)$ with $\residue(s_1) = \ldots = \residue(s_m) = i$, for some $i\in I$, and $\omega_0^\mathbf{d}\in\mathfrak{S}_\mathbf{d} = \mathfrak{S}_m$. By definition, $\mathfrak{t}_{\mathbf{d}}(s_j) < \mathfrak{t}_{\mathbf{d}}(s_{j+1})$ and $s_j \prec_\theta^a s_{j+1}$, for all $1\leq j<m$. Hence, $\mathfrak{t}_{\mathbf{d}}(s_j) = j$ and $\shape((\mathfrak{t}_{\mathbf{d}})_{\downarrow j}) = \{s_1,\ldots,s_j\}$, for $1\leq j\leq m$. By definition, the boxes of a staggered sequence are the rightmost removable boxes for a given residue. So $|\mathscr{A}_{\mathfrak{t}_\mathbf{d}}^\theta(j)| = m-j$ and $|\mathscr{R}_{\mathfrak{t}_\mathbf{d}}^\theta(j)| = 0$, for all $1\leq j\leq m$. Hence,
	\begin{equation*}
	\deg_\theta(\mathfrak{t}_\mathbf{d}) = \sum_{j=1}^{m} (m - j) = \sum_{j=1}^{m-1} (m - j) = \frac{m(m-1)}{2} = \ell(\omega_0^\mathbf{d}).
	\end{equation*}
	Now assume $k>1$, so $\mathbf{d} = \mathbf{s}_1\mathbf{s}_2\ldots\mathbf{s}_k$. Let $|\lambda| = n$, for some $n\in\mathbb{N}$. Let $\mathbf{s}_1 = (s_1,\ldots,s_m)$, for some $n>m\in\mathbb{N}$, and $\residue(s_1) = \ldots = \residue(s_m) = i$, for some $i\in I$. Let $\mathbf{d'} = \mathbf{s}_2\ldots\mathbf{s}_k$ and $\lambda' = \lambda\setminus\mathbf{s}_1$. Then by induction, $\deg_\theta(\mathfrak{t}_\mathbf{d'}) = \ell(\omega_0^\mathbf{d'})$, for $\omega_0^\mathbf{d'}\in\mathfrak{S}_\mathbf{d}$. But $\lambda' = \shape(\mathfrak{t}_\mathbf{d'}) = \shape((\mathfrak{t}_\mathbf{d})_{\downarrow n-m})$, so $\deg_\theta(\mathfrak{t}_\mathbf{d}) - \deg_\theta(\mathfrak{t}_\mathbf{d'}) = \sum_{j'=n-m+1}^{n} \left(\left|\mathscr{A}_{\mathfrak{t}_\mathbf{d}}^\theta(j')\right| - \left|\mathscr{R}_{\mathfrak{t}_\mathbf{d}}^\theta(j')\right|\right)$. Also $\ell(\omega_0^\mathbf{d}) - \ell(\omega_0^\mathbf{d'}) = \frac{m(m-1)}{2}$. Now similarly to the base case by definition, $\mathfrak{t}_{\mathbf{d}}(s_j) < \mathfrak{t}_{\mathbf{d}}(s_{j+1})$ and $s_j \prec_\theta^a s_{j+1}$, for all $1\leq j<m$. But also by definition, $\mathfrak{t}_{\mathbf{d}}(\alpha) < \mathfrak{t}_{\mathbf{d}}(s_j)$, for all $\alpha\in\lambda'$ and $1\leq j\leq m$. Hence, $\mathfrak{t}_{\mathbf{d}}(s_j) = n-m+j$ and $\shape((\mathfrak{t}_{\mathbf{d}})_{\downarrow n-m+j}) = \lambda' \cup \{s_1,\ldots,s_j\}$, for $1\leq j\leq m$. So $|\mathscr{A}_{\mathfrak{t}_\mathbf{d}}^\theta(n-m+j)| = m-j$ and $|\mathscr{R}_{\mathfrak{t}_\mathbf{d}}^\theta(n-m+j)| = 0$, for all $1\leq j\leq m$. Hence,
	\begin{equation*}
	\deg_\theta(\mathfrak{t}_\mathbf{d}) = \deg_\theta(\mathfrak{t}_\mathbf{d'}) + \sum_{j' = n-m+1}^{n} |\mathscr{A}_{\mathfrak{t}_\mathbf{d}}^\theta(j')| = \ell(\omega_0^\mathbf{d'}) + \sum_{j = 1}^{m} (m-j) = \ell(\omega_0^\mathbf{d}).
	\end{equation*}

	By a similar argument, $\deg(\mathfrak{t}_{\mathbf{d}}^\diamond) = - \ell(\omega_0^{\mathbf{d}}) = - \deg(\mathfrak{t}_{\mathbf{d}})$. So now determine how a simple reflection $(j,j+1)\in\mathfrak{S}_\mathbf{d}$, with $\ell(\sigma\cdot (j,j+1)) > \ell(\sigma)$, changes the degree of the element $\mathfrak{t}_\mathbf{d}\cdot\sigma$. $(j,j+1)\in\mathfrak{S}_\mathbf{d}$ implies that the boxes $\mathfrak{t}_\mathbf{d}^{-1}(j)$ and $\mathfrak{t}_\mathbf{d}^{-1}(j+1)$ are in the same staggered sequence. It is easy to see that $\deg(\mathfrak{t}_\mathbf{d}\cdot\sigma\cdot (j,j+1)) = \deg(\mathfrak{t}_\mathbf{d}\cdot\sigma) - 2$; see for example \cite[Proposition 3.13]{BKW:GradedSpecht}. The result now follows because $\ell(\omega_0^{\mathbf{d}}) = \deg(\mathfrak{t}_\mathbf{d})$ and $\deg(\mathfrak{t}_\mathbf{d}\cdot\omega_0^{\mathbf{d}}) = \deg(\mathfrak{t}_\mathbf{d}^\diamond) = - \deg(\mathfrak{t}_\mathbf{d})$.
\end{proof}

Let $\mathbf{d} = \mathbf{s}_1\mathbf{s}_2\ldots \mathbf{s}_k$ be a decomposition path for $\lambda\in\mathscr{U}_n^\ell(\theta)$, where $\mathbf{s}_1 = (s_1,\ldots,s_m)$. By definition, $\mathfrak{t}_\mathbf{d}^\diamond(\beta) \leq \mathfrak{t}_\mathbf{d}^\diamond(s_1) = n$, for all $\beta\in\lambda$. Then $\mathbf{s}_1' = \mathbf{s}_1\setminus s_1$ is a staggered sequence of $\lambda \setminus s_1$ or is empty if $m=1$ and so $\mathbf{d'} = \mathbf{s}_1'\mathbf{s}_2\ldots \mathbf{s}_k$ is a decomposition path for $\lambda \setminus s_1$ with negative decomposition tableau $(\mathfrak{t}_\mathbf{d}^\diamond)_{\downarrow (n-1)}$.

\begin{Lemma}\label{L: negative decomposition tableaux property}
	Let $\mathbf{d} = \mathbf{s}_1\mathbf{s}_2\ldots \mathbf{s}_k$ be a decomposition path for $\lambda\in\mathscr{U}_n^\ell(\theta)$. Let $\alpha\in\lambda$ such that $\mathfrak{t}_\mathbf{d}^\diamond(\alpha) = n$ and $\mathbf{d'} = \mathbf{s}_1'\mathbf{s}_2\ldots \mathbf{s}_k$, where $\mathbf{s}_1' = \mathbf{s}_1\setminus \alpha$. Then $\mathfrak{t}_{\mathbf{d'}}^\diamond = (\mathfrak{t}_\mathbf{d}^\diamond)_{\downarrow (n-1)}$.
\end{Lemma}

\section{Classification and Applications}

In this section we use the standard tableaux of shape $\lambda\in\mathscr{U}_n^\ell(\theta)$ defined in \autoref{S: Decomposition paths}. They allow us to state the key theorem of this section. This theorem then implies the classification of the graded simple $\mathscr{R}_n^\Lambda(\theta)$-modules and gives us information about the graded dimensions of the graded simple $\mathscr{R}_n^\Lambda(\theta)$-modules.

We have seen that $\mathscr{R}_n^\Lambda(\theta)$ is a graded cellular algebra with graded cellular basis $\{c_{\mathfrak{st}}^\lambda | \lambda\in\mathscr{P}_n^\ell(\theta) \text{ and } \mathfrak{s},\mathfrak{t}\in\Std(\lambda)\}$. Its graded cell modules are $S_\theta^\lambda$ with basis $\{c_\mathfrak{t}| \lambda\in\mathscr{P}_n^\ell(\theta) \text{ and } \mathfrak{t}\in\Std(\lambda)\}$ and its simple modules are $D_\theta^\lambda$, for $\lambda\in\mathscr{P}_n^\ell(\theta)$ with $D_\theta^\lambda \neq 0$. By \autoref{D: Simple}, if we find $x,y\in S_\theta^\lambda$ such that $\langle x,y\rangle_\lambda \neq 0$ then $D_\theta^\lambda \neq 0$.

We are now able to state the following theorem. The proof of this theorem, which requires diagram calculations, will be done in \autoref{S: Inner products}.

\begin{Theorem}\label{T: Inner Product}
	Let $\lambda\in\mathscr{U}^\ell(\theta)$ be an Uglov multipartition and $\mathbf{d}$ a decomposition path for $\lambda$. Then $|\langle c_{\mathfrak{t}_\mathbf{d}}^\lambda, c_{\mathfrak{t}_\mathbf{d}^\diamond}^\lambda\rangle_\lambda| = 1$. Consequently, $D_\theta^{\lambda} \neq 0$.
\end{Theorem}

\subsection{Classification and DJM Conjecture}

In this section we use \autoref{T: Inner Product} to classify the graded simple $\mathscr{R}_n^\Lambda(\theta)$-modules. We also use the tableaux defined in \autoref{S: Decomposition Tableaux} to give an alternative proof of a recent result by Jacon \cite[Corollary 10.1.2]{Jacon:Staggered} on a known conjecture.

\begin{Corollary}\label{T: Classification}
	Let $\lambda\in\mathscr{P}_n^\ell(\theta)$. Then $D_\theta^{\lambda}$ is simple if and only if $\lambda\in\mathscr{U}_n^\ell(\theta)$.
\end{Corollary}

\begin{proof}
	If $\lambda\in\mathscr{U}_n^\ell(\theta)$ then $D_\theta^{\lambda}$ is simple, by \autoref{T: Inner Product}. On the other hand, by \autoref{T: Number of simples}, the number of simple $\underline{\mathscr{R}}_n^\Lambda(\theta)$-modules is $|\mathscr{U}_n^\ell(\theta)|$. Hence, the result follows.
\end{proof}

\noindent Now we can prove our main result:

\begin{Main Theorem}\label{T: Main Theorem}
	Let $\lambda\in\mathscr{P}_n^\ell(\theta)$. Then the following are equivalent:
	\begin{enumerate}
		\item $D_\theta^\lambda\neq 0$
		\item $\lambda\in\mathscr{U}_n^\ell(\theta)$
		\item $|\langle c_{\mathfrak{t}_{\mathbf{d}}}^\lambda, c_{\mathfrak{t}_{\mathbf{d}}^\diamond}^\lambda\rangle_\lambda| = 1$, for any decomposition path $\mathbf{d}$ of $\lambda$
	\end{enumerate}
\end{Main Theorem}

\begin{proof}
	By \autoref{T: Classification}, parts (a) and (b) are equivalent and, by \autoref{T: Inner Product} (b) implies (c). On the other hand, (c) implies (a) by \autoref{D: Simple}.
\end{proof}

\begin{Corollary}[Classification]\label{C: Main Classification}
	The set $\{D_\theta^\lambda\langle x\rangle | \lambda\in\mathscr{U}_n^\ell(\theta), \text{ and } x\in\mathbb{Z}\}$ is a complete set of pairwise non-isomorphic simple $\mathscr{R}_n^\Lambda(\theta)$-modules.
\end{Corollary}

\begin{Definition}
	A standard tableau $\mathfrak{t}$ is a \emph{Dipper-James-Murphy tableau}, or a \emph{DJM-tableau}, if whenever there exists a standard tableau $\mathfrak{s}$ such that $\residue(\mathfrak{s}) = \residue(\mathfrak{t})$, then $\shape(\mathfrak{s}) \ntriangleleft_\theta \shape(\mathfrak{t})$.
\end{Definition}

In type B, Dipper, James, and Murphy \cite[Conjecture 8.13]{DJM} conjectured that $D^\lambda \neq 0$ if and only if there exists a DJM-tableau of shape $\lambda$. Graham and Lehrer \cite[Question 5.10 (i)]{GL} asked if this holds in general for the simple modules of the Ariki-Koike algebras (when $\theta$ is well-separated).

Let $\mathbf{d} = \mathbf{s}_1\ldots \mathbf{s}_k$ be a decomposition path for $\lambda\in\mathscr{U}_n^\ell(\theta)$. Assume there is a standard tableau $\mathfrak{s}$ such that $\residue(\mathfrak{s}) = \residue(\mathfrak{t}_{\mathbf{d}}^\diamond)$. Then for every staggered sequence $\mathbf{s}_j$ of $\mathbf{d}$, for $j=1,\ldots,k$, there exists a corresponding set of boxes $\mathbf{r}_j \subset \shape(\mathfrak{s})$ such that $\mathfrak{t}_{\mathbf{d}}^\diamond(\mathbf{s}_j) = \mathfrak{s}(\mathbf{r}_j)$, for $1\leq j\leq k$. The disjoint union of sets of boxes $\bigcup_{j = 1}^k \mathbf{r}_j = \shape(\mathfrak{s})$ is the \emph{partitioning of $\mathfrak{s}$ along $\mathbf{d}$}.

Unless stated otherwise, we assume for the rest of this section that $\theta$ is well-separated. Recall from \autoref{D: Uglov/Kleshchev} that the set of Kleshchev multipartitions is $\mathscr{K}_n^\ell = \mathscr{U}_n^\ell(\theta)$, in this case. 

\begin{Lemma}\label{L: Kleshchev all t_d are DJM}
	Let $\lambda\in\mathscr{K}_n^\ell$ and $\mathbf{d}$ a decomposition path for $\lambda$. Then $\mathfrak{t}_{\mathbf{d}}\cdot\sigma$ is a DJM-tableau, for all $\sigma\in\mathfrak{S}_\mathbf{d}$.
\end{Lemma}	

\begin{proof}
	By \autoref{L: Residue sequence}, $\residue(\mathfrak{t}_{\mathbf{d}}\cdot\sigma) = \residue(\mathfrak{t}_{\mathbf{d}}^\diamond)$, for all $\sigma\in\mathfrak{S}_\mathbf{d}$. Moreover, by definition, $\shape(\mathfrak{t}_{\mathbf{d}}\cdot\sigma) = \shape(\mathfrak{t}_{\mathbf{d}}^\diamond) = \lambda$, for all $\sigma\in\mathfrak{S}_\mathbf{d}$. Hence, if $\mathfrak{t}_{\mathbf{d}}^\diamond$ is a DJM-tableau then so is $\mathfrak{t}_{\mathbf{d}}\cdot\sigma$, for all $\sigma\in\mathfrak{S}_\mathbf{d}$.
	
	So let $\mathfrak{s}$ be a standard tableau such that $\residue(\mathfrak{s}) = \residue(\mathfrak{t}_{\mathbf{d}}^\diamond)$ and let $\mu = \shape(\mathfrak{s})$. If $\mathbf{d} = \mathbf{s}_1\ldots \mathbf{s}_k$ then let $\bigcup_{j = 1}^k \mathbf{r}_j$ be the partitioning of $\mathfrak{s}$ along $\mathbf{d}$. Now $\mathfrak{S}_\mathbf{d}$ operates on $\mathfrak{s}$ and the orbit $\mathfrak{s}\cdot\mathfrak{S}_\mathbf{d}$ is a set of standard tableaux of shape $\mu$ with residue $\residue(\mathfrak{s})$. WLOG we can assume $\mathfrak{s}\in \mathfrak{s}\cdot\mathfrak{S}_\mathbf{d}$ such that $\mathfrak{s}(\alpha) > \mathfrak{s}(\beta)$ if $\alpha,\beta\in\mathbf{r}_j$, for some $1\leq j\leq k$, and $\alpha \prec_\theta \beta$, because all tableaux in the set $\mathfrak{s}\cdot\mathfrak{S}_\mathbf{d}$ have the same shape. We are setting up a claim proven by induction on $|\lambda| = n$ and in light of \autoref{L: negative decomposition tableaux property} we have $\mathfrak{t}_{\mathbf{d'}}^\diamond = (\mathfrak{t}_\mathbf{d}^\diamond)_{\downarrow (n-1)}$ and $\mathfrak{s}_{\downarrow (n-1)}\in \mathfrak{s}_{\downarrow (n-1)}\cdot\mathfrak{S}_\mathbf{d}$ such that $\mathfrak{s}_{\downarrow (n-1)}(\alpha) > \mathfrak{s}_{\downarrow (n-1)}(\beta)$ if $\alpha,\beta\in\mathbf{r}'_j$, for some $1\leq j\leq k$, and $\alpha \prec_\theta \beta$. This allows us the inductive step. 
	
	Now recall \autoref{D: Dominance order} and the following note. So $\mu \trianglelefteq_\theta \lambda$ if and only if $\mu\setminus\lambda \trianglelefteq_\theta \lambda\setminus\mu$. In particular, $\mu \trianglelefteq_\theta \lambda$ if and only if $(\mu\setminus\lambda)_i \trianglelefteq_\theta (\lambda\setminus\mu)_i$, for all $i\in I$. Because $\residue(\mathfrak{s}) = \residue(\mathfrak{t}_{\mathbf{d}}^\diamond)$ we have $|(\mu\setminus\lambda)_i| = |(\lambda\setminus\mu)_i|$, for all $i\in I$. So we also have $\mu \trianglelefteq_\theta \lambda$ if and only if for all $\alpha\in\mathscr{N}_n^\ell$
	\begin{equation*}
	|\{\beta\in\mu\setminus\lambda | \residue(\beta) = \residue(\alpha), \alpha \preceq_\theta \beta\}| \geq |\{\beta\in\lambda\setminus\mu | \residue(\beta) = \residue(\alpha), \alpha \preceq_\theta \beta\}|.
	\end{equation*}
	
	If $\mu = \lambda$ then there is nothing to show, hence, let $\mu \neq \lambda$. We now claim that if $\mu \neq \lambda$ then there exists $\gamma\in\lambda\setminus\mu$ such that $\delta \prec_\theta \gamma$, for all $\delta\in\mu\setminus\lambda$. If this claim is true then there exists $\residue(\gamma) = i\in I$ such that $(\mu\setminus\lambda)_i \ntriangleleft_\theta (\lambda\setminus\mu)_i$ and, thus, $\mu \ntriangleleft_\theta \lambda$. So to finish the proof we show this claim by induction on $|\lambda| = n$.
	
	If $n = 1$ then $\mathbf{d} = \mathbf{s}_1$ with $\mathbf{s}_1 = \{\alpha\}$ and $\mathbf{r}_1 = \{\beta\}$, and $\alpha \neq \beta$ with $\residue(\alpha) = \residue(\beta)$. By \autoref{D: Staggered sequence}, $\beta \prec_\theta \alpha$ and the claim follows. So now let $n > 1$ and assume the claim be true for up to $n-1$. Set $\alpha = (\mathfrak{t}_{\mathbf{d}}^\diamond)^{-1}(n)$ and $\beta = \mathfrak{s}^{-1}(n)$ then $\residue(\alpha) = \residue(\beta)$. If $\lambda\setminus\alpha = \mu\setminus\beta$ then again, by \autoref{D: Staggered sequence}, $\beta \prec_\theta \alpha$ and the claim follows. Hence, let $\lambda\setminus\alpha \neq \mu\setminus\beta$ then by induction there exists $\gamma\in (\lambda\setminus\alpha) \setminus (\mu\setminus\beta)$ such that $\delta \prec_\theta \gamma$, for all $\delta\in (\mu\setminus\beta) \setminus (\lambda\setminus\alpha)$. WLOG we can assume $\gamma' \preceq_\theta \gamma$, for all $\gamma' \in (\lambda\setminus\alpha) \setminus (\mu\setminus\beta)$. Then 
	\begin{equation}\label{E: DJM Lemma}
	\{\delta\in\lambda\setminus\alpha | \gamma \prec_\theta \delta\} = \{\delta\in\mu\setminus\beta | \gamma \prec_\theta \delta\}.
	\end{equation}
	Let $(r,c,l) = \gamma$, for some $r,c\in\mathbb{N}$ and $1\leq l\leq \ell$. If $\beta \prec_\theta \gamma$ then there is nothing else to show. If $\gamma \prec_\theta \beta$ then $\beta = (r',c',l')$ with $l \leq l'$ and $\beta$ is an addable box in $\lambda\setminus\alpha$, as $\theta$ is well-separated. Then, by \autoref{D: Staggered sequence}, $\beta \preceq_\theta \alpha$ and the claim follows. So it remains to check $\gamma = \beta$. Then if $\beta \prec_\theta \alpha$ the claim follows again, so let $\alpha \prec_\theta \beta = \gamma$.
	
	If $\gamma\in\lambda$ is not removable in $\lambda$ then $(r,c+1,l)\in\lambda\setminus\mu$, by \autoref{E: DJM Lemma} and $\gamma\notin\mu\setminus\beta$. The claim then follows with $(r,c+1,l)$. So assume $\gamma\in\lambda$ is removable in $\lambda$. Then, as $\residue(\gamma) = \residue(\alpha)$, we have $\gamma\in\mathbf{s}_1$ and $|\mathbf{s}_1| = m > 1$, by \autoref{D: Staggered sequence}. But then $|\mathbf{r}_1| = m > 1$. However, this implies $|\{\delta\in\mathbf{s}_1 | \gamma \prec_\theta \delta\}| \leq m-2$ but $|\{\delta\in\mathbf{r}_1 | \gamma \prec_\theta \delta\}| = m-1$. So
	\begin{equation*}
	|\{\delta\in\mathbf{s}_1 | \gamma \prec_\theta \delta\}| < |\{\delta\in\mathbf{r}_1 | \gamma \prec_\theta \delta\}|
	\end{equation*} 
	which is a contradiction to \autoref{E: DJM Lemma}, as it implies a different amount of removable $\residue(\alpha)$-boxes in $\lambda$ and $\mu$ to the right of $\gamma$. This completes the proof.
\end{proof}

\noindent Before we get to prove the actual conjecture we will prove a useful Corollary which is also implied by the previous Lemma.

\begin{Corollary}\label{C: Kleshchev t_d orbit is full}
	Let $\lambda\in\mathscr{K}_n^\ell$, $\mathbf{d}$ a decomposition path for $\lambda$, and $\mathfrak{s}$ a standard tableau with $\residue(\mathfrak{s}) = \residue(\mathfrak{t}_\mathbf{d}^\diamond)$. If $\shape(\mathfrak{s}) = \lambda$ then there exists $\sigma\in\mathfrak{S}_\mathbf{d}$ such that $\mathfrak{s} = \mathfrak{t}_\mathbf{d}^\diamond\cdot \sigma$.
\end{Corollary}

\begin{proof}
	We prove by induction on $|\lambda| = n$. The case $n=1$ is trivial. So let $n > 1$ and assume the statement is true for up to $n-1$. If $\mathbf{d} = \mathbf{s}_1\ldots \mathbf{s}_k$ then let $\bigcup_{j = 1}^k \mathbf{r}_j$ be the partitioning of $\mathfrak{s}$ along $\mathbf{d}$. Let $|\mathbf{s}_1| = |\mathbf{r}_1| = m\in\mathbb{N}$. WLOG we assume, similarly as in the previous proof, $\mathfrak{s}\in \mathfrak{s}\cdot\mathfrak{S}_\mathbf{d}$ such that $\mathfrak{s}(\gamma) > \mathfrak{s}(\delta)$ if $\gamma,\delta\in\mathbf{r}_1$, for $\gamma \prec_\theta \delta$.
	
	Let $\alpha = (\mathfrak{t}_{\mathbf{d}}^\diamond)^{-1}(n)$ and $\beta = \mathfrak{s}^{-1}(n)$. Then $\beta \preceq_\theta \alpha$ because $|\mathbf{s}_1| = |\mathbf{r}_1|$, as $\alpha \prec_\theta \beta$ implies
	\begin{equation*}
	m-1 = |\{\delta\in\mathbf{s}_1 | \alpha\prec_\theta\delta\}| < |\{\delta\in\mathbf{r}_1 | \alpha\prec_\theta\delta\}| = m,
	\end{equation*}
	which is a contradiction because $|\{\delta\in\mathbf{s}_1 | \alpha\prec_\theta\delta\}|$ is also the number of removable $\residue(\alpha)$-boxes to the right of $\alpha$ in $\lambda$, by \autoref{D: Staggered sequence}. 
	
	If $\beta = \alpha$ then $\shape(\mathfrak{s}_{\downarrow (n-1)}) = \lambda\setminus\beta = \lambda\setminus\alpha = \shape((\mathfrak{t}_\mathbf{d}^\diamond)_{\downarrow (n-1)})$ and $\residue(\mathfrak{s}_{\downarrow (n-1)}) = \residue((\mathfrak{t}_\mathbf{d}^\diamond)_{\downarrow (n-1)})$. So, by induction and \autoref{L: negative decomposition tableaux property}, $\mathfrak{s}_{\downarrow (n-1)}\in\{\mathfrak{t}_{\mathbf{d}'}^\diamond\cdot \sigma' | \sigma'\in\mathfrak{S}_{\mathbf{d}'}\}$, where $\mathbf{d}' = \mathbf{s}_1'\mathbf{s}_2\ldots \mathbf{s}_k$ with $\mathbf{s}_1' = \mathbf{s}_1\setminus\alpha$. Let $\sigma'\in\mathfrak{S}_{\mathbf{d}'}$ such that $\mathfrak{s}_{\downarrow (n-1)} = \mathfrak{t}_{\mathbf{d}'}\cdot \sigma'$. Set $\sigma\in\mathfrak{S}_n$ with $\sigma(j) = \sigma'(j)$, for $1\leq j\leq n-1$, and $\sigma(n) = n$. Then $\sigma\in\mathfrak{S}_\mathbf{d}$ and $\mathfrak{t}_\mathbf{d}^\diamond\cdot \sigma = \mathfrak{s}$.
	
	Now let $\beta \prec_\theta \alpha$. Then $\shape(\mathfrak{s}_{\downarrow (n-1)}) = \lambda\setminus\beta \vartriangleleft_\theta \lambda\setminus\alpha = \shape((\mathfrak{t}_\mathbf{d}^\diamond)_{\downarrow (n-1)})$ and $\residue(\mathfrak{s}_{\downarrow (n-1)}) = \residue((\mathfrak{t}_\mathbf{d}^\diamond)_{\downarrow (n-1)})$. However, this is a contradiction to \autoref{L: Kleshchev all t_d are DJM}.
\end{proof}

\noindent Jacon \cite{Jacon:Staggered} proved the following result. We give an alternative proof here.

\begin{Theorem}[Generalized Dipper-James-Murphy Conjecture; see {\cite[Corollary 10.1.2]{Jacon:Staggered}}]\label{T: DJM}
	Let $\lambda\in\mathscr{P}_n^\ell(\theta)$ and $\theta$ well-separated. Then $D_\theta^\lambda \neq 0$ if and only if there exists a DJM-tableau of shape $\lambda$.
\end{Theorem}

\begin{proof}
	Let $\mathfrak{t}$ be a DJM-tableau. By \cite[Theorem 5.10]{BK:GradedDecomp}, the only composition factors in $S_\theta^\lambda$ are simple modules $D_\theta^\mu$ with $\mu \trianglelefteq_\theta \lambda$. Assuming $D_\theta^\lambda = 0$ then implies there exists $\mu \trianglelefteq_\theta \lambda$ and $\mathfrak{s}\in\Std(\mu)$ with $\residue(\mathfrak{s}) = \residue(\mathfrak{t})$, a contradiction to our assumption. Hence, $D_\theta^\lambda \neq 0$.
	
	Now let $D_\theta^\lambda \neq 0$. By \autoref{T: Classification}, this implies $\lambda\in\mathscr{K}_n^\ell$. So there exists a decomposition path $\mathbf{d}$ for $\lambda$ and, by \autoref{L: Kleshchev all t_d are DJM}, $\mathfrak{t}_\mathbf{d}$ is a DJM-tableau of shape $\lambda$.
\end{proof}

\noindent Unfortunately, \autoref{L: Kleshchev all t_d are DJM} does not hold for general $\theta$, as the following example shows.

\begin{Example}
	Let $e = 5, \ell = 3, n = 11, \theta = (0,20,31)$, and $\kappa = (2,3,0)$. Let $\lambda = ((1,1,1),(1),(5,1,1))$ and $\mathbf{d}$ such that $\residue(\mathfrak{t}_\mathbf{d}^\diamond) = (0,1,2,3,2,1,0,4,4,3,3)$. Then there exists $\mathfrak{s}$ with $\residue(\mathfrak{s}) = \residue(\mathfrak{t}_\mathbf{d}^\diamond)$:
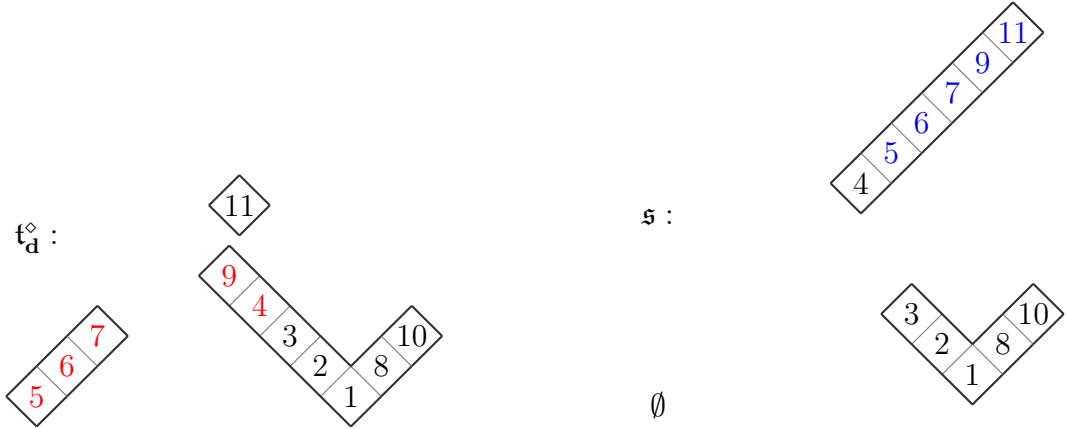
	\[
	\begin{tikzpicture}
	\node at (0,2.5) {$\mathfrak{t}_\mathbf{d}^\diamond:$};
	\draw (0,0) pic{russiantableau={scale=.4}{{\color{red}{5}},{\color{red}{6}},{\color{red}{7}}}};
	\draw (4.1333,0) pic{russiantableau={scale=.4}{{1,2,3,\color{red}{4},\color{red}{9}},{8},{10}}};
	\draw (2.667,2.5333) pic{russiantableau={scale=.4}{{11}}};
	\end{tikzpicture}\qquad\qquad\qquad
	\begin{tikzpicture}
	\node at (0,2.5) {$\mathfrak{s}:$};
	\node at (0,0) {$\emptyset$};
	\draw (4.1333,0) pic{russiantableau={scale=.4}{{1,2,3},{8},{10}}};
	\draw (2.667,2.5333) pic{russiantableau={scale=.4}{{4},{\color{blue}{5}},{\color{blue}{6}},{\color{blue}{7}},{\color{blue}{9}},{\color{blue}{11}}}};
	\end{tikzpicture}
	\]
	Highlighted in red and blue are the numbers in the boxes of $\lambda\setminus\shape(\mathfrak{s})$ and $\shape(\mathfrak{s})\setminus\lambda$ respectively. This clearly shows $(\shape(\mathfrak{s})\setminus\lambda)_i \vartriangleleft_\theta (\lambda\setminus\shape(\mathfrak{s}))_i$, for $i = 0,\ldots,4$. Thus, $\shape(\mathfrak{s}) \vartriangleleft_\theta \lambda$.
\end{Example}

\noindent Based on the example above one can also construct a counterexample that shows \autoref{C: Kleshchev t_d orbit is full} does not hold for general $\theta$ either. Thus, for the purposes of this paper, as we are more interested in \autoref{C: Kleshchev t_d orbit is full} and its implication on the graded dimensions of the graded simple modules than in \autoref{T: DJM} itself, we do not explore the conjecture for general $\theta$ any further here.

\subsection{Graded dimension of simple modules}

In this section we compute a non-trivial lower bound on the graded dimension of $D_\theta^\lambda$, as an application of \autoref{C: Kleshchev t_d orbit is full}. First we need the following technical Proposition.

\begin{Proposition}\label{P: Diagonal Gram matrix}
	Let $\lambda\in\mathscr{U}_n^\ell(\theta)$, $\mathbf{d}$ a decomposition path for $\lambda$, and $\sigma,\tau\in\mathfrak{S}_{\mathbf{d}}$. Then $|\langle c_{\mathfrak{t}_{\mathbf{d}}\cdot\sigma}^\lambda, c_{\mathfrak{t}_{\mathbf{d}}^\diamond\cdot\tau}^\lambda\rangle_\lambda| = \delta_{\sigma^{-1},\tau}$.
\end{Proposition}

The proof of this statement requires technical diagram calculations, so we defer the proof to \autoref{S: Inner products}. \autoref{P: Diagonal Gram matrix} implies the following result.

\begin{Corollary}\label{L: linear independence}
	Let $\lambda\in\mathscr{U}_n^\ell(\theta)$ and $\mathbf{d}$ a decomposition path for $\lambda$. Then the elements in the set $\{c_{\mathfrak{t}_{\mathbf{d}}^\diamond\cdot \sigma}^\lambda + \radic(S_\theta^\lambda) \in D_\theta^\lambda | \sigma\in\mathfrak{S}_{\mathbf{d}}\}$ are linearly independent.
\end{Corollary}

\begin{proof}
	Assume we have scalars $s_\sigma\in K$ such that
	\begin{equation*}
	\sum_{\sigma\in\mathfrak{S}_{\mathbf{d}}} s_\sigma c_{\mathfrak{t}_{\mathbf{d}}^\diamond\cdot \sigma}^\lambda \in \radic(S_\theta^\lambda).
	\end{equation*}
	Note that $c_{\mathfrak{t}_{\mathbf{d}}^\diamond\cdot \sigma}^\lambda \notin \radic(S_\theta^\lambda)$, by \autoref{P: Diagonal Gram matrix}. Let $\tau\in\mathfrak{S}_{\mathbf{d}}$. Taking the inner product with $c_{\mathfrak{t}_{\mathbf{d}}\cdot \tau^{-1}}^\lambda$ shows that $s_\tau = 0$, by \autoref{P: Diagonal Gram matrix}. Hence, $s_\sigma = 0$, for all $\sigma\in\mathfrak{S}_{\mathbf{d}}$, and the elements are linearly independent.
\end{proof}

Let $q\in\mathbb{N}_0$ and $t$ an indeterminate over $\mathbb{N}_0$. We define the polynomial
\begin{equation*}
[q]_t = \sum_{i=0}^{q-1} t^{-(q-1)+2i} = t^{-(q-1)} + t^{-(q-1)+2} + \ldots + t^{(q-1)-2} + t^{q-1}
\end{equation*}
and set $[q]_t^! = \prod_{i=1}^{q} [i]_t = [1]_t[2]_t\ldots [q]_t$. Let $p,\hat{p}\in\mathbb{N}_0[t]$ and write $p\leq \hat{p}$ if $\hat{p}-p\in\mathbb{N}_0[t]$. Then $p = \hat{p}$ if and only if $p\leq \hat{p}$ and $p\geq \hat{p}$. Note that $[q]_1 = q$ and $[q]_1^! = q!$.

\begin{Theorem}\label{T: Lower bound}
	Let $\lambda\in\mathscr{U}_n^\ell(\theta)$ and $\mathbf{d} = \mathbf{s}_1\ldots \mathbf{s}_k$ a decomposition path for $\lambda$ with $|\mathbf{s}_j| = m_j$. Then
	\begin{equation*}
	\dim_t(D_\theta^\lambda \cdot e_\omega(\residue(\mathfrak{t}_{\mathbf{d}}^\diamond))) \geq \prod_{j=1}^{k} [m_j]_t^!
	\end{equation*}
	Moreover, if $\lambda\in\mathscr{K}_n^\ell$ then we have equality.
\end{Theorem}

\begin{proof}
	By \autoref{L: linear independence}, for any decomposition path $\mathbf{d} = \mathbf{s}_1\ldots \mathbf{s}_k$ of $\lambda$ there is a set $\{c_{\mathfrak{t}_\mathbf{d}\cdot \sigma}^\lambda + \radic(S_\theta^\lambda) |\sigma\in\mathfrak{S}_\mathbf{d}\}$ of linearly independent elements in $D_\theta^\lambda$ of size $|\mathfrak{S}_{\mathbf{d}}| = \prod_{j=1}^{k} m_j!$. These elements have the same residue sequence as $\mathfrak{t}_{\mathbf{d}}$, by \autoref{L: Residue sequence}. By \autoref{L: degree of tableaux} and \autoref{C: R is Cellular Alg}, the degree of these elements is $\deg_\theta(c_{\mathfrak{t}_\mathbf{d}\cdot \sigma}^\lambda + \radic(S_\theta^\lambda)) = \deg_\theta(\mathfrak{t}_\mathbf{d}\cdot \sigma) = \ell(\omega_0^\mathbf{d}) - 2\ell(\sigma)$.

	The following equation is a variation of a well-known result about Poincare polynomials; see for example \cite[Section 1.11]{Hump}
	\begin{equation*}
	\sum_{\sigma\in\mathfrak{S}_m} t^{\ell(\omega_0) - 2\ell(\sigma)} = [m]_t^!,
	\end{equation*}
	where $\omega_0\in\mathfrak{S}_m$ is the longest element. But then the result follows because $\mathfrak{S}_\mathbf{d} = \mathfrak{S}_{m_1} \times \ldots \times \mathfrak{S}_{m_k}$ is a direct product.
	
	Now if $\theta$ is well-separated and $\lambda\in\mathscr{K}_n^\ell$ then, by \autoref{C: Kleshchev t_d orbit is full} and \autoref{L: linear independence},
	\begin{equation*}
	\prod_{j=1}^{k} [m_j]_t^! = \dim_t(S_\theta^\lambda \cdot e_\omega(\residue(\mathfrak{t}_{\mathbf{d}}^\diamond))) \geq \dim_t(D_\theta^\lambda \cdot e_\omega(\residue(\mathfrak{t}_{\mathbf{d}}^\diamond)))
	\end{equation*}
	and the result follows. 
\end{proof}

\emph{Note} that we did not make a specific choice for the field $K$ and its characteristic. But the graded dimension of the simple module does depend on the characteristic. Hence, \autoref{T: Lower bound} implies that we are looking at a specific piece of the simple module which is independent of the characteristic of the field $K$.

\section{Inner products}\label{S: Inner products}

This section gives the proof of \autoref{T: Inner Product} and \autoref{P: Diagonal Gram matrix} using Webster's diagram calculus. This completes the proof of \autoref{T: Main Theorem}, which gives a classification of the simple $\mathscr{R}_n^\Lambda(\theta)$-modules.

\subsection{Diagram Calculus}

In this section we prove some technical results that hold locally for the diagrammatic Cherednik algebra.

Partly following Bowman \cite[Section 3]{Bowman:ManyCellular}, we first introduce some useful notation. Recall the equivalence relation on boxes with its equivalence classes called diagonals. Let $\Delta = \{\delta_0,\ldots,\delta_k\}$ be an $i$-diagonal in some multipartition $\lambda$ such that $\delta_j = (r_0 + j, c_0 + j,l)$ for some $r_0,c_0\in\mathbb{N}$ and $1\leq l\leq\ell$. The diagonal \emph{starts} at $\delta_0 = (r_0,c_0,l)$ and \emph{ends} at $\delta_k$. An $i$-diagonal is \emph{addable} if $(r_0+(k+1),c_0+(k+1),l)$ is an addable box, \emph{removable} if $\delta_k$ is a removable box, and \emph{negligible} otherwise. There are two types of negligible $i$-diagonals. A negligible $i$-diagonal is \emph{$(i+1)$-negligible} if $(r_0+k,c_0+(k+1),l)\in\lambda$ and \emph{$(i-1)$-negligible} if $(r_0+k,c_0+(k+1),l)\notin\lambda$. An $i$-diagonal is a \emph{left diagonal} if $x_\delta < \theta_l$, for all $\delta\in\Delta$ . Similarly, it is a \emph{centred diagonal} if $x_\delta - \ell < \theta_l < x_\delta$ and a \emph{right diagonal} if $\theta_l < x_\delta - \ell$. Consequently, there are $3$ kinds of addable and removable $i$-diagonals and $6$ kinds of negligible $i$-diagonals.

Let $\lambda\in\mathscr{P}_n^\ell(\theta)$ and $\Delta = \{\delta_0,\ldots,\delta_k\}$ an $i$-diagonal in $\lambda$, with $\delta_j = (r_0 + j, c_0 + j,l)$ for some $r_0,c_0\in\mathbb{N}$ and $1\leq l\leq\ell$. Let $\Delta_L = \{(r,c,l)\in\lambda | (r,c,l) \sim (r_0,c_0+1,l)\}$ and $\Delta_R = \{(r,c,l)\in\lambda | (r,c,l) \sim (r_0+1,c_0,l)\}$. Define the \emph{extended $i$-diagonal of $\Delta$}, as $\overline{\Delta} = \Delta_L \cup \Delta \cup \Delta_R$.

\begin{Example}\label{E: Diagonals}
	Let $e=3$ and $\kappa=(0)$. Then for $\lambda = (9,7,6,6,5,3,2)$ the $0$-diagonals are:

	\[
	\begin{tikzpicture}[scale = 1]
	\draw (0,0) pic{russiantableau={scale=0.5}{{\color{red}{0},\color{red}{1},\phantom{1},0,\phantom{1},\phantom{1},0,\phantom{1},\phantom{1}},{\color{red}{2},\color{red}{0},\color{red}{1},\phantom{1},0,\phantom{1},\phantom{1}},{\phantom{1},\color{red}{2},\color{red}{0},\color{red}{1},\phantom{1},0},{0,\phantom{1},\color{red}{2},\color{red}{0},\color{red}{1},\phantom{1}},{\phantom{1},0,\phantom{1},\color{red}{2},\color{red}{0}},{\phantom{1},\phantom{1},0},{0,\phantom{1}}}};
	\end{tikzpicture}
	\]

	\noindent From left to right we have an addable left diagonal, a $2$-negligible left diagonal, a removable centred diagonal, a removable right diagonal, and a $1$-negligible right diagonal. The boxes marked with coloured numbers is the extended $0$-diagonal of the removable centred diagonal.
\end{Example}

\noindent \emph{Note} that if $\Delta = \{\delta_0,\ldots,\delta_k\}$ is a removable diagonal in $\lambda$ then $\Delta\setminus\delta_k$ is an addable diagonal in $\lambda\setminus\delta_k$.

\begin{Proposition}\label{P: Crossing to Dot}
	Let $i\in I$ then locally the following relation holds:
	\[

	\]
\end{Proposition}

\begin{proof}
	Applying relations \ref{D: Cherednik Algebra}\ref{D: Cherednik Dots and crossings} and then \ref{D: Cherednik Algebra}\ref{D: Cherednik double crossings}\ref{D: Cherednik double i crossing} gives us
	\[
	%
	\]
\end{proof}

\begin{Lemma}\label{L: Loop over wall}
	Let $i\in I$ then the following three local relations hold:
	\[\begin{aligned}
	%
	\end{aligned}\]
\end{Lemma}

\begin{proof}
	For $(a)$ we apply \autoref{P: Crossing to Dot} and then relation \ref{D: Cherednik Algebra}\ref{D: Cherednik double crossings}\ref{D: Cherednik red solid double}:
	\[
	\scriptsize
	%
	\]
	For $(b)$ we apply again \autoref{P: Crossing to Dot} and then relations \ref{D: Cherednik Algebra}\ref{D: Cherednik double crossings}\ref{D: Cherednik double i crossing} and \ref{D: Cherednik Algebra}\ref{D: Cherednik double crossings}\ref{D: Cherednik ghost solid double}:
	\[\begin{aligned}
	\scriptsize
	%
	\end{aligned}\]
	The proof for $(c)$ is similar to $(b)$ except that we use relation \ref{D: Cherednik Algebra}\ref{D: Cherednik double crossings}\ref{D: Cherednik solid ghost double} instead of \ref{D: Cherednik Algebra}\ref{D: Cherednik double crossings}\ref{D: Cherednik ghost solid double}.
\end{proof}

\begin{Lemma}\label{L: Decomposition}
	For $i\in I$ we have the following decomposition:
	\[\begin{aligned}
	%
	\end{aligned}\]
\end{Lemma}

\begin{proof}
	First, apply relation \ref{D: Cherednik Algebra}\ref{D: Cherednik triple crossings}\ref{D: Cherednik 2-solid triple} to get:
	\[\scriptsize
	%
	\]
	Next, apply relation \ref{D: Cherednik Algebra}\ref{D: Cherednik triple crossings}\ref{D: Cherednik 2-ghost triple} to the first summand giving
	\[\scriptsize
	%
	\]
	Now applying relations \ref{D: Cherednik Algebra}\ref{D: Cherednik double crossings}\ref{D: Cherednik solid ghost double} and \ref{D: Cherednik Algebra}\ref{D: Cherednik double crossings}\ref{D: Cherednik ghost solid double} to the second summands of the two equations above gives us the result.
\end{proof}

\begin{Lemma}\label{L: Crossing over B1 Brick}
	For $i\in I$ we have the following decomposition:
	\[\scriptsize
	%
	\]
\end{Lemma}

\begin{proof}
	First we apply \autoref{L: Decomposition} to get:
	\[\begin{aligned}
	%
	\end{aligned}\]
	By relation \ref{D: Cherednik Algebra}\ref{D: Cherednik double crossings}\ref{D: Cherednik double i crossing}, the second, fourth, and fifth summand are zero and by relation \ref{D: Cherednik Algebra}\ref{D: Cherednik triple crossings} the first summand becomes:
	\[
	%
	\]
	Finally, by relation \ref{D: Cherednik Algebra}\ref{D: Cherednik Dots and crossings}, the remaining summand, which is the third, becomes:
	\[
	%
	\]
	Collecting the surviving terms completes the proof.
\end{proof}

\noindent Now recall the definition of $e_\mu(\mathbf{i})$, for some $\mu\in\mathscr{C}_n^\ell(\theta)$ and $\mathbf{i}\in I^n$. Let $\Delta$ be an $i$-diagonal and $\overline{\Delta}$ its extended $i$-diagonal. For the following Lemmata it is useful to consider the diagram $e_{\overline{\Delta}}(\residue(\overline{\Delta}))$. To shorten the notation we just write $e_{\overline{\Delta}}$.

\begin{Example}
Continue \autoref{E: Diagonals} and let $\Delta$ be the centred removable $0$-diagonal. Then the diagram $e_{\overline{\Delta}}$ is:
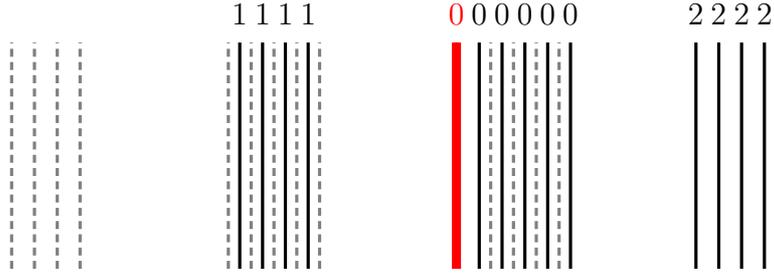
\[

\]
\end{Example}

\noindent \emph{Note} that in the following Lemma the straight vertical strings in the diagrams on the left-hand side of the equations correspond to diagrams $e_{\overline{\Delta}}$, where $\Delta$ is an addable $i$-diagonal.

\begin{Lemma}\label{L: Crossing over removable diagonal}
	Let $i\in I$ and $k\in\mathbb{N}$. Locally the following hold:
	\[\begin{aligned}
	%
	\end{aligned}\]
\end{Lemma}

\begin{proof}
	All three diagrams on the left hand side of the equations share the following subdiagram:
	\[
	%
	\end{aligned}\]
	We see that the second summand on the right-hand side of the equation contains the diagram on the left-hand side of the equation as a subdiagram, except that the number of vertical strings is decreased by $1$. Hence, we can apply \autoref{L: Crossing over B1 Brick} again on this subdiagram. As for the first summand, we can apply relations \ref{D: Cherednik Algebra}\ref{D: Cherednik triple crossings}\ref{D: Cherednik 2-ghost triple} and \ref{D: Cherednik 2-solid triple} to get
	\[
	%
	\]
	Note that relations \ref{D: Cherednik Algebra}\ref{D: Cherednik triple crossings}\ref{D: Cherednik 2-ghost triple} and \ref{D: Cherednik 2-solid triple} have a second summand with vertical strings. Replacing the $i$-crossing in the diagram with vertical strings yields a double $i$-crossing, hence, this diagram vanishes, by \ref{D: Cherednik Algebra}\ref{D: Cherednik double crossings}\ref{D: Cherednik double i crossing}.

	So consequently, we can repeat applying \autoref{L: Crossing over B1 Brick} until we have done so $k$-times. In total, we then collect the following terms:
	\[\begin{aligned}
	\begin{tikzpicture}[very thick,baseline,scale=0.2,rounded corners]\scriptsize
	\draw (-7.5,0) +(0,0) -- +(0,8);
	\draw (-3.5,0) +(0,0) -- +(0,8);

	\node at (-5.5,9) {$i+1$};
	\node at (-5.5,4) {$\cdots$};
	\draw [|-|,thin] (-3.5,-1) -- (-7.5,-1);
	\node at (-5.5,-2) {$k$};

	\draw (1,0) +(0,0) -- +(0,8);
	\draw[ghost] (-8,0) +(0,0) -- +(0,8);

	\node at (3.5,9) {$i$};
	\node at (3.5,4) {$\cdots$};
	\draw [|-|,thin] (6,-1) -- (1,-1);
	\node at (3.5,-2) {$k+1$};

	\draw (5,0) +(0,0) -- +(0,8);
	\draw[ghost] (-4,0) +(0,0) -- +(0,8);

	\draw (6,0) +(0,0) -- +(0,3) -- +(-0.5,5) -- +(-16,6);
	\draw (6,8) +(0,0) -- +(0,-3) -- +(-0.5,-5) -- +(-16,-6);
	\draw[ghost] (-3,0) +(0,0) -- +(0,3) -- +(-0.5,5) -- +(-7,5.5);
	\draw[ghost] (-3,8) +(0,0) -- +(0,-3) -- +(-0.5,-5) -- +(-7,-5.5);

	\draw[ghost] (1.5,0) +(0,0) -- +(0,8);
	\draw[ghost] (5.5,0) +(0,0) -- +(0,8);
	\draw (10.5,0) +(0,0) -- +(0,8);
	\draw (14.5,0) +(0,0) -- +(0,8);

	\node at (12.5,9) {$i-1$};
	\node at (12.5,4) {$\cdots$};
	\draw [|-|,thin] (14.5,-1) -- (10.5,-1);
	\node at (12.5,-2) {$k$};
	\end{tikzpicture}
	&\begin{tikzpicture}[very thick,baseline,scale=0.2,rounded corners]\scriptsize
	\node at (-20,4) {$=$};

	\normalsize
	\node at (-15,4) {$(-1)^k\cdot$};
	\scriptsize

	\draw (-7.5,0) +(0,0) -- +(0,8);
	\draw (-3.5,0) +(0,0) -- +(0,8);

	\node at (-5.5,9) {$i+1$};
	\node at (-5.5,4) {$\cdots$};
	\draw [|-|,thin] (-3.5,-1) -- (-7.5,-1);
	\node at (-5.5,-2) {$k$};

	\draw (2,0) +(0,0) -- +(0,8);
	\draw[ghost] (-7,0) +(0,0) -- +(0,8);

	\node at (3.5,9) {$i$};
	\node at (3.5,4) {$\cdots$};
	\draw [|-|,thin] (6,-1) -- (1,-1);
	\node at (3.5,-2) {$k+1$};

	\draw (5,0) +(0,0) -- +(0,8);
	\draw[ghost] (-4,0) +(0,0) -- +(0,8);

	\draw (6,0) +(0,0) -- +(0,8);
	\draw[ghost] (-3,0) +(0,0) -- +(0,8);

	\draw (1,0) +(0,0) -- +(0,3) -- +(-0.5,5) -- +(-11.5,6);
	\draw (1,8) +(0,0) -- +(0,-3) -- +(-0.5,-5) -- +(-11.5,-6);
	\draw[ghost] (-8,0) +(0,0) -- +(0,3) -- +(-0.5,5) -- +(-2.5,5.5);
	\draw[ghost] (-8,8) +(0,0) -- +(0,-3) -- +(-0.5,-5) -- +(-2.5,-5.5);

	\draw[ghost] (1.5,0) +(0,0) -- +(0,8);
	\draw[ghost] (5.5,0) +(0,0) -- +(0,8);
	\draw (10.5,0) +(0,0) -- +(0,8);
	\draw (14.5,0) +(0,0) -- +(0,8);

	\node at (12.5,9) {$i-1$};
	\node at (12.5,4) {$\cdots$};
	\draw [|-|,thin] (14.5,-1) -- (10.5,-1);
	\node at (12.5,-2) {$k$};
	\end{tikzpicture}\\
	&\begin{tikzpicture}[very thick,baseline,scale=0.2,rounded corners]\scriptsize
	\node at (-20,4) {$+$};

	\normalsize
	\node at (-15,4) {$\sum\limits_{j=0}^{k-1} (-1)^j\cdot$};
	\scriptsize

	\draw (-7.5,0) +(0,0) -- +(0,8);
	\draw (-3.5,0) +(0,0) -- +(0,8);

	\node at (-5.5,9) {$i+1$};
	\node at (-5.5,4) {$\cdots$};
	\draw [|-|,thin] (-3.5,-1) -- (-7.5,-1);
	\node at (-5.5,-2) {$k$};

	\draw (1,0) +(0,0) -- +(0,8);
	\draw[ghost] (-8,0) +(0,0) -- +(0,8);

	\node at (3.5,9) {$i$};
	\node at (3.5,4) {$\cdots$};
	\draw[|-|,thin] (6,-1) -- (4,-1);
	\node at (5,-2) {$j$};

	\draw (5,0) +(0,0) -- +(0,8);
	\draw[ghost] (-4,0) +(0,0) -- +(0,8);
	\draw (6,0) +(0,0) -- +(0,8);
	\draw[ghost] (-3,0) +(0,0) -- +(0,8);

	\draw (3.5,0) +(0,0) -- +(0,2) -- +(-5,5) -- +(-14,6);
	\draw[ghost] (-5.5,0) +(0,0) -- +(0,2) -- +(-5,5);
	\draw (3.5,8) +(0,0) -- +(0,-2) -- +(-5,-5) -- +(-14,-6);
	\draw[ghost] (-5.5,8) +(0,0) -- +(0,-2) -- +(-5,-5);

	\draw[ghost] (1.5,0) +(0,0) -- +(0,8);
	\draw[ghost] (5.5,0) +(0,0) -- +(0,8);
	\draw (10.5,0) +(0,0) -- +(0,8);
	\draw (14.5,0) +(0,0) -- +(0,8);

	\node at (12.5,9) {$i-1$};
	\node at (12.5,4) {$\cdots$};
	\draw [|-|,thin] (14.5,-1) -- (10.5,-1);
	\node at (12.5,-2) {$k$};
	\end{tikzpicture}
	\end{aligned}\]
	The difference for the three equations of the Lemma lies in either an extra solid $(i+1)$-, ghost $(i-1)$-, or red $i$-string to the left of the $i$-crossings in the decomposition above. The following is the case of the additional red $i$-string:
	\[\begin{aligned}
	\begin{tikzpicture}[very thick,baseline,scale=0.2,rounded corners]\scriptsize
	\draw[wei] (0,0) +(0,0) -- +(0,8);

	\draw (-7.5,0) +(0,0) -- +(0,8);
	\draw (-3.5,0) +(0,0) -- +(0,8);

	\node at (-5.5,9) {$i+1$};
	\node at (-5.5,4) {$\cdots$};
	\draw [|-|,thin] (-3.5,-1) -- (-7.5,-1);
	\node at (-5.5,-2) {$k$};

	\draw (1,0) +(0,0) -- +(0,8);
	\draw[ghost] (-8,0) +(0,0) -- +(0,8);

	\node at (3.5,9) {$i$};
	\node at (3.5,4) {$\cdots$};
	\draw [|-|,thin] (6,-1) -- (1,-1);
	\node at (3.5,-2) {$k+1$};

	\draw (5,0) +(0,0) -- +(0,8);
	\draw[ghost] (-4,0) +(0,0) -- +(0,8);

	\draw (6,0) +(0,0) -- +(0,3) -- +(-0.5,5) -- +(-16,6);
	\draw (6,8) +(0,0) -- +(0,-3) -- +(-0.5,-5) -- +(-16,-6);
	\draw[ghost] (-3,0) +(0,0) -- +(0,3) -- +(-0.5,5) -- +(-7,5.5);
	\draw[ghost] (-3,8) +(0,0) -- +(0,-3) -- +(-0.5,-5) -- +(-7,-5.5);

	\draw[ghost] (1.5,0) +(0,0) -- +(0,8);
	\draw[ghost] (5.5,0) +(0,0) -- +(0,8);
	\draw (10.5,0) +(0,0) -- +(0,8);
	\draw (14.5,0) +(0,0) -- +(0,8);

	\node at (12.5,9) {$i-1$};
	\node at (12.5,4) {$\cdots$};
	\draw [|-|,thin] (14.5,-1) -- (10.5,-1);
	\node at (12.5,-2) {$k$};
	\end{tikzpicture}
	&\begin{tikzpicture}[very thick,baseline,scale=0.2,rounded corners]\scriptsize
	\node at (-20,4) {$=$};

	\draw[wei] (-.5,0) +(0,0) -- +(0,8);

	\normalsize
	\node at (-15,4) {$(-1)^k\cdot$};
	\scriptsize

	\draw (-7.5,0) +(0,0) -- +(0,8);
	\draw (-3.5,0) +(0,0) -- +(0,8);

	\node at (-5.5,9) {$i+1$};
	\node at (-5.5,4) {$\cdots$};
	\draw [|-|,thin] (-3.5,-1) -- (-7.5,-1);
	\node at (-5.5,-2) {$k$};

	\draw (2,0) +(0,0) -- +(0,8);
	\draw[ghost] (-7,0) +(0,0) -- +(0,8);

	\node at (3.5,9) {$i$};
	\node at (3.5,4) {$\cdots$};
	\draw [|-|,thin] (6,-1) -- (1,-1);
	\node at (3.5,-2) {$k+1$};

	\draw (5,0) +(0,0) -- +(0,8);
	\draw[ghost] (-4,0) +(0,0) -- +(0,8);

	\draw (6,0) +(0,0) -- +(0,8);
	\draw[ghost] (-3,0) +(0,0) -- +(0,8);

	\draw (1,0) +(0,0) -- +(0,3) -- +(-0.5,5) -- +(-11.5,6);
	\draw (1,8) +(0,0) -- +(0,-3) -- +(-0.5,-5) -- +(-11.5,-6);
	\draw[ghost] (-8,0) +(0,0) -- +(0,3) -- +(-0.5,5) -- +(-2.5,5.5);
	\draw[ghost] (-8,8) +(0,0) -- +(0,-3) -- +(-0.5,-5) -- +(-2.5,-5.5);

	\draw[ghost] (1.5,0) +(0,0) -- +(0,8);
	\draw[ghost] (5.5,0) +(0,0) -- +(0,8);
	\draw (10.5,0) +(0,0) -- +(0,8);
	\draw (14.5,0) +(0,0) -- +(0,8);

	\node at (12.5,9) {$i-1$};
	\node at (12.5,4) {$\cdots$};
	\draw [|-|,thin] (14.5,-1) -- (10.5,-1);
	\node at (12.5,-2) {$k$};
	\end{tikzpicture}\\
	&\begin{tikzpicture}[very thick,baseline,scale=0.2,rounded corners]\scriptsize
	\node at (-20,4) {$+$};

	\draw[wei] (-.5,0) +(0,0) -- +(0,8);

	\normalsize
	\node at (-15,4) {$\sum\limits_{j=0}^{k-1} (-1)^j\cdot$};
	\scriptsize

	\draw (-7.5,0) +(0,0) -- +(0,8);
	\draw (-3.5,0) +(0,0) -- +(0,8);

	\node at (-5.5,9) {$i+1$};
	\node at (-5.5,4) {$\cdots$};
	\draw [|-|,thin] (-3.5,-1) -- (-7.5,-1);
	\node at (-5.5,-2) {$k$};

	\draw (1,0) +(0,0) -- +(0,8);
	\draw[ghost] (-8,0) +(0,0) -- +(0,8);

	\node at (3.5,9) {$i$};
	\node at (3.5,4) {$\cdots$};
	\draw[|-|,thin] (6,-1) -- (4,-1);
	\node at (5,-2) {$j$};

	\draw (5,0) +(0,0) -- +(0,8);
	\draw[ghost] (-4,0) +(0,0) -- +(0,8);
	\draw (6,0) +(0,0) -- +(0,8);
	\draw[ghost] (-3,0) +(0,0) -- +(0,8);

	\draw (3.5,0) +(0,0) -- +(0,2) -- +(-5,5) -- +(-14,6);
	\draw[ghost] (-5.5,0) +(0,0) -- +(0,2) -- +(-5,5);
	\draw (3.5,8) +(0,0) -- +(0,-2) -- +(-5,-5) -- +(-14,-6);
	\draw[ghost] (-5.5,8) +(0,0) -- +(0,-2) -- +(-5,-5);

	\draw[ghost] (1.5,0) +(0,0) -- +(0,8);
	\draw[ghost] (5.5,0) +(0,0) -- +(0,8);
	\draw (10.5,0) +(0,0) -- +(0,8);
	\draw (14.5,0) +(0,0) -- +(0,8);

	\node at (12.5,9) {$i-1$};
	\node at (12.5,4) {$\cdots$};
	\draw [|-|,thin] (14.5,-1) -- (10.5,-1);
	\node at (12.5,-2) {$k$};
	\end{tikzpicture}
	\end{aligned}\]
	Applying relation \ref{D: Cherednik Algebra}\ref{D: Cherednik triple crossings}\ref{D: Cherednik red triple} to all diagrams in the sum above yields
	\[\begin{aligned}
	\begin{tikzpicture}[very thick,baseline,scale=0.2,rounded corners]\scriptsize
	\draw[wei] (0,0) +(0,0) -- +(0,8);

	\draw (-7.5,0) +(0,0) -- +(0,8);
	\draw (-3.5,0) +(0,0) -- +(0,8);

	\node at (-5.5,9) {$i+1$};
	\node at (-5.5,4) {$\cdots$};
	\draw [|-|,thin] (-3.5,-1) -- (-7.5,-1);
	\node at (-5.5,-2) {$k$};

	\draw (2,0) +(0,0) -- +(0,8);
	\draw[ghost] (-7,0) +(0,0) -- +(0,8);

	\node at (3.5,9) {$i$};
	\node at (3.5,4) {$\cdots$};
	\draw [|-|,thin] (6,-1) -- (1,-1);
	\node at (3.5,-2) {$k+1$};

	\draw (5,0) +(0,0) -- +(0,8);
	\draw[ghost] (-4,0) +(0,0) -- +(0,8);

	\draw (6,0) +(0,0) -- +(0,8);
	\draw[ghost] (-3,0) +(0,0) -- +(0,8);

	\draw (1,0) +(0,0) -- +(0,3) -- +(-0.5,5) -- +(-11.5,6);
	\draw (1,8) +(0,0) -- +(0,-3) -- +(-0.5,-5) -- +(-11.5,-6);
	\draw[ghost] (-8,0) +(0,0) -- +(0,3) -- +(-0.5,5) -- +(-2.5,5.5);
	\draw[ghost] (-8,8) +(0,0) -- +(0,-3) -- +(-0.5,-5) -- +(-2.5,-5.5);

	\draw[ghost] (1.5,0) +(0,0) -- +(0,8);
	\draw[ghost] (5.5,0) +(0,0) -- +(0,8);
	\draw (10.5,0) +(0,0) -- +(0,8);
	\draw (14.5,0) +(0,0) -- +(0,8);

	\node at (12.5,9) {$i-1$};
	\node at (12.5,4) {$\cdots$};
	\draw [|-|,thin] (14.5,-1) -- (10.5,-1);
	\node at (12.5,-2) {$k$};
	\end{tikzpicture}
	&\begin{tikzpicture}[very thick,baseline,scale=0.2,rounded corners]\scriptsize
	\node at (-14,4) {$=$};

	\draw[wei] (0,0) +(0,0) -- +(0,8);

	\draw (-7.5,0) +(0,0) -- +(0,8);
	\draw (-3.5,0) +(0,0) -- +(0,8);

	\node at (-5.5,9) {$i+1$};
	\node at (-5.5,4) {$\cdots$};
	\draw [|-|,thin] (-3.5,-1) -- (-7.5,-1);
	\node at (-5.5,-2) {$k$};

	\draw (2,0) +(0,0) -- +(0,8);
	\draw[ghost] (-7,0) +(0,0) -- +(0,8);

	\node at (3.5,9) {$i$};
	\node at (3.5,4) {$\cdots$};
	\draw [|-|,thin] (6,-1) -- (1,-1);
	\node at (3.5,-2) {$k+1$};

	\draw (5,0) +(0,0) -- +(0,8);
	\draw[ghost] (-4,0) +(0,0) -- +(0,8);

	\draw (6,0) +(0,0) -- +(0,8);
	\draw[ghost] (-3,0) +(0,0) -- +(0,8);

	\draw (1,0) +(0,0) -- +(0,2) -- +(-3,5) -- +(-11.5,6);
	\draw (1,8) +(0,0) -- +(0,-2) -- +(-3,-5) -- +(-11.5,-6);
	\draw[ghost] (-8,0) +(0,0) -- +(0,3) -- +(-0.5,5) -- +(-2.5,5.5);
	\draw[ghost] (-8,8) +(0,0) -- +(0,-3) -- +(-0.5,-5) -- +(-2.5,-5.5);

	\draw[ghost] (1.5,0) +(0,0) -- +(0,8);
	\draw[ghost] (5.5,0) +(0,0) -- +(0,8);
	\draw (10.5,0) +(0,0) -- +(0,8);
	\draw (14.5,0) +(0,0) -- +(0,8);

	\node at (12.5,9) {$i-1$};
	\node at (12.5,4) {$\cdots$};
	\draw [|-|,thin] (14.5,-1) -- (10.5,-1);
	\node at (12.5,-2) {$k$};
	\end{tikzpicture}\\
	&\begin{tikzpicture}[very thick,baseline,scale=0.2,rounded corners]\scriptsize
	\node at (-14,4) {$+$};

	\node[red] at (0,9) {$i$};
	\draw[wei] (0,0) +(0,0) -- +(0,8);

	\draw (-7.5,0) +(0,0) -- +(0,8);
	\draw (-3.5,0) +(0,0) -- +(0,8);

	\node at (-5.5,9) {$i+1$};
	\node at (-5.5,4) {$\cdots$};
	\draw [|-|,thin] (-3.5,-1) -- (-7.5,-1);
	\node at (-5.5,-2) {$k$};

	\draw (1,0) +(0,0) -- +(0,8);
	\draw[ghost] (-8,0) +(0,0) -- +(0,8);

	\node at (3.5,9) {$i$};
	\node at (3.5,4) {$\cdots$};
	\draw [|-|,thin] (6,-1) -- (1,-1);
	\node at (3.5,-2) {$k+1$};

	\draw (5,0) +(0,0) -- +(0,8);
	\draw[ghost] (-4,0) +(0,0) -- +(0,8);

	\draw (6,0) +(0,0) -- +(0,8);
	\draw[ghost] (-3,0) +(0,0) -- +(0,8);

	\draw (-12,1) +(0,0) -- +(11,2) -- +(11,4) -- +(0,6);
	\draw[ghost] (-12,2) +(0,0) -- +(3,1) -- +(3,3) -- +(0,4);

	\draw[ghost] (1.5,0) +(0,0) -- +(0,8);
	\draw[ghost] (5.5,0) +(0,0) -- +(0,8);
	\draw (10.5,0) +(0,0) -- +(0,8);
	\draw (14.5,0) +(0,0) -- +(0,8);

	\node at (12.5,9) {$i-1$};
	\node at (12.5,4) {$\cdots$};
	\draw [|-|,thin] (14.5,-1) -- (10.5,-1);
	\node at (12.5,-2) {$k$};
	\end{tikzpicture}\\
	\begin{tikzpicture}[very thick,baseline,scale=0.2,rounded corners]\scriptsize
	\draw[wei] (-.5,0) +(0,0) -- +(0,8);

	\draw (-7.5,0) +(0,0) -- +(0,8);
	\draw (-3.5,0) +(0,0) -- +(0,8);

	\node at (-5.5,9) {$i+1$};
	\node at (-5.5,4) {$\cdots$};
	\draw [|-|,thin] (-3.5,-1) -- (-7.5,-1);
	\node at (-5.5,-2) {$k$};

	\draw (1,0) +(0,0) -- +(0,8);
	\draw[ghost] (-8,0) +(0,0) -- +(0,8);

	\node at (3.5,9) {$i$};
	\node at (3.5,4) {$\cdots$};
	\draw[|-|,thin] (6,-1) -- (4,-1);
	\node at (5,-2) {$j$};

	\draw (5,0) +(0,0) -- +(0,8);
	\draw[ghost] (-4,0) +(0,0) -- +(0,8);
	\draw (6,0) +(0,0) -- +(0,8);
	\draw[ghost] (-3,0) +(0,0) -- +(0,8);

	\draw (3.5,0) +(0,0) -- +(0,2) -- +(-5,5) -- +(-14,6);
	\draw[ghost] (-5.5,0) +(0,0) -- +(0,2) -- +(-5,5);
	\draw (3.5,8) +(0,0) -- +(0,-2) -- +(-5,-5) -- +(-14,-6);
	\draw[ghost] (-5.5,8) +(0,0) -- +(0,-2) -- +(-5,-5);

	\draw[ghost] (1.5,0) +(0,0) -- +(0,8);
	\draw[ghost] (5.5,0) +(0,0) -- +(0,8);
	\draw (10.5,0) +(0,0) -- +(0,8);
	\draw (14.5,0) +(0,0) -- +(0,8);

	\node at (12.5,9) {$i-1$};
	\node at (12.5,4) {$\cdots$};
	\draw [|-|,thin] (14.5,-1) -- (10.5,-1);
	\node at (12.5,-2) {$k$};
	\end{tikzpicture}
	&\begin{tikzpicture}[very thick,baseline,scale=0.2,rounded corners]\scriptsize
	\node at (-14,4) {$=$};

	\draw[wei] (0,0) +(0,0) -- +(0,8);

	\draw (-7.5,0) +(0,0) -- +(0,8);
	\draw (-3.5,0) +(0,0) -- +(0,8);

	\node at (-5.5,9) {$i+1$};
	\node at (-5.5,4) {$\cdots$};
	\draw [|-|,thin] (-3.5,-1) -- (-7.5,-1);
	\node at (-5.5,-2) {$k$};

	\draw (1,0) +(0,0) -- +(0,8);
	\draw[ghost] (-8,0) +(0,0) -- +(0,8);

	\node at (3.5,9) {$i$};
	\node at (3.5,4) {$\cdots$};
	\draw[|-|,thin] (6,-1) -- (4,-1);
	\node at (5,-2) {$j$};

	\draw (5,0) +(0,0) -- +(0,8);
	\draw[ghost] (-4,0) +(0,0) -- +(0,8);
	\draw (6,0) +(0,0) -- +(0,8);
	\draw[ghost] (-3,0) +(0,0) -- +(0,8);

	\draw (3.5,0) +(0,0) -- +(0,2) -- +(-7,5) -- +(-14,6);
	\draw[ghost] (-5.5,0) +(0,0) -- +(0,2) -- +(-5,5);
	\draw (3.5,8) +(0,0) -- +(0,-2) -- +(-7,-5) -- +(-14,-6);
	\draw[ghost] (-5.5,8) +(0,0) -- +(0,-2) -- +(-5,-5);

	\draw[ghost] (1.5,0) +(0,0) -- +(0,8);
	\draw[ghost] (5.5,0) +(0,0) -- +(0,8);
	\draw (10.5,0) +(0,0) -- +(0,8);
	\draw (14.5,0) +(0,0) -- +(0,8);

	\node at (12.5,9) {$i-1$};
	\node at (12.5,4) {$\cdots$};
	\draw [|-|,thin] (14.5,-1) -- (10.5,-1);
	\node at (12.5,-2) {$k$};
	\end{tikzpicture}\\
	&\begin{tikzpicture}[very thick,baseline,scale=0.2,rounded corners]\scriptsize
	\node at (-14,4) {$+$};

	\draw[wei] (-.5,0) +(0,0) -- +(0,8);

	\draw (-7.5,0) +(0,0) -- +(0,8);
	\draw (-3.5,0) +(0,0) -- +(0,8);

	\node at (-5.5,9) {$i+1$};
	\node at (-5.5,4) {$\cdots$};
	\draw [|-|,thin] (-3.5,-1) -- (-7.5,-1);
	\node at (-5.5,-2) {$k$};

	\draw (1,0) +(0,0) -- +(0,8);
	\draw[ghost] (-8,0) +(0,0) -- +(0,8);

	\node at (3.5,9) {$i$};
	\node at (3.5,4) {$\cdots$};
	\draw[|-|,thin] (6,-1) -- (4,-1);
	\node at (5,-2) {$j$};

	\draw (5,0) +(0,0) -- +(0,8);
	\draw[ghost] (-4,0) +(0,0) -- +(0,8);
	\draw (6,0) +(0,0) -- +(0,8);
	\draw[ghost] (-3,0) +(0,0) -- +(0,8);

	\draw (3.5,0) +(0,0) -- +(0,2) -- +(-3.25,3) -- +(-3.25,5) -- +(0,6) -- +(0,8);
	\draw[ghost] (-5.5,0) +(0,0) -- +(0,2) -- +(-3.25,3) -- +(-3.25,5) -- +(0,6) -- +(0,8);
	\draw (-12,1) +(0,0) -- +(10.5,2) -- +(10.5,4) -- +(0,6);
	\draw[ghost] (-12,2) +(0,0) -- +(2,1) -- +(2,3) -- +(0,4);

	\draw[ghost] (1.5,0) +(0,0) -- +(0,8);
	\draw[ghost] (5.5,0) +(0,0) -- +(0,8);
	\draw (10.5,0) +(0,0) -- +(0,8);
	\draw (14.5,0) +(0,0) -- +(0,8);

	\node at (12.5,9) {$i-1$};
	\node at (12.5,4) {$\cdots$};
	\draw [|-|,thin] (14.5,-1) -- (10.5,-1);
	\node at (12.5,-2) {$k$};
	\end{tikzpicture}
	\end{aligned}\]
	Where the last summand above vanishes by relation \ref{D: Cherednik Algebra}\ref{D: Cherednik double crossings}\ref{D: Cherednik double i crossing}. This proves the result.

	The other two cases are similar, except that \ref{D: Cherednik Algebra}\ref{D: Cherednik triple crossings}\ref{D: Cherednik red triple} is replaced by \ref{D: Cherednik Algebra}\ref{D: Cherednik triple crossings}\ref{D: Cherednik 2-ghost triple} or \ref{D: Cherednik 2-solid triple} respectively.
\end{proof}

\noindent \emph{Note} that in the following Lemma the straight vertical strings in the diagrams on the left- and right-hand side of the equations correspond to diagrams $e_{\overline{\Delta}}$, where $\Delta$ is a negligible $i$-diagonal.

\begin{Lemma}\label{L: Crossing over negligible diagonal}
	Let $i\in I$ and $k\in\mathbb{N}$. Locally, for a crossing of solid $i$-strings, the following local identities hold:
	
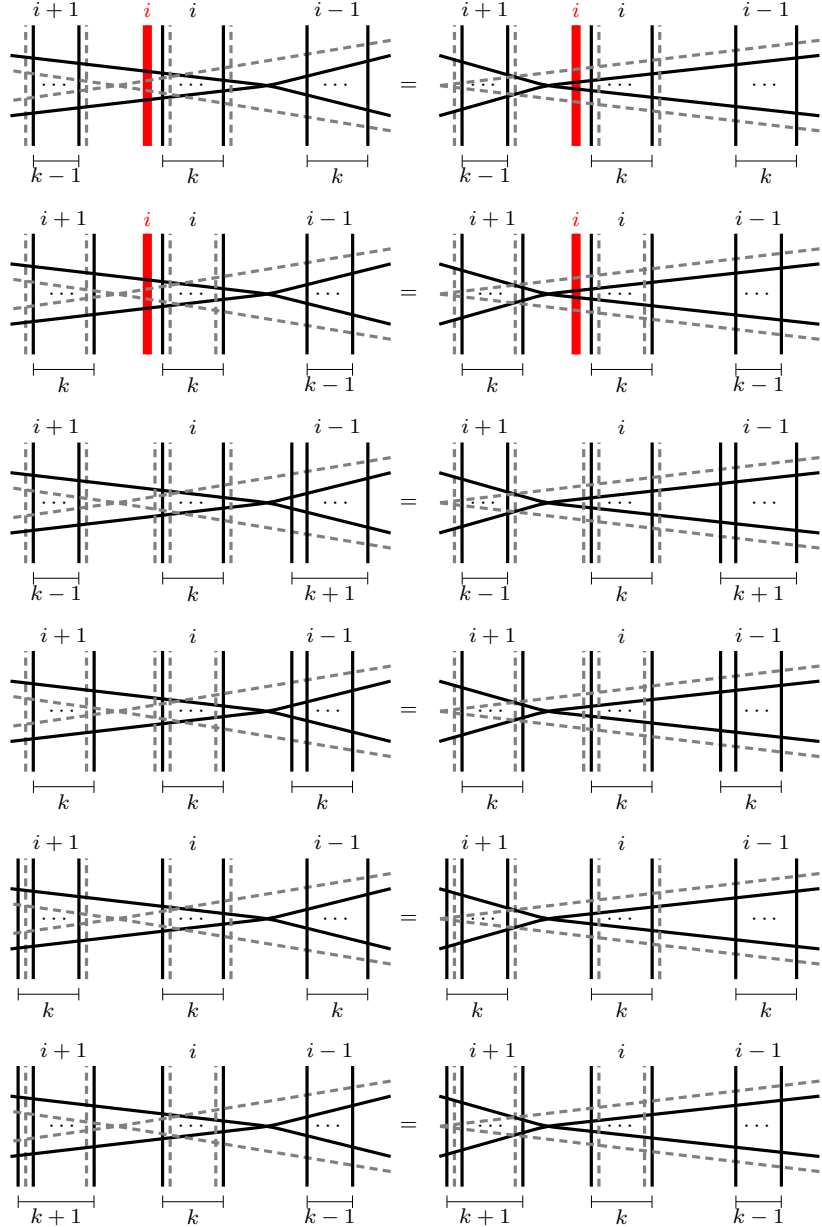
\[\begin{aligned}
	\begin{tikzpicture}[very thick,baseline,scale=0.2,rounded corners]\scriptsize
	\node[red] at (0,9) {$i$};
	\draw[wei] (0,0) +(0,0) -- +(0,8);

	\draw (-7.5,0) +(0,0) -- +(0,8);
	\draw (-4.5,0) +(0,0) -- +(0,8);

	\node at (-6,9) {$i+1$};
	\node at (-6,4) {$\cdots$};
	\draw [|-|,thin] (-4.5,-1) -- (-7.5,-1);
	\node at (-6,-2) {$k-1$};

	\draw (1,0) +(0,0) -- +(0,8);
	\draw[ghost] (-8,0) +(0,0) -- +(0,8);

	\node at (3,9) {$i$};
	\node at (3,4) {$\cdots$};
	\draw [|-|,thin] (5,-1) -- (1,-1);
	\node at (3,-2) {$k$};

	\draw (5,0) +(0,0) -- +(0,8);
	\draw[ghost] (-4,0) +(0,0) -- +(0,8);

	\draw (16,6) -- (8,4) -- (-9,2);
	\draw (16,2) -- (8,4) -- (-9,6);
	\draw[ghost] (16,7) -- (-2,4) -- (-9,3);
	\draw[ghost] (16,1) -- (-2,4) -- (-9,5);

	\draw[ghost] (1.5,0) +(0,0) -- +(0,8);
	\draw[ghost] (5.5,0) +(0,0) -- +(0,8);
	\draw (10.5,0) +(0,0) -- +(0,8);
	\draw (14.5,0) +(0,0) -- +(0,8);

	\node at (12.5,9) {$i-1$};
	\node at (12.5,4) {$\cdots$};
	\draw [|-|,thin] (14.5,-1) -- (10.5,-1);
	\node at (12.5,-2) {$k$};
	\end{tikzpicture}
	&\begin{tikzpicture}[very thick,baseline,scale=0.2,rounded corners]\scriptsize
	\node at (-11,4) {$=$};

	\node[red] at (0,9) {$i$};
	\draw[wei] (0,0) +(0,0) -- +(0,8);

	\draw (-7.5,0) +(0,0) -- +(0,8);
	\draw (-4.5,0) +(0,0) -- +(0,8);

	\node at (-6,9) {$i+1$};
	\node at (-6,4) {$\cdots$};
	\draw [|-|,thin] (-4.5,-1) -- (-7.5,-1);
	\node at (-6,-2) {$k-1$};

	\draw (1,0) +(0,0) -- +(0,8);
	\draw[ghost] (-8,0) +(0,0) -- +(0,8);

	\node at (3,9) {$i$};
	\node at (3,4) {$\cdots$};
	\draw [|-|,thin] (5,-1) -- (1,-1);
	\node at (3,-2) {$k$};

	\draw (5,0) +(0,0) -- +(0,8);
	\draw[ghost] (-4,0) +(0,0) -- +(0,8);

	\draw (16,6) -- (-2,4) -- (-9,2);
	\draw (16,2) -- (-2,4) -- (-9,6);
	\draw[ghost] (16,7) -- (-9,4);
	\draw[ghost] (16,1) -- (-9,4);

	\draw[ghost] (1.5,0) +(0,0) -- +(0,8);
	\draw[ghost] (5.5,0) +(0,0) -- +(0,8);
	\draw (10.5,0) +(0,0) -- +(0,8);
	\draw (14.5,0) +(0,0) -- +(0,8);

	\node at (12.5,9) {$i-1$};
	\node at (12.5,4) {$\cdots$};
	\draw [|-|,thin] (14.5,-1) -- (10.5,-1);
	\node at (12.5,-2) {$k$};
	\end{tikzpicture}\\
	\begin{tikzpicture}[very thick,baseline,scale=0.2,rounded corners]\scriptsize
	\node[red] at (0,9) {$i$};
	\draw[wei] (0,0) +(0,0) -- +(0,8);

	\draw (-7.5,0) +(0,0) -- +(0,8);
	\draw (-3.5,0) +(0,0) -- +(0,8);

	\node at (-5.5,9) {$i+1$};
	\node at (-5.5,4) {$\cdots$};
	\draw [|-|,thin] (-3.5,-1) -- (-7.5,-1);
	\node at (-5.5,-2) {$k$};

	\draw (1,0) +(0,0) -- +(0,8);
	\draw[ghost] (-8,0) +(0,0) -- +(0,8);

	\node at (3,9) {$i$};
	\node at (3,4) {$\cdots$};
	\draw [|-|,thin] (5,-1) -- (1,-1);
	\node at (3,-2) {$k$};

	\draw (5,0) +(0,0) -- +(0,8);
	\draw[ghost] (-4,0) +(0,0) -- +(0,8);

	\draw (16,6) -- (8,4) -- (-9,2);
	\draw (16,2) -- (8,4) -- (-9,6);
	\draw[ghost] (16,7) -- (-2,4) -- (-9,3);
	\draw[ghost] (16,1) -- (-2,4) -- (-9,5);

	\draw[ghost] (1.5,0) +(0,0) -- +(0,8);
	\draw[ghost] (4.5,0) +(0,0) -- +(0,8);
	\draw (10.5,0) +(0,0) -- +(0,8);
	\draw (13.5,0) +(0,0) -- +(0,8);

	\node at (12,9) {$i-1$};
	\node at (12,4) {$\cdots$};
	\draw [|-|,thin] (13.5,-1) -- (10.5,-1);
	\node at (12,-2) {$k-1$};
	\end{tikzpicture}
	&\begin{tikzpicture}[very thick,baseline,scale=0.2,rounded corners]\scriptsize
	\node at (-11,4) {$=$};

	\node[red] at (0,9) {$i$};
	\draw[wei] (0,0) +(0,0) -- +(0,8);

	\draw (-7.5,0) +(0,0) -- +(0,8);
	\draw (-3.5,0) +(0,0) -- +(0,8);

	\node at (-5.5,9) {$i+1$};
	\node at (-5.5,4) {$\cdots$};
	\draw [|-|,thin] (-3.5,-1) -- (-7.5,-1);
	\node at (-5.5,-2) {$k$};

	\draw (1,0) +(0,0) -- +(0,8);
	\draw[ghost] (-8,0) +(0,0) -- +(0,8);

	\node at (3,9) {$i$};
	\node at (3,4) {$\cdots$};
	\draw [|-|,thin] (5,-1) -- (1,-1);
	\node at (3,-2) {$k$};

	\draw (5,0) +(0,0) -- +(0,8);
	\draw[ghost] (-4,0) +(0,0) -- +(0,8);

	\draw (16,6) -- (-2,4) -- (-9,2);
	\draw (16,2) -- (-2,4) -- (-9,6);
	\draw[ghost] (16,7) -- (-9,4);
	\draw[ghost] (16,1) -- (-9,4);

	\draw[ghost] (1.5,0) +(0,0) -- +(0,8);
	\draw[ghost] (4.5,0) +(0,0) -- +(0,8);
	\draw (10.5,0) +(0,0) -- +(0,8);
	\draw (13.5,0) +(0,0) -- +(0,8);

	\node at (12,9) {$i-1$};
	\node at (12,4) {$\cdots$};
	\draw [|-|,thin] (13.5,-1) -- (10.5,-1);
	\node at (12,-2) {$k-1$};
	\end{tikzpicture}\\
	\begin{tikzpicture}[very thick,baseline,scale=0.2,rounded corners]\scriptsize
	\draw (-7.5,0) +(0,0) -- +(0,8);
	\draw (-4.5,0) +(0,0) -- +(0,8);

	\node at (-6,9) {$i+1$};
	\node at (-6,4) {$\cdots$};
	\draw [|-|,thin] (-4.5,-1) -- (-7.5,-1);
	\node at (-6,-2) {$k-1$};

	\draw (1,0) +(0,0) -- +(0,8);
	\draw[ghost] (-8,0) +(0,0) -- +(0,8);

	\node at (3,9) {$i$};
	\node at (3,4) {$\cdots$};
	\draw [|-|,thin] (5,-1) -- (1,-1);
	\node at (3,-2) {$k$};

	\draw (5,0) +(0,0) -- +(0,8);
	\draw[ghost] (-4,0) +(0,0) -- +(0,8);

	\draw (16,6) -- (8,4) -- (-9,2);
	\draw (16,2) -- (8,4) -- (-9,6);
	\draw[ghost] (16,7) -- (-2,4) -- (-9,3);
	\draw[ghost] (16,1) -- (-2,4) -- (-9,5);

	\draw[ghost] (0.5,0) +(0,0) -- +(0,8);
	\draw[ghost] (1.5,0) +(0,0) -- +(0,8);
	\draw[ghost] (5.5,0) +(0,0) -- +(0,8);
	\draw (9.5,0) +(0,0) -- +(0,8);
	\draw (10.5,0) +(0,0) -- +(0,8);
	\draw (14.5,0) +(0,0) -- +(0,8);

	\node at (12.5,9) {$i-1$};
	\node at (12.5,4) {$\cdots$};
	\draw [|-|,thin] (14.5,-1) -- (9.5,-1);
	\node at (12,-2) {$k+1$};
	\end{tikzpicture}
	&\begin{tikzpicture}[very thick,baseline,scale=0.2,rounded corners]\scriptsize
	\node at (-11,4) {$=$};

	\draw (-7.5,0) +(0,0) -- +(0,8);
	\draw (-4.5,0) +(0,0) -- +(0,8);

	\node at (-6,9) {$i+1$};
	\node at (-6,4) {$\cdots$};
	\draw [|-|,thin] (-4.5,-1) -- (-7.5,-1);
	\node at (-6,-2) {$k-1$};

	\draw (1,0) +(0,0) -- +(0,8);
	\draw[ghost] (-8,0) +(0,0) -- +(0,8);

	\node at (3,9) {$i$};
	\node at (3,4) {$\cdots$};
	\draw [|-|,thin] (5,-1) -- (1,-1);
	\node at (3,-2) {$k$};

	\draw (5,0) +(0,0) -- +(0,8);
	\draw[ghost] (-4,0) +(0,0) -- +(0,8);

	\draw (16,6) -- (-2,4) -- (-9,2);
	\draw (16,2) -- (-2,4) -- (-9,6);
	\draw[ghost] (16,7) -- (-9,4);
	\draw[ghost] (16,1) -- (-9,4);

	\draw[ghost] (0.5,0) +(0,0) -- +(0,8);
	\draw[ghost] (1.5,0) +(0,0) -- +(0,8);
	\draw[ghost] (5.5,0) +(0,0) -- +(0,8);
	\draw (9.5,0) +(0,0) -- +(0,8);
	\draw (10.5,0) +(0,0) -- +(0,8);
	\draw (14.5,0) +(0,0) -- +(0,8);

	\node at (12.5,9) {$i-1$};
	\node at (12.5,4) {$\cdots$};
	\draw [|-|,thin] (14.5,-1) -- (9.5,-1);
	\node at (12,-2) {$k+1$};
	\end{tikzpicture}\\
	\begin{tikzpicture}[very thick,baseline,scale=0.2,rounded corners]\scriptsize
	\draw (-7.5,0) +(0,0) -- +(0,8);
	\draw (-3.5,0) +(0,0) -- +(0,8);

	\node at (-5.5,9) {$i+1$};
	\node at (-5.5,4) {$\cdots$};
	\draw [|-|,thin] (-3.5,-1) -- (-7.5,-1);
	\node at (-5.5,-2) {$k$};

	\draw (1,0) +(0,0) -- +(0,8);
	\draw[ghost] (-8,0) +(0,0) -- +(0,8);

	\node at (3,9) {$i$};
	\node at (3,4) {$\cdots$};
	\draw [|-|,thin] (5,-1) -- (1,-1);
	\node at (3,-2) {$k$};

	\draw (5,0) +(0,0) -- +(0,8);
	\draw[ghost] (-4,0) +(0,0) -- +(0,8);

	\draw (16,6) -- (8,4) -- (-9,2);
	\draw (16,2) -- (8,4) -- (-9,6);
	\draw[ghost] (16,7) -- (-2,4) -- (-9,3);
	\draw[ghost] (16,1) -- (-2,4) -- (-9,5);

	\draw[ghost] (0.5,0) +(0,0) -- +(0,8);
	\draw[ghost] (1.5,0) +(0,0) -- +(0,8);
	\draw[ghost] (4.5,0) +(0,0) -- +(0,8);
	\draw (9.5,0) +(0,0) -- +(0,8);
	\draw (10.5,0) +(0,0) -- +(0,8);
	\draw (13.5,0) +(0,0) -- +(0,8);

	\node at (12,9) {$i-1$};
	\node at (12,4) {$\cdots$};
	\draw [|-|,thin] (13.5,-1) -- (9.5,-1);
	\node at (11.5,-2) {$k$};
	\end{tikzpicture}
	&\begin{tikzpicture}[very thick,baseline,scale=0.2,rounded corners]\scriptsize
	\node at (-11,4) {$=$};

	\draw (-7.5,0) +(0,0) -- +(0,8);
	\draw (-3.5,0) +(0,0) -- +(0,8);

	\node at (-5.5,9) {$i+1$};
	\node at (-5.5,4) {$\cdots$};
	\draw [|-|,thin] (-3.5,-1) -- (-7.5,-1);
	\node at (-5.5,-2) {$k$};

	\draw (1,0) +(0,0) -- +(0,8);
	\draw[ghost] (-8,0) +(0,0) -- +(0,8);

	\node at (3,9) {$i$};
	\node at (3,4) {$\cdots$};
	\draw [|-|,thin] (5,-1) -- (1,-1);
	\node at (3,-2) {$k$};

	\draw (5,0) +(0,0) -- +(0,8);
	\draw[ghost] (-4,0) +(0,0) -- +(0,8);

	\draw (16,6) -- (-2,4) -- (-9,2);
	\draw (16,2) -- (-2,4) -- (-9,6);
	\draw[ghost] (16,7) -- (-9,4);
	\draw[ghost] (16,1) -- (-9,4);

	\draw[ghost] (0.5,0) +(0,0) -- +(0,8);
	\draw[ghost] (1.5,0) +(0,0) -- +(0,8);
	\draw[ghost] (4.5,0) +(0,0) -- +(0,8);
	\draw (9.5,0) +(0,0) -- +(0,8);
	\draw (10.5,0) +(0,0) -- +(0,8);
	\draw (13.5,0) +(0,0) -- +(0,8);

	\node at (12,9) {$i-1$};
	\node at (12,4) {$\cdots$};
	\draw [|-|,thin] (13.5,-1) -- (9.5,-1);
	\node at (11.5,-2) {$k$};
	\end{tikzpicture}\\
	\begin{tikzpicture}[very thick,baseline,scale=0.2,rounded corners]\scriptsize
	\draw (-8.5,0) +(0,0) -- +(0,8);
	\draw (-7.5,0) +(0,0) -- +(0,8);
	\draw (-4.5,0) +(0,0) -- +(0,8);

	\node at (-6,9) {$i+1$};
	\node at (-6,4) {$\cdots$};
	\draw [|-|,thin] (-4.5,-1) -- (-8.5,-1);
	\node at (-6.5,-2) {$k$};

	\draw (1,0) +(0,0) -- +(0,8);
	\draw[ghost] (-8,0) +(0,0) -- +(0,8);

	\node at (3,9) {$i$};
	\node at (3,4) {$\cdots$};
	\draw [|-|,thin] (5,-1) -- (1,-1);
	\node at (3,-2) {$k$};

	\draw (5,0) +(0,0) -- +(0,8);
	\draw[ghost] (-4,0) +(0,0) -- +(0,8);

	\draw (16,6) -- (8,4) -- (-9,2);
	\draw (16,2) -- (8,4) -- (-9,6);
	\draw[ghost] (16,7) -- (-2,4) -- (-9,3);
	\draw[ghost] (16,1) -- (-2,4) -- (-9,5);

	\draw[ghost] (1.5,0) +(0,0) -- +(0,8);
	\draw[ghost] (5.5,0) +(0,0) -- +(0,8);
	\draw (10.5,0) +(0,0) -- +(0,8);
	\draw (14.5,0) +(0,0) -- +(0,8);

	\node at (12.5,9) {$i-1$};
	\node at (12.5,4) {$\cdots$};
	\draw [|-|,thin] (14.5,-1) -- (10.5,-1);
	\node at (12.5,-2) {$k$};
	\end{tikzpicture}
	&\begin{tikzpicture}[very thick,baseline,scale=0.2,rounded corners]\scriptsize
	\node at (-11,4) {$=$};

	\draw (-8.5,0) +(0,0) -- +(0,8);
	\draw (-7.5,0) +(0,0) -- +(0,8);
	\draw (-4.5,0) +(0,0) -- +(0,8);

	\node at (-6,9) {$i+1$};
	\node at (-6,4) {$\cdots$};
	\draw [|-|,thin] (-4.5,-1) -- (-8.5,-1);
	\node at (-6.5,-2) {$k$};

	\draw (1,0) +(0,0) -- +(0,8);
	\draw[ghost] (-8,0) +(0,0) -- +(0,8);

	\node at (3,9) {$i$};
	\node at (3,4) {$\cdots$};
	\draw [|-|,thin] (5,-1) -- (1,-1);
	\node at (3,-2) {$k$};

	\draw (5,0) +(0,0) -- +(0,8);
	\draw[ghost] (-4,0) +(0,0) -- +(0,8);

	\draw (16,6) -- (-2,4) -- (-9,2);
	\draw (16,2) -- (-2,4) -- (-9,6);
	\draw[ghost] (16,7) -- (-9,4);
	\draw[ghost] (16,1) -- (-9,4);

	\draw[ghost] (1.5,0) +(0,0) -- +(0,8);
	\draw[ghost] (5.5,0) +(0,0) -- +(0,8);
	\draw (10.5,0) +(0,0) -- +(0,8);
	\draw (14.5,0) +(0,0) -- +(0,8);

	\node at (12.5,9) {$i-1$};
	\node at (12.5,4) {$\cdots$};
	\draw [|-|,thin] (14.5,-1) -- (10.5,-1);
	\node at (12.5,-2) {$k$};
	\end{tikzpicture}\\
	\begin{tikzpicture}[very thick,baseline,scale=0.2,rounded corners]\scriptsize
	\draw (-8.5,0) +(0,0) -- +(0,8);
	\draw (-7.5,0) +(0,0) -- +(0,8);
	\draw (-3.5,0) +(0,0) -- +(0,8);

	\node at (-5.5,9) {$i+1$};
	\node at (-5.5,4) {$\cdots$};
	\draw [|-|,thin] (-3.5,-1) -- (-8.5,-1);
	\node at (-6,-2) {$k+1$};

	\draw (1,0) +(0,0) -- +(0,8);
	\draw[ghost] (-8,0) +(0,0) -- +(0,8);

	\node at (3,9) {$i$};
	\node at (3,4) {$\cdots$};
	\draw [|-|,thin] (5,-1) -- (1,-1);
	\node at (3,-2) {$k$};

	\draw (5,0) +(0,0) -- +(0,8);
	\draw[ghost] (-4,0) +(0,0) -- +(0,8);

	\draw (16,6) -- (8,4) -- (-9,2);
	\draw (16,2) -- (8,4) -- (-9,6);
	\draw[ghost] (16,7) -- (-2,4) -- (-9,3);
	\draw[ghost] (16,1) -- (-2,4) -- (-9,5);

	\draw[ghost] (1.5,0) +(0,0) -- +(0,8);
	\draw[ghost] (4.5,0) +(0,0) -- +(0,8);
	\draw (10.5,0) +(0,0) -- +(0,8);
	\draw (13.5,0) +(0,0) -- +(0,8);

	\node at (12,9) {$i-1$};
	\node at (12,4) {$\cdots$};
	\draw [|-|,thin] (13.5,-1) -- (10.5,-1);
	\node at (12,-2) {$k-1$};
	\end{tikzpicture}
	&\begin{tikzpicture}[very thick,baseline,scale=0.2,rounded corners]\scriptsize
	\node at (-11,4) {$=$};

	\draw (-8.5,0) +(0,0) -- +(0,8);
	\draw (-7.5,0) +(0,0) -- +(0,8);
	\draw (-3.5,0) +(0,0) -- +(0,8);

	\node at (-5.5,9) {$i+1$};
	\node at (-5.5,4) {$\cdots$};
	\draw [|-|,thin] (-3.5,-1) -- (-8.5,-1);
	\node at (-6,-2) {$k+1$};

	\draw (1,0) +(0,0) -- +(0,8);
	\draw[ghost] (-8,0) +(0,0) -- +(0,8);

	\node at (3,9) {$i$};
	\node at (3,4) {$\cdots$};
	\draw [|-|,thin] (5,-1) -- (1,-1);
	\node at (3,-2) {$k$};

	\draw (5,0) +(0,0) -- +(0,8);
	\draw[ghost] (-4,0) +(0,0) -- +(0,8);

	\draw (16,6) -- (-2,4) -- (-9,2);
	\draw (16,2) -- (-2,4) -- (-9,6);
	\draw[ghost] (16,7) -- (-9,4);
	\draw[ghost] (16,1) -- (-9,4);

	\draw[ghost] (1.5,0) +(0,0) -- +(0,8);
	\draw[ghost] (4.5,0) +(0,0) -- +(0,8);
	\draw (10.5,0) +(0,0) -- +(0,8);
	\draw (13.5,0) +(0,0) -- +(0,8);

	\node at (12,9) {$i-1$};
	\node at (12,4) {$\cdots$};
	\draw [|-|,thin] (13.5,-1) -- (10.5,-1);
	\node at (12,-2) {$k-1$};
	\end{tikzpicture}
	\end{aligned}\]
\end{Lemma}

\begin{proof}
	Each of these equations is proved using a similar argument to that of \autoref{L: Crossing over removable diagonal}. When the $i$-crossing is dragged leftward over the previous vertical solid $i+1$-string, ghost $i-1$-string, or red $i$-string, respectively, we apply relations \ref{D: Cherednik Algebra}\ref{D: Cherednik triple crossings}\ref{D: Cherednik 2-ghost triple}, \ref{D: Cherednik 2-solid triple}, or \ref{D: Cherednik red triple} respectively. Similarly, as argued in the proof of \autoref{L: Crossing over removable diagonal} by application of the relations there is one diagram with the crossing dragged over the vertical string and one diagram where the crossing is replaced by vertical $i$-strings locally. Because of the double $i$-crossings the later diagram vanishes, by \ref{D: Cherednik Algebra}\ref{D: Cherednik double crossings}\ref{D: Cherednik double i crossing}. Hence, in each case, an $i$-crossing can be dragged leftwards through all the vertical strings, as claimed.
\end{proof}

\noindent \emph{Note} that in the following Lemma the straight vertical strings in the diagrams on the left-hand side of the equations correspond to diagrams $e_{\overline{\Delta}}$, where $\Delta$ is a negligible $i$-diagonal.

\begin{Lemma}\label{L: Loop over negligible diagonal}
	Let $i\in I$ and $k\in\mathbb{N}$. Locally, for a curved solid $i$-string, the following hold:
	\[\begin{aligned}

	\end{aligned}\]
\end{Lemma}

\begin{proof}
	The proof for all six equations is similar, so we prove the first equation carefully and then sketch how to prove the other five. The diagram on the left-hand side of the first equation contains the following subdiagram
	\[
	%
	\]
	Applying \ref{D: Cherednik Algebra}\ref{D: Cherednik double crossings}\ref{D: Cherednik solid ghost double} to this subdiagram gives us
	\[
	%
	\]
	The first summand vanishes by relation \ref{D: Cherednik Algebra}\ref{D: Cherednik double crossings}\ref{D: Cherednik double i crossing}. Applying \ref{D: Cherednik Algebra}\ref{D: Cherednik Dots and crossings} to the second summand gives us then
	\[
	%
	\]
	Again here the first summand vanishes by relation \ref{D: Cherednik Algebra}\ref{D: Cherednik double crossings}\ref{D: Cherednik double i crossing}. So combining both equations gives us
	\[
	%
	\]
	Replacing the subdiagram in with the one on the right-hand side above in the original diagram gives us
	\[
	%
	\]
	Now we see that the diagram on the right-hand side above contains a subdiagram as in the first equation of \autoref{L: Crossing over removable diagonal}. Applying \autoref{L: Crossing over removable diagonal} on this subdiagram completes the proof for the first equation.

	The proof for the third and fifth equation follows the same argument. For the proof of the second, fourth, and sixth equation we instead look at the subdiagram
	\[
	%
	\]
	By a similar chain of arguments we get
	\[
	%
	\]
	The result then follows again by applying \autoref{L: Crossing over removable diagonal}.
\end{proof}

\subsection{Computing inner products}

This section proves the key Theorem for the classification of all simple $\mathscr{R}_n^\Lambda(\theta)$-modules. Before giving the proof we explain how to construct certain diagrams that we need to compute inner products. We give two different constructions of these diagrams and show that these different diagrams represent the same element in $\mathscr{R}_n^\Lambda(\theta)$.

Fix a small positive quantity $\delta < \frac{1}{2n}$. Let $\lambda\in\mathscr{U}_n^\ell(\theta)$ and $\mathbf{d} = \mathbf{s}_1\ldots \mathbf{s}_k$ a decomposition path for $\lambda$. Recall $\omega = (\emptyset,\ldots,\emptyset,(1^n)) \in\mathscr{P}_n^\ell(\theta)$ and let $\omega = \{\beta_1,\ldots,\beta_n\}$ considered as a set of boxes, such that $\beta_a \prec_\theta \beta_b$, for $1\leq a < b\leq n$. Let $\mathbf{s}_j = (s_1^{(j)},\ldots,s_{m_j}^{(j)})$, for $1\leq j\leq k$. Let $\sigma\in\mathfrak{S}_{\mathbf{d}}$ and consider the tableau $\mathfrak{t}_{\mathbf{d}}\cdot\sigma$. We now describe a construction of a diagram of the element $c_{\mathfrak{t}_{\mathbf{d}}\cdot\sigma}^\lambda$ satisfying the definition of $\theta$-diagrams; see beginning of \autoref{S: Cherednik algebra section}.

\emph{Construction A:} First we draw the red strings straight vertically then we need to draw the solid and their ghost strings. We do this recursively and assume we have already drawn the solid and ghost strings for $(\mathfrak{t}_{\mathbf{d}}\cdot\sigma)_{\downarrow r-1}$ already, for $1\leq r-1 < n$. Let $\alpha\in\lambda$ be the box containing the number $r$ in $\mathfrak{t}_{\mathbf{d}}\cdot\sigma$ then we need to draw a solid string of type $(x_\alpha,x_{\beta_r})$. This path is a union of straight line segments between the following four points:
\begin{equation*}
(x_\alpha,1) \longrightarrow (x_\alpha,\frac{r}{n+1}) \longrightarrow (x_{\beta_r},\frac{r}{n+1}-\delta) \longrightarrow (x_{\beta_r},0)
\end{equation*}
Together with its ghost string we get the following strings:
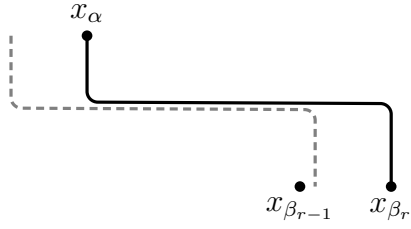
\[
\begin{tikzpicture}[very thick,baseline,rounded corners]
\draw[ghost,yscale=2] (-3,1)  +(0,0) -- +(0,-0.48) -- +(4,-0.49) -- +(4,-1);
\draw[solid,yscale=2] (-2,1)  +(0,0) -- +(0,-0.44) -- +(4,-0.45) -- +(4,-1);

\node at (-2,2.3) {$x_\alpha$};
\node at (2,-0.3) {$x_{\beta_r}$};
\node at (0.8,-0.3) {$x_{\beta_{r-1}}$};

\fill (-2,2) circle (2pt);
\fill (2,0) circle (2pt);
\fill (0.8,0) circle (2pt);
\end{tikzpicture}
\]
Define $A_{\mathfrak{t}_{\mathbf{d}}\cdot\sigma}^\lambda$ to be the resulting diagram for $c_{\mathfrak{t}_{\mathbf{d}}\cdot\sigma}^\lambda$ using Construction A. The benefit of this construction is that in the horizontal strip between $y_1 = \frac{r_1-1}{n+1}$ and $y_2 = \frac{r_2}{n+1}$, for $1 \leq r_1 < r_2 \leq n$, all solid strings of type $(x_\alpha,x_{\beta_j})$ and their ghost strings are straight vertical, for $1\leq j\leq n$ such that $j\notin [r_1,r_2]$. By construction of $\mathfrak{t}_{\mathbf{d}}\cdot\sigma$, this allows us to partition the diagram into $k$ horizontal strips according to the $k$ staggered sequences in $\mathbf{d}$.

\begin{Example}\label{E: Presentation 1}
	Let $\theta=(0,8)$ and $\lambda = ((2,1),(3,1))$. Let $e=3$ and $\kappa = (0,0)$. Then there is only one staggered sequence $\mathbf{d}$; see \autoref{E: decomposition tableaux}. The diagram $A_{\mathfrak{t}_{\mathbf{d}}}^\lambda$ of $c_{\mathfrak{t}_{\mathbf{d}}}^\lambda$ is:
	
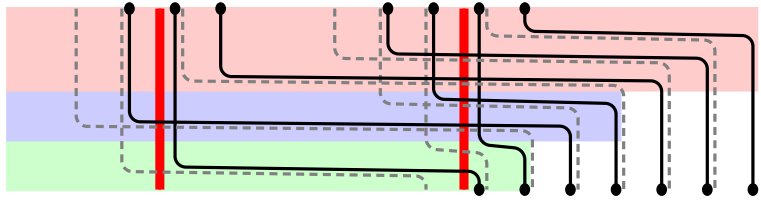
\[
	\begin{tikzpicture}[very thick,baseline,rounded corners,yscale = 1.2]
	\filldraw[color=green!20,sharp corners] (-2,0) +(0,0) rectangle +(6.85,0.55);
	\filldraw[color=blue!20,sharp corners] (-2,0.55) +(0,0) rectangle +(8.05,0.55);
	\filldraw[color=red!20,sharp corners] (-2,1.1) +(0,0) rectangle +(9.85,0.9);

	\draw[wei](0,0) -- (0,2);
	\draw[wei](4,0) -- (4,2);

	\fill (-0.4,2) circle (2pt);
	\fill (0.2,2) circle (2pt);
	\fill (0.8,2) circle (2pt);

	\fill (4.2,2) circle (2pt);
	\fill (3.6,2) circle (2pt);
	\fill (3,2) circle (2pt);
	\fill (4.8,2) circle (2pt);

	\fill (4.2,0) circle (2pt);
	\fill (4.8,0) circle (2pt);
	\fill (5.4,0) circle (2pt);
	\fill (6,0) circle (2pt);
	\fill (6.6,0) circle (2pt);
	\fill (7.2,0) circle (2pt);
	\fill (7.8,0) circle (2pt);

	\draw[ghost] (-0.5,2) +(0,0) -- +(0,-1.8) -- +(4,-1.85) -- +(4,-2);
	\draw[solid] (0.2,2) +(0,0) -- +(0,-1.75) -- +(4,-1.8) -- +(4,-2);

	\draw[ghost] (3.5,2) +(0,0) -- +(0,-1.55) -- +(0.8,-1.6) -- +(0.8,-2);
	\draw[solid] (4.2,2) +(0,0) -- +(0,-1.5) -- +(0.6,-1.55) -- +(0.6,-2);

	\draw[ghost] (-1.1,2) +(0,0) -- +(0,-1.3) -- +(6,-1.35) -- +(6,-2);
	\draw[solid] (-0.4,2) +(0,0) -- +(0,-1.25) -- +(5.8,-1.3) -- +(5.8,-2);

	\draw[ghost] (2.9,2) +(0,0) -- +(0,-1.05) -- +(2.6,-1.1) -- +(2.6,-2);
	\draw[solid] (3.6,2) +(0,0) -- +(0,-1) -- +(2.4,-1.05) -- +(2.4,-2);

	\draw[ghost] (0.3,2) +(0,0) -- +(0,-0.8) -- +(5.8,-0.85) -- +(5.8,-2);
	\draw[solid] (0.8,2) +(0,0) -- +(0,-0.75) -- +(5.8,-0.8) -- +(5.8,-2);

	\draw[ghost] (2.3,2) +(0,0) -- +(0,-0.55) -- +(4.4,-0.6) -- +(4.4,-2);
	\draw[solid] (3,2) +(0,0) -- +(0,-0.5) -- +(4.2,-0.55) -- +(4.2,-2);

	\draw[ghost] (4.3,2) +(0,0) -- +(0,-0.3) -- +(3,-0.35) -- +(3,-2);
	\draw[solid] (4.8,2) +(0,0) -- +(0,-0.25) -- +(3,-0.3) -- +(3,-2);
	\end{tikzpicture}
	\]
	The three different shadings highlight the three horizontal strips
\end{Example}

We can now prove the main result of this section. Recall from \autoref{D: Decomposition Tableaux} that we have a pair of decomposition tableaux $\mathfrak{t}_\mathbf{d}$ and $\mathfrak{t}_\mathbf{d}^\diamond$, for each decomposition path $\mathbf{d}$.

\begin{Theorem}\label{T: Chap 5 main}
	Let $\lambda\in\mathscr{U}^\ell(\theta)$ be an Uglov multipartition and $\mathbf{d}$ any decomposition path for $\lambda$. Then $|\langle A_{\mathfrak{t}_\mathbf{d}}^\lambda, A_{\mathfrak{t}_\mathbf{d}^\diamond}^\lambda\rangle_\lambda| = 1$.
\end{Theorem}

\begin{proof}
	Let $\lambda\in\mathscr{U}^\ell(\theta)$ and $\mathbf{d} = \mathbf{s}_1\ldots \mathbf{s}_k$ a decomposition path for $\lambda$. Recall \ref{eq: def inner product}:
	\begin{equation*}
	\langle c_{\mathfrak{t}_\mathbf{d}}^\lambda, c_{\mathfrak{t}_\mathbf{d}^\diamond}^\lambda\rangle_\lambda c_{\mathfrak{t}_{\mathbf{d}}}^\lambda = c_{\mathfrak{t}_{\mathbf{d}}}^\lambda c_{\mathfrak{t}_{\mathbf{d}}^\diamond \mathfrak{t}_{\mathbf{d}}}^\lambda = c_{\mathfrak{t}_{\mathbf{d}}}^\lambda \circ (c_{\mathfrak{t}_{\mathbf{d}}^\diamond}^\lambda)^\ast \circ c_{\mathfrak{t}_{\mathbf{d}}}^\lambda
	\end{equation*}
	Hence, the inner product is determined by the diagram $A_{\mathfrak{t}_{\mathbf{d}}}^\lambda \circ (A_{\mathfrak{t}_{\mathbf{d}}^\diamond}^\lambda)^\ast$. We show that this diagram `straightens out' to $\pm e_\lambda(\residue(\mathfrak{t}_{\mathbf{d}}))$ modulo terms of more dominant shape. This will show that the inner product is $\pm 1$, giving the result.

	We can describe the diagram $A_{\mathfrak{t}_{\mathbf{d}}}^\lambda \circ (A_{\mathfrak{t}_{\mathbf{d}}^\diamond}^\lambda)^\ast$ recursively as the diagram
	
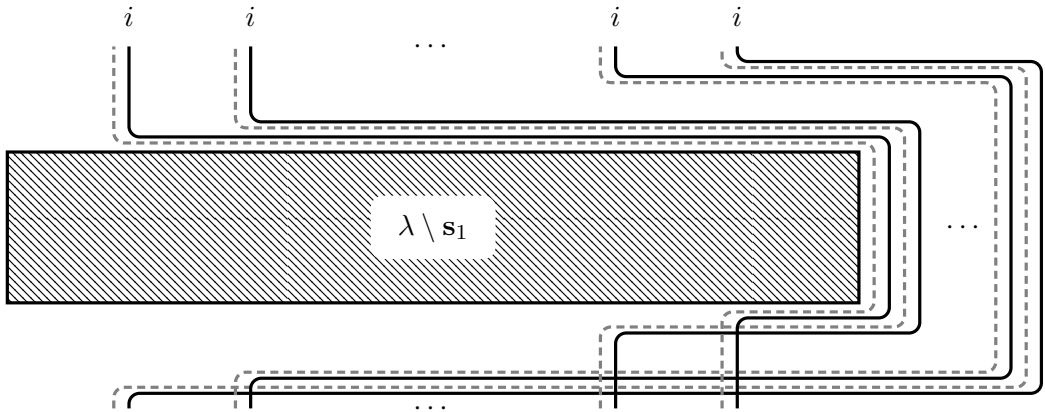
\[
	\begin{tikzpicture}[very thick,baseline,scale=0.4,rounded corners]\small
	\node at (0,13) {$i$};
	\node at (4,13) {$i$};
	\node at (16,13) {$i$};
	\node at (20,13) {$i$};

	\draw[ghost] (-0.5,0) +(0,0) -- +(0,0.7) -- +(30,0.7) -- +(30,11.3) -- +(20,11.3) -- +(20,12);
	\draw (0,0) +(0,0) -- +(0,0.5) -- +(30,0.5) -- +(30,11.5) -- +(20,11.5) -- +(20,12);

	\draw[ghost] (3.5,0) +(0,0) -- +(0,1.2) -- +(25,1.2) -- +(25,10.8) -- +(12,10.8) -- +(12,12);
	\draw (4,0) +(0,0) -- +(0,1) -- +(25,1) -- +(25,11) -- +(12,11) -- +(12,12);

	\draw[ghost] (15.5,0) +(0,0) -- +(0,2.7) -- +(10,2.7) -- +(10,9.3) -- +(-12,9.3) -- +(-12,12);
	\draw (16,0) +(0,0) -- +(0,2.5) -- +(10,2.5) -- +(10,9.5) -- +(-12,9.5) -- +(-12,12);

	\draw[ghost] (19.5,0) +(0,0) -- +(0,3.2) -- +(5,3.2) -- +(5,8.8) -- +(-20,8.8) -- +(-20,12);
	\draw (20,0) +(0,0) -- +(0,3) -- +(5,3) -- +(5,9) -- +(-20,9) -- +(-20,12);

	\node at (10,12) {$\cdots$};
	\node at (10,0) {$\cdots$};
	\node at (27.5,6) {$\cdots$};

	\draw[pattern = north west lines,sharp corners] (-4,3.5) +(0,0) rectangle +(28,5);
	\draw[fill,color = white] (8,5) +(0,0) rectangle +(4,2);
	\node at (10,6) {$\lambda\setminus \mathbf{s}_1$};
	\end{tikzpicture}
	\]
	where $\mu = \lambda\setminus \mathbf{s}_1$ is a smaller Uglov multipartition. The shaded part of the diagram above is the diagram $A_{\mathfrak{t}_{\mathbf{d}'}}^\mu \circ (A_{\mathfrak{t}_{\mathbf{d}'}^\diamond}^\mu)^\ast$, where $\mathbf{d}' = \mathbf{s}_2\ldots\mathbf{s}_k$. By construction, in the diagram above all strings not belonging to the staggered sequence $\mathbf{s}_1$ are straight vertical when they are above and below the shaded component. This recursive description of this diagram allows us to argue by induction on $|\mathbf{d}| = k$.

	Let $k=0$ then $\lambda = \emptyset$ and there is nothing to prove. So let $k>0$ and assume the statement is true for $k-1$. Hence, by induction, the diagram above simplifies to:
	
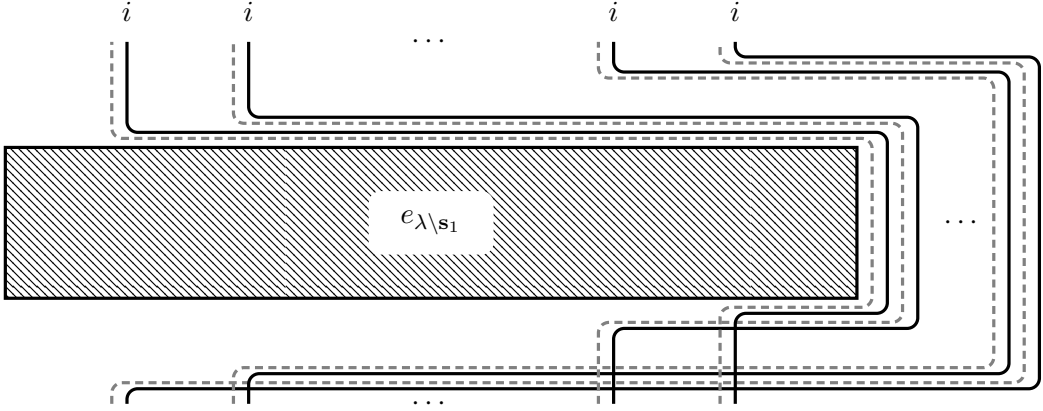
\[
	\begin{tikzpicture}[very thick,baseline,scale=0.4,rounded corners]\small
	\node at (0,13) {$i$};
	\node at (4,13) {$i$};
	\node at (16,13) {$i$};
	\node at (20,13) {$i$};

	\draw[ghost] (-0.5,0) +(0,0) -- +(0,0.7) -- +(30,0.7) -- +(30,11.3) -- +(20,11.3) -- +(20,12);
	\draw (0,0) +(0,0) -- +(0,0.5) -- +(30,0.5) -- +(30,11.5) -- +(20,11.5) -- +(20,12);

	\draw[ghost] (3.5,0) +(0,0) -- +(0,1.2) -- +(25,1.2) -- +(25,10.8) -- +(12,10.8) -- +(12,12);
	\draw (4,0) +(0,0) -- +(0,1) -- +(25,1) -- +(25,11) -- +(12,11) -- +(12,12);

	\draw[ghost] (15.5,0) +(0,0) -- +(0,2.7) -- +(10,2.7) -- +(10,9.3) -- +(-12,9.3) -- +(-12,12);
	\draw (16,0) +(0,0) -- +(0,2.5) -- +(10,2.5) -- +(10,9.5) -- +(-12,9.5) -- +(-12,12);

	\draw[ghost] (19.5,0) +(0,0) -- +(0,3.2) -- +(5,3.2) -- +(5,8.8) -- +(-20,8.8) -- +(-20,12);
	\draw (20,0) +(0,0) -- +(0,3) -- +(5,3) -- +(5,9) -- +(-20,9) -- +(-20,12);

	\node at (10,12) {$\cdots$};
	\node at (10,0) {$\cdots$};
	\node at (27.5,6) {$\cdots$};

	\draw[pattern = north west lines,sharp corners] (-4,3.5) +(0,0) rectangle +(28,5);
	\draw[fill,color = white] (8,5) +(0,0) rectangle +(4,2);
	\node at (10,6) {$e_{\lambda\setminus \mathbf{s}_1}$};
	\end{tikzpicture}
	\]
	where all of the strings inside the shaded area are now straight vertical. So it remains to show that the non-straight vertical strings, which correspond to $\mathbf{s}_1$, can be straightened. We assume that $\mathbf{s}_1 = (s_1,\ldots,s_m)$ is a sequence of residue $i\in I$ and that $\mathbf{s}_1$ has $|\mathbf{s}_1| = m$ solid strings.

	We now proceed by induction on $m$. Let $m = 1$ then the diagram is:
	
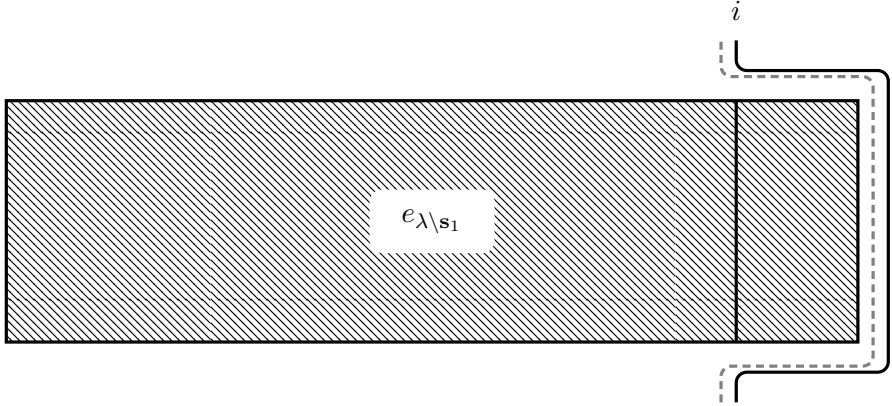
\[
	\begin{tikzpicture}[very thick,baseline,scale=0.4,rounded corners]\small
	\node at (20,13) {$i$};

	\draw[ghost] (19.5,0) +(0,0) -- +(0,1.2) -- +(5,1.2) -- +(5,10.8) -- +(0,10.8) -- +(0,12);
	\draw (20,0) +(0,0) -- +(0,1) -- +(5,1) -- +(5,11) -- +(0,11) -- +(0,12);

	\draw[pattern = north west lines,sharp corners] (-4,2) +(0,0) rectangle +(28,8);
	\draw[fill,color = white] (8,5) +(0,0) rectangle +(4,2);
	\node at (10,6) {$e_{\lambda\setminus \mathbf{s}_1}$};

	\draw (20,2) -- (20,10);
	\end{tikzpicture}
	\]
	where the black line is the straight vertical line at coordinate $x_{s_1}$. This is the place where the solid string becomes straight vertical. We refer to the curved $i$-string in the diagram above as an $i$-loop for our argument. Now because $\mathbf{s}_1$ is a staggered sequence we have $x_\alpha < x_{s_1}$, for every addable or removable $i$-box $\alpha\in\lambda$. Hence, to the right of this line there can only be straight vertical solid $i$-strings as in diagrams $e_{\overline{\Delta}}$, where $\Delta$ is a negligible $i$-diagonal. If there are no such strings then there are also no vertical solid $i+1$-strings and ghost $i-1$-strings as in diagrams $e_{\overline{\Delta}}$, where $\Delta$ is a negligible $i$-diagonal. So the $i$-loop can be straightened with no extra terms, by relation \ref{D: Cherednik Algebra}\ref{D: Cherednik double crossings}.

	If there are such strings, then the diagram does contain a subdiagram $e_{\overline{\Delta}}$, where $\Delta$ is a negligible $i$-diagonal, to the right of the line. By \autoref{L: Loop over negligible diagonal}, this subdiagram can be replaced with a sum of terms. One term contains the $i$-loop to the left of all the straight vertical $i$-strings and all other terms contain an $i$-crossing to the left of all the vertical $i$-strings such that the two left exiting strings of the crossing are connected to the points $(x_{s_1},1)$ and $(x_{s_1},0)$. Denote the first term as $D_1$ and the other terms $X_j^{(1)}$, for $j\in J_1$ some finite indexing set. So $A_{\mathfrak{t}_{\mathbf{d}}}^\lambda \circ (A_{\mathfrak{t}_{\mathbf{d}}^\diamond}^\lambda)^\ast = D_1 + \sum_{j\in J_1} X_j^{(1)}$. If in the diagram $D_1$ there are still subdiagrams $e_{\overline{\Delta}}$, where $\Delta$ is a negligible $i$-diagonal, to the left of the $i$-loop and to the right of the vertical line, then we can apply \autoref{L: Loop over negligible diagonal} to $D_1$ again. Hence, $D_1 = D_2 + \sum_{j\in J_2} X_j^{(2)}$, where $D_2$ is the diagram where the $i$-loop is to the left of the vertical $i$-strings and $X_j^{(2)}$, for $j\in J_2$ a finite indexing set, a diagram with an $i$-crossing to the left of the vertical $i$-strings such that the two left exiting strings of the crossing are connected to the points $(x_{s_1},1)$ and $(x_{s_1},0)$. So if in $\lambda$ there are $c$ negligible $i$-diagonals to the right of $s_1$ then $A_{\mathfrak{t}_{\mathbf{d}}}^\lambda \circ (A_{\mathfrak{t}_{\mathbf{d}}^\diamond}^\lambda)^\ast = D_c + \sum_{l=1}^{c} \sum_{j\in J_l} X_j^{(l)}$, by \autoref{L: Loop over negligible diagonal}. By the argument used in the case $c=0$, the $i$-loop in the diagram $D_c$ can be straightened with no extra terms, by \ref{D: Cherednik Algebra}\ref{D: Cherednik double crossings}. So $D_c$ can be transformed into
	
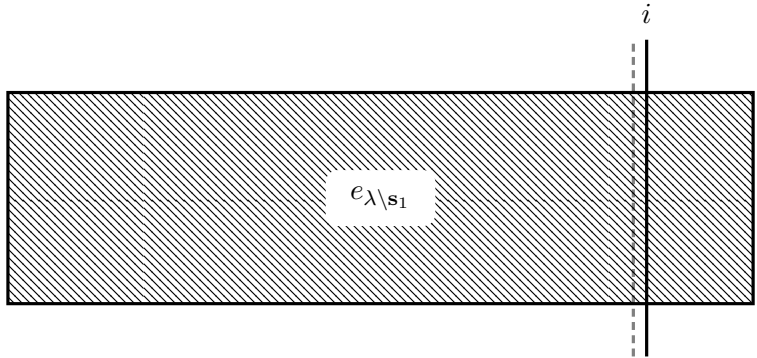
\[
	\begin{tikzpicture}[very thick,baseline,scale=0.35,rounded corners]\small
	\node at (20,13) {$i$};

	\draw[ghost] (19.5,0) +(0,0) -- +(0,12);
	\draw (20,0) +(0,0) -- +(0,12);

	\draw[pattern = north west lines,sharp corners] (-4,2) +(0,0) rectangle +(28,8);
	\draw[fill,color = white] (8,5) +(0,0) rectangle +(4,2);
	\node at (10,6) {$e_{\lambda\setminus \mathbf{s}_1}$};
	\end{tikzpicture}
	\]
	It remains to show that the terms $X_j^{(l)}$ `vanish'. In the diagram $X_j^{(l)}$, for $1\leq l\leq c$ and $j\in J_l$, there are $c-l$ subdiagrams of the form $e_{\overline{\Delta}}$, where $\Delta$ is a negligible $i$-diagonal, to the left of the position of the $i$-crossing and to the right of the vertical line at $x_{s_1}$. By \autoref{L: Crossing over negligible diagonal}, we can drag this crossing through all these subdiagrams and, by relation \ref{D: Cherednik Algebra}\ref{D: Cherednik triple crossings}, through any other straight vertical string to the right of the line at $x_{s_1}$ without creating any extra terms. Locally, the diagram then contains a subdiagram that is equal to one of the diagrams in \autoref{L: Loop over wall}. By \autoref{L: Loop over wall}, we can replace this subdiagram with one that contains two $i$-crossings. In the resulting diagram we take a small horizontal strip between the two $i$-crossings, where all strings are straight vertical. Consider this strip as a subdiagram. Then this subdiagram is a diagram of type $(\lambda',\lambda')$, where $\lambda \vartriangleleft_\theta \lambda'$, because $\lambda$ and $\lambda'$ are equal up to one $i$-string that is to the right of $x_{s_1}$ in $\lambda$ and to the left of $x_{s_1}$ in $\lambda'$. Hence, the whole diagram $X_j^{(l)}$ vanishes because it is of more dominant shape. Implying that all these terms vanish. So $A_{\mathfrak{t}_{\mathbf{d}}}^\lambda \circ (A_{\mathfrak{t}_{\mathbf{d}}^\diamond}^\lambda)^\ast = D_c = \pm e_\lambda(\residue(\mathfrak{t}_{\mathbf{d}}))$, completing the proof of the base case $m=1$ of the induction.

	Now let $m > 1$ and assume the statement is true for $m-1$. Because $m > 1$ we have at least one $i$-crossing of the solid strings of the staggered sequence and the diagram is:

	
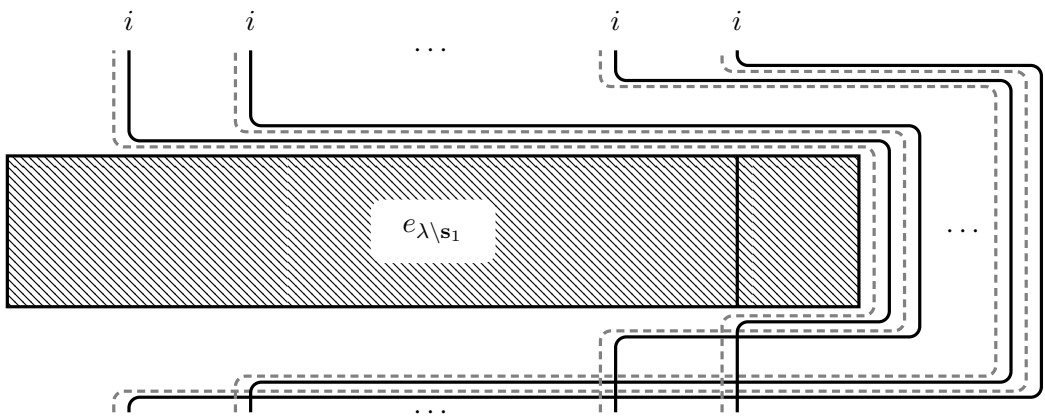
\[
	\begin{tikzpicture}[very thick,baseline,scale=0.4,rounded corners]\small
	\node at (0,13) {$i$};
	\node at (4,13) {$i$};
	\node at (16,13) {$i$};
	\node at (20,13) {$i$};

	\draw[ghost] (-0.5,0) +(0,0) -- +(0,0.7) -- +(30,0.7) -- +(30,11.3) -- +(20,11.3) -- +(20,12);
	\draw (0,0) +(0,0) -- +(0,0.5) -- +(30,0.5) -- +(30,11.5) -- +(20,11.5) -- +(20,12);

	\draw[ghost] (3.5,0) +(0,0) -- +(0,1.2) -- +(25,1.2) -- +(25,10.8) -- +(12,10.8) -- +(12,12);
	\draw (4,0) +(0,0) -- +(0,1) -- +(25,1) -- +(25,11) -- +(12,11) -- +(12,12);

	\draw[ghost] (15.5,0) +(0,0) -- +(0,2.7) -- +(10,2.7) -- +(10,9.3) -- +(-12,9.3) -- +(-12,12);
	\draw (16,0) +(0,0) -- +(0,2.5) -- +(10,2.5) -- +(10,9.5) -- +(-12,9.5) -- +(-12,12);

	\draw[ghost] (19.5,0) +(0,0) -- +(0,3.2) -- +(5,3.2) -- +(5,8.8) -- +(-20,8.8) -- +(-20,12);
	\draw (20,0) +(0,0) -- +(0,3) -- +(5,3) -- +(5,9) -- +(-20,9) -- +(-20,12);

	\node at (10,12) {$\cdots$};
	\node at (10,0) {$\cdots$};
	\node at (27.5,6) {$\cdots$};

	\draw[pattern = north west lines,sharp corners] (-4,3.5) +(0,0) rectangle +(28,5);
	\draw[fill,color = white] (8,5) +(0,0) rectangle +(4,2);
	\node at (10,6) {$e_{\lambda\setminus \mathbf{s}_1}$};

	\draw (20,3.5) -- (20,8.5);
	\end{tikzpicture}
	\]
	However, by the same argument as for case $m = 1$, this diagram can be transformed into:
	
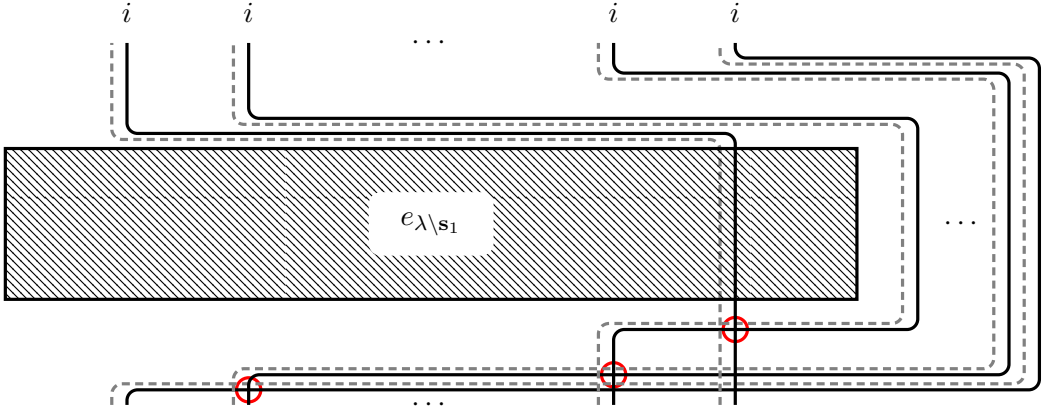
\[
	\begin{tikzpicture}[very thick,baseline,scale=0.4,rounded corners]\small
	\draw[red] (4,0.5) circle [radius=0.4];
	\draw[red] (16,1) circle [radius=0.4];
	\draw[red] (20,2.5) circle [radius=0.4];

	\node at (0,13) {$i$};
	\node at (4,13) {$i$};
	\node at (16,13) {$i$};
	\node at (20,13) {$i$};

	\draw[ghost] (-0.5,0) +(0,0) -- +(0,0.7) -- +(30,0.7) -- +(30,11.3) -- +(20,11.3) -- +(20,12);
	\draw (0,0) +(0,0) -- +(0,0.5) -- +(30,0.5) -- +(30,11.5) -- +(20,11.5) -- +(20,12);

	\draw[ghost] (3.5,0) +(0,0) -- +(0,1.2) -- +(25,1.2) -- +(25,10.8) -- +(12,10.8) -- +(12,12);
	\draw (4,0) +(0,0) -- +(0,1) -- +(25,1) -- +(25,11) -- +(12,11) -- +(12,12);

	\draw[ghost] (15.5,0) +(0,0) -- +(0,2.7) -- +(10,2.7) -- +(10,9.3) -- +(-12,9.3) -- +(-12,12);
	\draw (16,0) +(0,0) -- +(0,2.5) -- +(10,2.5) -- +(10,9.5) -- +(-12,9.5) -- +(-12,12);

	\draw[ghost] (19.5,0) +(0,0) -- +(0,8.8) -- +(-20,8.8) -- +(-20,12);
	\draw (20,0) +(0,0) -- +(0,9) -- +(-20,9) -- +(-20,12);

	\node at (10,12) {$\cdots$};
	\node at (10,0) {$\cdots$};
	\node at (27.5,6) {$\cdots$};

	\draw[pattern = north west lines,sharp corners] (-4,3.5) +(0,0) rectangle +(28,5);
	\draw[fill,color = white] (8,5) +(0,0) rectangle +(4,2);
	\node at (10,6) {$e_{\lambda\setminus \mathbf{s}_1}$};
	\end{tikzpicture}
	\]
	where the circles indicate important $i$-crossings. Note that their location is at the respective coordinate $x_{s_l}$, for $2\leq l\leq m$. Because $\mathbf{s}_1$ is a staggered sequence there are no addable $i$-boxes to the right of $s_1$ in $\lambda$ and $\{s_l | 2\leq l\leq m\}$ is the complete set of removable $i$-boxes to the right of $s_1$ in $\lambda$. So the circles indicate the locations of all removable $i$-diagonals $\Delta_l$ to the right of $x_{s_1}$ in $\lambda$, for $2\leq l\leq m$. Since every string in the diagram, except for the strings of the staggered sequence, are straight vertical, we have locally at those locations subdiagrams $e_{\overline{\Delta_l\setminus s_l}}$, where $\Delta_l\setminus s_l$ is an addable $i$-diagonal.

	The two left exiting strings of the rightmost circled $i$-crossing are connected to the points $(x_{s_1},1)$ and $(x_{s_{m-1}},0)$. At its position we have the subdiagram $e_{\overline{\Delta_m\setminus s_m}}$. By \autoref{L: Crossing over removable diagonal}, this subdiagram can be replaced with a sum of terms. One term contains an $i$-loop to the left of all the straight vertical $i$-strings in $e_{\overline{\Delta_m\setminus s_m}}$ such that it is connected to the points $(x_{s_1},1)$ and $(x_{s_{m-1}},0)$. Each of the other terms contain an $i$-crossing to the left of all the straight vertical $i$-strings such that the two left exiting strings of the crossing are connected to the points $(x_{s_1},1)$ and $(x_{s_{m-1}},0)$. Denote the first term as $D_0^{(1)}$ and the other terms $X_j^{(0,1)}$, for $j\in J_{(0,1)}$ a finite indexing set. So $A_{\mathfrak{t}_{\mathbf{d}}}^\lambda \circ (A_{\mathfrak{t}_{\mathbf{d}}^\diamond}^\lambda)^\ast = D_0^{(1)} + \sum_{j\in J_{(0,1)}} X_j^{(0,1)}$. 
	
	Let the number of negligible $i$-diagonals between the coordinates $x_{s_{m-t}}$ and $x_{s_{m-(t-1)}}$ be $c_t$. Then in the diagram $D_0^{(1)}$ there are $c_1$ many subdiagrams of the form $e_{\overline{\Delta}}$, where $\Delta$ is a negligible $i$-diagonal, to the left of the $i$-loop and right of $x_{s_{m-1}}$. Then, similar to the base case of the induction, there exist diagrams such that $D_0^{(1)} = D_{c_1}^{(1)} + \sum_{l = 1}^{c_1} \sum_{j\in J_{(l,1)}} X_j^{(l,1)}$, for finite indexing sets $J_{(l,1)}$ with $1\leq l\leq c_1$, by \autoref{L: Loop over negligible diagonal}. The diagram $D_{c_1}^{(1)}$ is the diagram
	
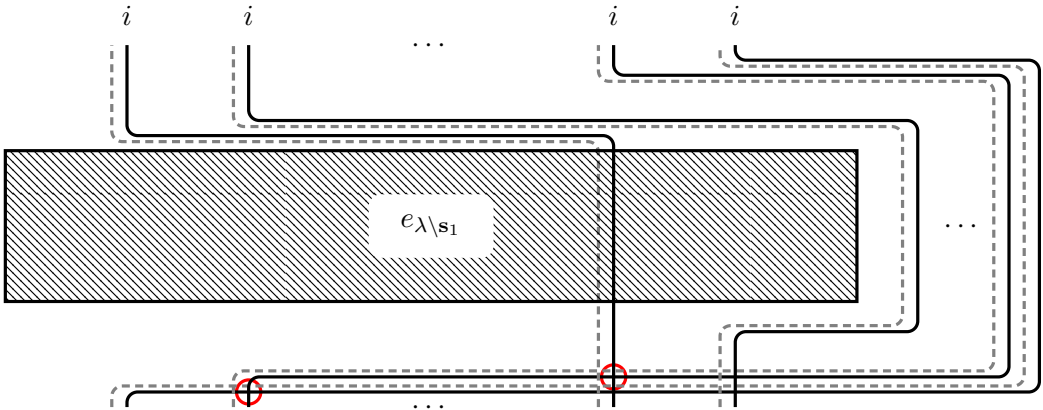
\[
	\begin{tikzpicture}[very thick,baseline,scale=0.4,rounded corners]\small
	\draw[red] (4,0.5) circle [radius=0.4];
	\draw[red] (16,1) circle [radius=0.4];

	\node at (0,13) {$i$};
	\node at (4,13) {$i$};
	\node at (16,13) {$i$};
	\node at (20,13) {$i$};

	\draw[ghost] (-0.5,0) +(0,0) -- +(0,0.7) -- +(30,0.7) -- +(30,11.3) -- +(20,11.3) -- +(20,12);
	\draw (0,0) +(0,0) -- +(0,0.5) -- +(30,0.5) -- +(30,11.5) -- +(20,11.5) -- +(20,12);

	\draw[ghost] (3.5,0) +(0,0) -- +(0,1.2) -- +(25,1.2) -- +(25,10.8) -- +(12,10.8) -- +(12,12);
	\draw (4,0) +(0,0) -- +(0,1) -- +(25,1) -- +(25,11) -- +(12,11) -- +(12,12);

	\draw[ghost] (15.5,0) +(0,0) --  +(0,8.8) -- +(-16,8.8) -- +(-16,12);
	\draw (16,0) +(0,0) --  +(0,9) -- +(-16,9) -- +(-16,12);

	\draw[ghost] (19.5,0) +(0,0) -- +(0,2.7) -- +(6,2.7) -- +(6,9.3) -- +(-16,9.3) -- +(-16,12);
	\draw (20,0) +(0,0) -- +(0,2.5) -- +(6,2.5) -- +(6,9.5) -- +(-16,9.5) -- +(-16,12);

	\node at (10,12) {$\cdots$};
	\node at (10,0) {$\cdots$};
	\node at (27.5,6) {$\cdots$};

	\draw[pattern = north west lines,sharp corners] (-4,3.5) +(0,0) rectangle +(28,5);
	\draw[fill,color = white] (8,5) +(0,0) rectangle +(4,2);
	\node at (10,6) {$e_{\lambda\setminus \mathbf{s}_1}$};
	\end{tikzpicture}
	\]
	In $D_{c_1}^{(1)}$ at the rightmost circled crossing we have the subdiagram $e_{\overline{\Delta_{m-1}\setminus s_{m-1}}}$. By \autoref{L: Crossing over removable diagonal}, we get $D_{c_1}^{(1)} = D_{0}^{(2)} + \sum_{j\in J_{(0,2)}} X_j^{(0,2)}$, for a finite indexing set $J_{(0,2)}$. By repetition of the argument we get
	\begin{equation*}
	A_{\mathfrak{t}_{\mathbf{d}}}^\lambda \circ (A_{\mathfrak{t}_{\mathbf{d}}^\diamond}^\lambda)^\ast = D_{c_{m-1}}^{(m-1)} + \sum_{b=1}^{m-1} \sum_{a=0}^{c_b} \sum_{j\in J_{(a,b)}} X_j^{(a,b)},
	\end{equation*}
	where $J_{(a,b)}$ is a finite indexing set, for $1\leq b\leq m-1$ and $0\leq a\leq c_b$. The diagram $D_{c_{m-1}}^{(m-1)}$ is
	
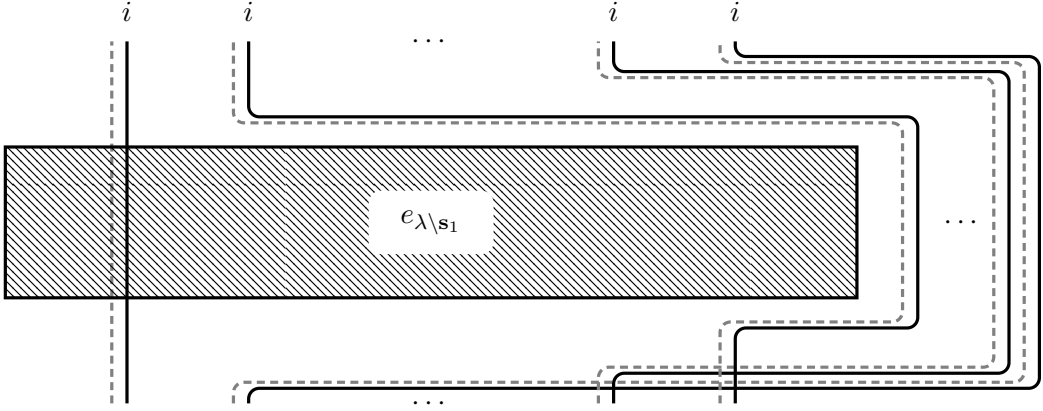
\[
	\begin{tikzpicture}[very thick,baseline,scale=0.4,rounded corners]\small
	\node at (0,13) {$i$};
	\node at (4,13) {$i$};
	\node at (16,13) {$i$};
	\node at (20,13) {$i$};

	\draw[ghost] (-0.5,0) +(0,0) -- +(0,12);
	\draw (0,0) +(0,0) -- +(0,12);

	\draw[ghost] (3.5,0) +(0,0) -- +(0,0.7) -- +(26,0.7) -- +(26,11.3) -- +(16,11.3) -- +(16,12);
	\draw (4,0) +(0,0) -- +(0,0.5) -- +(26,0.5) -- +(26,11.5) -- +(16,11.5) -- +(16,12);

	\draw[ghost] (15.5,0) +(0,0) -- +(0,1.2) -- +(13,1.2) -- +(13,10.8) -- +(0,10.8) -- +(0,12);
	\draw (16,0) +(0,0) -- +(0,1) -- +(13,1) -- +(13,11) -- +(0,11) -- +(0,12);

	\draw[ghost] (19.5,0) +(0,0) -- +(0,2.7) -- +(6,2.7) -- +(6,9.3) -- +(-16,9.3) -- +(-16,12);
	\draw (20,0) +(0,0) -- +(0,2.5) -- +(6,2.5) -- +(6,9.5) -- +(-16,9.5) -- +(-16,12);

	\node at (10,12) {$\cdots$};
	\node at (10,0) {$\cdots$};
	\node at (27.5,6) {$\cdots$};

	\draw[pattern = north west lines,sharp corners] (-4,3.5) +(0,0) rectangle +(28,5);
	\draw[fill,color = white] (8,5) +(0,0) rectangle +(4,2);
	\node at (10,6) {$e_{\lambda\setminus \mathbf{s}_1}$};
	\end{tikzpicture}
	\]
	Without the straightened vertical $i$-string this diagram corresponds to strings of the staggered sequence $\mathbf{s}_1' = (s_2,\ldots,s_m)$. So $|\mathbf{s}_1'| = m-1$. By induction, $D_{c_{m-1}}^{(m-1)}$ is equal to

	
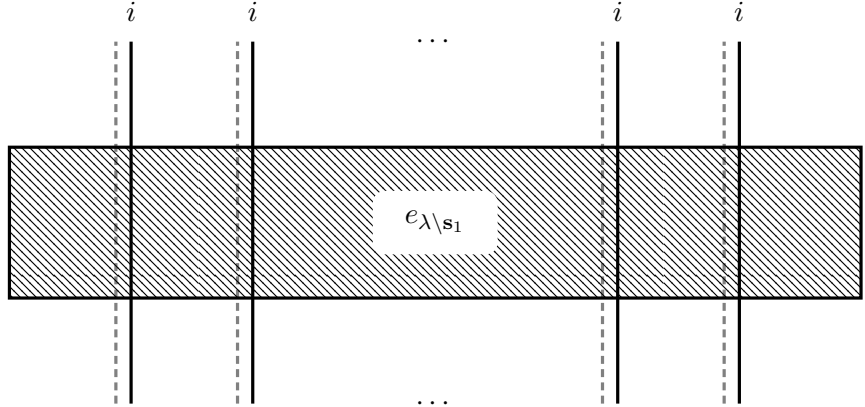
\[
	\begin{tikzpicture}[very thick,baseline,scale=0.4,rounded corners]\small
	\node at (0,13) {$i$};
	\node at (4,13) {$i$};
	\node at (16,13) {$i$};
	\node at (20,13) {$i$};

	\draw[ghost] (-0.5,0) +(0,0) -- +(0,12);
	\draw (0,0) +(0,0) -- +(0,12);

	\draw[ghost] (3.5,0) +(0,0) -- +(0,12);
	\draw (4,0) +(0,0) -- +(0,12);

	\draw[ghost] (15.5,0) +(0,0) -- +(0,12);
	\draw (16,0) +(0,0) -- +(0,12);

	\draw[ghost] (19.5,0) +(0,0) -- +(0,12);
	\draw (20,0) +(0,0) -- +(0,12);

	\node at (10,12) {$\cdots$};
	\node at (10,0) {$\cdots$};

	\draw[pattern = north west lines,sharp corners] (-4,3.5) +(0,0) rectangle +(28,5);
	\draw[fill,color = white] (8,5) +(0,0) rectangle +(4,2);
	\node at (10,6) {$e_{\lambda\setminus \mathbf{s}_1}$};
	\end{tikzpicture}
	\]

	It remains to show that the terms $X_j^{(a,b)}$, for $1\leq b\leq m-1$, $0\leq a\leq c_b$, and $j\in J_{(a,b)}$ a finite indexing set, `vanish'. First let $b=m-1$. The diagram $X_j^{(a,m-1)}$, for $1\leq a\leq c_{m-1}$ and $j\in J_{(a,m-1)}$ a finite indexing set, contains an $i$-crossing such that its position is to the left of $x_{s_2}$ and that its two exiting strings at the left are connected to the points $(x_{s_1},1)$ and $(x_{s_1},0)$. Between its position and $x_{s_1}$ there are $c_{m-1}-a$ many subdiagrams of the form $e_{\overline{\Delta}}$, where $\Delta$ is a negligible $i$-diagonal. Then following a similar argument to the one used in the base case, $X_j^{(a,m-1)}$ vanishes.

	So let $1\leq b < m-1$ and assume all terms with $b+1$ vanish. The diagram $X_j^{(a,b)}$, for $1\leq a\leq c_{b}$ and $j\in J_{(a,b)}$ a finite indexing set, contains an $i$-crossing such that its position is to the left of $x_{s_{m-(b-1)}}$ and that its two exiting strings at the left are connected to the points $(x_{s_1},1)$ and $(x_{s_{m-b}},0)$. Between its position and $x_{s_{m-b}}$ there are $c_{b}-a$ many subdiagrams of the form $e_{\overline{\Delta}}$, where $\Delta$ is a negligible $i$-diagonal. Also between its position and $x_{s_1}$ there are the subdiagrams $e_{\overline{\Delta_l\setminus s_l}}$, for $2\leq l \leq m-b$. By \autoref{L: Crossing over negligible diagonal}, we can drag the crossing up to the position $x_{s_{m-b}}$ without creating any extra terms. In $X_j^{(a,b)}$ we then have the following subdiagram
	
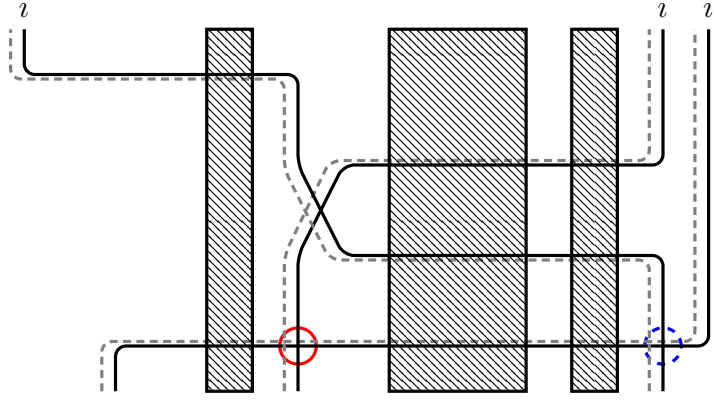
\[
	\begin{tikzpicture}[very thick,baseline,scale=0.6,rounded corners]
	\draw[red] (4,1) circle [radius=0.4];
	\draw[blue, dashed] (12,1) circle [radius=0.4];

	\node at (-2,8.5) {$i$};
	\node at (12,8.5) {$i$};
	\node at (13,8.5) {$i$};

	\draw[ghost] (-0.3,0) +(0,0) -- +(0,1.1) -- +(13,1.1) -- +(13,8);
	\draw (0,0) +(0,0) -- +(0,1) -- +(13,1) -- +(13,8);

	\draw[ghost] (3.7,0) +(0,0) -- +(0,3) -- +(1,5.1) -- +(8,5.1) -- +(8,8);
	\draw (4,0) +(0,0) -- +(0,3) -- +(1,5) -- +(8,5) -- +(8,8);

	\draw[ghost] (11.7,0) +(0,0) -- +(0,2.9) -- +(-7,2.9) -- +(-8,4.9) -- +(-8,6.9) -- +(-14,6.9) -- +(-14,8);
	\draw (12,0) +(0,0) -- +(0,3) -- +(-7,3) -- +(-8,5) -- +(-8,7) -- +(-14,7) -- +(-14,8);

	\draw[pattern = north west lines,sharp corners] (2,0) +(0,0) rectangle +(1,8);

	\draw[pattern = north west lines,sharp corners] (6,0) +(0,0) rectangle +(3,8);
	\draw[pattern = north west lines,sharp corners] (10,0) +(0,0) rectangle +(1,8);
	\end{tikzpicture}
	\]
	where the left circled crossing is the same circled crossing as in the diagram $A_{\mathfrak{t}_{\mathbf{d}}}^\lambda \circ (A_{\mathfrak{t}_{\mathbf{d}}^\diamond}^\lambda)^\ast$ at the position $x_{s_{m-b}}$. The shaded area on the left corresponds to the subdiagram $e_{\overline{\Delta_{m-b}\setminus s_{m-b}}}$. The shaded area in the middle contains all the $c_b - a$ many subdiagrams of the form $e_{\overline{\Delta}}$, where $\Delta$ is a negligible $i$-diagonal. And the shaded area on the right corresponds to a diagram $e_{\overline{\Delta'}}$, where $\Delta'$ is an addable $i$-diagonal. Note that $\Delta'$ is not an addable $i$-diagonal in $\lambda$, however, see the diagrams in the sums in \autoref{L: Crossing over removable diagonal} and \autoref{L: Loop over negligible diagonal} to check that the straight vertical strings are equivalent to the diagram for an addable $i$-diagonal.
	
	So locally with $e_{\overline{\Delta'}}$ this diagram has a subdiagram at the right circled $i$-crossing that is equal to one in \autoref{L: Crossing over removable diagonal}. By \autoref{L: Crossing over removable diagonal}, $X_j^{(a,b)} = E_{(j,0)}^{(a,b)} + \sum_{w\in W_{(j,0)}^{(a,b)}} Y_{(j,w,0)}^{(a,b)}$, for a finite indexing set $W_{(j,0)}^{(a,b)}$. The diagram $E_{(j,0)}^{(a,b)}$ is the one containing an $i$-loop to the left of the vertical $i$-strings in the shaded area on the right such that it is connected to the points $(x_{s_1},1)$ and $(x_{s_{m-(b+1)}},0)$. Each diagram $Y_{(j,w,0)}^{(a,b)}$ contains an $i$-crossing to the left of the vertical $i$-strings in the shaded area on the right such that its two exiting strings on the left are connected to the points $(x_{s_1},1)$ and $(x_{s_{m-(b+1)}},0)$. Similarly as before, $E_{(j,0)}^{(a,b)} = E_{(j,c_b - a)}^{(a,b)} + \sum_{v = 1}^{c_b - a} \sum_{w\in W_{(j,v)}^{(a,b)}} Y_{(j,w,v)}^{(a,b)}$, by \autoref{L: Loop over negligible diagonal}. The subdiagram in $E_{(j,c_b - a)}^{(a,b)}$ is
	\[
	\begin{tikzpicture}[very thick,baseline,scale=0.6,rounded corners]
	\draw[red] (4,1) circle [radius=0.4];

	\node at (-2,8.5) {$i$};
	\node at (12,8.5) {$i$};
	\node at (13,8.5) {$i$};

	\draw[ghost] (-0.3,0) +(0,0) -- +(0,1.1) -- +(5,1.1) -- +(5,3) -- +(4,5) -- +(4,6.9) -- +(-2,6.9) -- +(-2,8);
	\draw (0,0) +(0,0) -- +(0,1) -- +(5,1) -- +(5,3) -- +(4,5) -- +(4,7) -- +(-2,7) -- +(-2,8);

	\draw[ghost] (3.7,0) +(0,0) -- +(0,3) -- +(1,5.1) -- +(8,5.1) -- +(8,8);
	\draw (4,0) +(0,0) -- +(0,3) -- +(1,5) -- +(8,5) -- +(8,8);

	\draw[ghost] (11.7,0) +(0,0) -- +(0,1.1) -- +(1,1.1) -- +(1,8);
	\draw (12,0) +(0,0) -- +(0,1) -- +(1,1) -- +(1,8);

	\draw[pattern = north west lines,sharp corners] (2,0) +(0,0) rectangle +(1,8);

	\draw[pattern = north west lines,sharp corners] (6,0) +(0,0) rectangle +(3,8);
	\draw[pattern = north west lines,sharp corners] (10,0) +(0,0) rectangle +(1,8);
	\end{tikzpicture}
	\]
	which vanishes, by relation \ref{D: Cherednik Algebra}\ref{D: Cherednik double crossings}\ref{D: Cherednik double i crossing}. By \autoref{L: Crossing over negligible diagonal}, the $i$-crossing in the diagram $Y_{(j,w,v)}^{(a,b)}$, for $0\leq v\leq c_b-a$ and $w\in W_{(j,v)}^{(a,b)}$ a finite indexing set, can be dragged through the shaded area in the middle up to the position $x_{s_{m-b}}$. So the subdiagram in $Y_{(j,w,v)}^{(a,b)}$ is

	
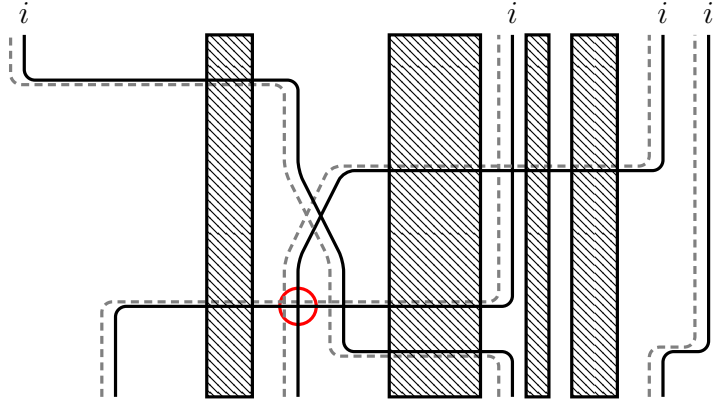
\[
	\begin{tikzpicture}[very thick,baseline,scale=0.6,rounded corners]
	\draw[red] (4,2) circle [radius=0.4];

	\node at (-2,8.5) {$i$};
	\node at (8.7,8.5) {$i$};
	\node at (12,8.5) {$i$};
	\node at (13,8.5) {$i$};

	\draw[ghost] (-0.3,0) +(0,0) -- +(0,2.1) -- +(8.7,2.1) -- +(8.7,8);
	\draw (0,0) +(0,0) -- +(0,2) -- +(8.7,2) -- +(8.7,8);

	\draw[ghost] (3.7,0) +(0,0) -- +(0,3) -- +(1,5.1) -- +(8,5.1) -- +(8,8);
	\draw (4,0) +(0,0) -- +(0,3) -- +(1,5) -- +(8,5) -- +(8,8);

	\draw[ghost] (8.4,0) +(0,0) -- +(0,0.9) -- +(-3.7,0.9) -- +(-3.7,3) -- +(-4.7,5) -- +(-4.7,6.9) -- +(-10.7,6.9) -- +(-10.7,8);
	\draw (8.7,0) +(0,0) -- +(0,1) -- +(-3.7,1) -- +(-3.7,3) -- +(-4.7,5) -- +(-4.7,7) -- +(-10.7,7) -- +(-10.7,8);

	\draw[ghost] (11.7,0) +(0,0) -- +(0,1.1) -- +(1,1.1) -- +(1,8);
	\draw (12,0) +(0,0) -- +(0,1) -- +(1,1) -- +(1,8);

	\draw[pattern = north west lines,sharp corners] (2,0) +(0,0) rectangle +(1,8);

	\draw[pattern = north west lines,sharp corners] (6,0) +(0,0) rectangle +(2,8);

	\draw[pattern = north west lines,sharp corners] (9,0) +(0,0) rectangle +(0.5,8);
	\draw[pattern = north west lines,sharp corners] (10,0) +(0,0) rectangle +(1,8);
	\end{tikzpicture}
	\]
	Applying relation \ref{D: Cherednik Algebra}\ref{D: Cherednik triple crossings} gives the subdiagram
	
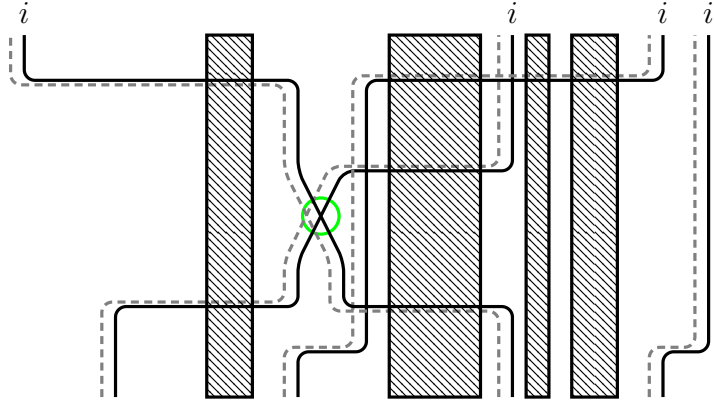
\[
	\begin{tikzpicture}[very thick,baseline,scale=0.6,rounded corners]
	\draw[green] (4.5,4) circle [radius=0.4];

	\node at (-2,8.5) {$i$};
	\node at (8.7,8.5) {$i$};
	\node at (12,8.5) {$i$};
	\node at (13,8.5) {$i$};

	\draw[ghost] (-0.3,0) +(0,0) -- +(0,2.1) -- +(4,2.1) -- +(4,3) -- +(5,5.1) -- +(8.7,5.1) -- +(8.7,8);
	\draw (0,0) +(0,0) -- +(0,2) -- +(4,2) -- +(4,3) -- +(5,5) -- +(8.7,5) -- +(8.7,8);

	\draw[ghost] (3.7,0) +(0,0) -- +(0,1.1) -- +(1.5,1.1) -- +(1.5,7.1) -- +(8,7.1) -- +(8,8);
	\draw (4,0) +(0,0) -- +(0,1) -- +(1.5,1) -- +(1.5,7) -- +(8,7) -- +(8,8);

	\draw[ghost] (8.4,0) +(0,0) -- +(0,1.9) -- +(-3.7,1.9) -- +(-3.7,3) -- +(-4.7,5) -- +(-4.7,6.9) -- +(-10.7,6.9) -- +(-10.7,8);
	\draw (8.7,0) +(0,0) -- +(0,2) -- +(-3.7,2) -- +(-3.7,3) -- +(-4.7,5) -- +(-4.7,7) -- +(-10.7,7) -- +(-10.7,8);

	\draw[ghost] (11.7,0) +(0,0) -- +(0,1.1) -- +(1,1.1) -- +(1,8);
	\draw (12,0) +(0,0) -- +(0,1) -- +(1,1) -- +(1,8);

	\draw[pattern = north west lines,sharp corners] (2,0) +(0,0) rectangle +(1,8);

	\draw[pattern = north west lines,sharp corners] (6,0) +(0,0) rectangle +(2,8);

	\draw[pattern = north west lines,sharp corners] (9,0) +(0,0) rectangle +(0.5,8);
	\draw[pattern = north west lines,sharp corners] (10,0) +(0,0) rectangle +(1,8);
	\end{tikzpicture}
	\]
	The shaded area on the left corresponds to the subdiagram $e_{\overline{\Delta_{m-b}\setminus s_{m-b}}}$. Hence, considering the circled $i$-crossing we have locally a subdiagram that is equal to one of \autoref{L: Crossing over removable diagonal}. By \autoref{L: Crossing over removable diagonal}, we can replace this subdiagram with a sum of terms such that in this sum there is only one term in which the rightmost $i$-string in the subdiagram is not straight vertical; see \autoref{L: Crossing over removable diagonal}. All terms where the rightmost $i$-string is straight vertical correspond to a diagram where we can replace the circled $i$-crossing above with two straight vertical $i$-strings. By relation \ref{D: Cherednik Algebra}\ref{D: Cherednik double crossings}\ref{D: Cherednik double i crossing}, all of these terms vanish. So the remaining term and, thus, the subdiagram in $Y_{(j,w,v)}^{(a,b)}$ is
	
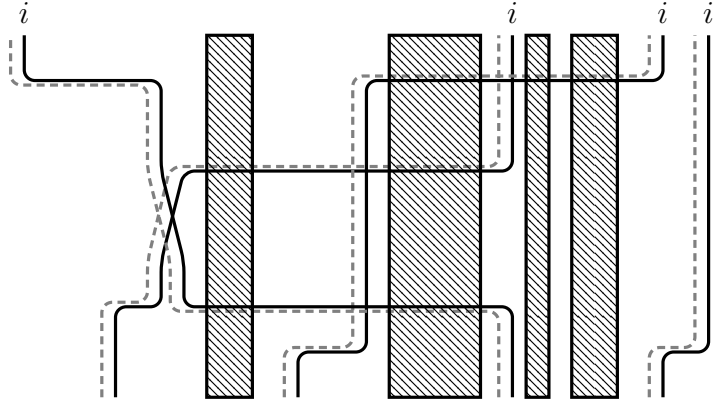
\[
	\begin{tikzpicture}[very thick,baseline,scale=0.6,rounded corners]
	\node at (-2,8.5) {$i$};
	\node at (8.7,8.5) {$i$};
	\node at (12,8.5) {$i$};
	\node at (13,8.5) {$i$};

	\draw[ghost] (-0.3,0) +(0,0) -- +(0,2.1) -- +(1,2.1) -- +(1,3) -- +(1.5,5.1) -- +(8.7,5.1) -- +(8.7,8);
	\draw (0,0) +(0,0) -- +(0,2) -- +(1,2) -- +(1,3) -- +(1.5,5) -- +(8.7,5) -- +(8.7,8);

	\draw[ghost] (3.7,0) +(0,0) -- +(0,1.1) -- +(1.5,1.1) -- +(1.5,7.1) -- +(8,7.1) -- +(8,8);
	\draw (4,0) +(0,0) -- +(0,1) -- +(1.5,1) -- +(1.5,7) -- +(8,7) -- +(8,8);

	\draw[ghost] (8.4,0) +(0,0) -- +(0,1.9) -- +(-7.2,1.9) -- +(-7.2,3) -- +(-7.7,5) -- +(-7.7,6.9) -- +(-10.7,6.9) -- +(-10.7,8);
	\draw (8.7,0) +(0,0) -- +(0,2) -- +(-7.2,2) -- +(-7.2,3) -- +(-7.7,5) -- +(-7.7,7) -- +(-10.7,7) -- +(-10.7,8);

	\draw[ghost] (11.7,0) +(0,0) -- +(0,1.1) -- +(1,1.1) -- +(1,8);
	\draw (12,0) +(0,0) -- +(0,1) -- +(1,1) -- +(1,8);

	\draw[pattern = north west lines,sharp corners] (2,0) +(0,0) rectangle +(1,8);

	\draw[pattern = north west lines,sharp corners] (6,0) +(0,0) rectangle +(2,8);

	\draw[pattern = north west lines,sharp corners] (9,0) +(0,0) rectangle +(0.5,8);
	\draw[pattern = north west lines,sharp corners] (10,0) +(0,0) rectangle +(1,8);
	\end{tikzpicture}
	\]
	But this is a diagram with an $i$-crossing to the left of $x_{s_{m-b}}$ and only the $m-(b+2)$ many subdiagrams $e_{\overline{\Delta_l\setminus s_l}}$, for $2\leq l \leq m-(b+1)$, between the crossing and $x_{s_1}$. Such a diagram is equivalent to one of the form $X_{j'}^{(a',b+1)}$, hence, it vanishes by induction. So
	\begin{equation*}
	X_{j}^{(a,b)} = E_{(j,c_b - a)}^{(a,b)} + \sum_{v = 0}^{c_b - a} \sum_{w\in W_{(j,v)}^{(a,b)}} Y_{(j,w,v)}^{(a,b)} = 0 + \sum_{v = 0}^{c_b - a} \sum_{w\in W_{(j,v)}^{(a,b)}} 0 = 0
	\end{equation*}
	Therefore,
	\begin{equation*}
	A_{\mathfrak{t}_{\mathbf{d}}}^\lambda \circ (A_{\mathfrak{t}_{\mathbf{d}}^\diamond}^\lambda)^\ast = D_{c_{m-1}}^{(m-1)} + \sum_{b=1}^{m-1} \sum_{a=0}^{c_b} \sum_{j\in J_{(a,b)}} X_j^{(a,b)} = D_{c_{m-1}}^{(m-1)} = \pm e_\lambda(\residue(\mathfrak{t}_{\mathbf{d}})),
	\end{equation*}
	which is what we wanted to show.
\end{proof}

We now describe a second, a priori, different construction of a diagram of the element $c_{\mathfrak{t}_\mathbf{d}\cdot\sigma}^\lambda$.

\emph{Construction B:} Again draw the red strings straight vertically. For the solid and ghost strings we first take the diagram $A_{\mathfrak{t}_{\mathbf{d}}}^\lambda$ of $c_{\mathfrak{t}_{\mathbf{d}}}^\lambda$ and scale it vertically by the factor $\frac{1}{2}$. Then we use this diagram as the top half of this construction for $c_{\mathfrak{t}_{\mathbf{d}}\cdot \sigma}^\lambda$. In the bottom half we add the permutation $\sigma = (\sigma_1,\ldots,\sigma_k)\in\mathfrak{S}_{\mathbf{d}}$ as a braid diagram with minimal crossings such that the corresponding solid and ghost strings stay within the boxed areas below. Note, that $\sigma_j$ is the permutation operating on the staggered sequence $\mathbf{s}_j$, hence, the boxed areas only contain solid and ghost strings of the same residue as $\mathbf{s}_j$:
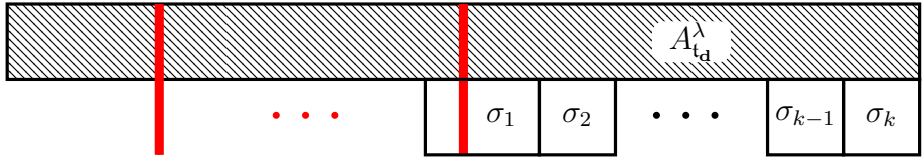
\[
\begin{tikzpicture}[very thick,baseline,rounded corners]
\draw[pattern = north west lines,sharp corners] (-4,1) +(0,0) rectangle +(12,1);
\draw[fill,color = white] (4.5,1.25) +(0,0) rectangle +(1,0.5);
\node at (5,1.5) {$A_{\mathfrak{t}_{\mathbf{d}}}^\lambda$};

\draw[wei] (-2,0) +(0,0) -- +(0,2);
\draw[wei] (2,0) +(0,0) -- +(0,2);

\Huge
\node[red] at (0,0.5) {$\cdots$};
\normalsize

\draw[sharp corners] (1.5,0) +(0,0) rectangle +(1.5,1);
\draw[sharp corners] (3,0) +(0,0) rectangle +(1,1);
\draw[sharp corners] (6,0) +(0,0) rectangle +(1,1);
\draw[sharp corners] (7,0) +(0,0) rectangle +(1,1);

\Huge
\node at (5,0.5) {$\cdots$};
\normalsize

\node at (2.5,0.5) {$\sigma_1$};
\node at (3.5,0.5) {$\sigma_2$};
\node at (6.5,0.5) {$\sigma_{k-1}$};
\node at (7.5,0.5) {$\sigma_k$};
\end{tikzpicture}
\]
Define $B_{\mathfrak{t}_{\mathbf{d}}\cdot\sigma}^\lambda$ to be the diagram for $c_{\mathfrak{t}_{\mathbf{d}}\cdot\sigma}^\lambda$ obtained using Construction B. The benefit of this construction is that it highlights the action of the permutation $\sigma\in\mathfrak{S}_{\mathbf{d}}$. Note that the minimal braid diagram for the permutation $\sigma_j$, for $1\leq j\leq k$, is not unique up to isotopy but recall that a reduced expression is unique up to braid moves. Since $\sigma_j$ only permutes boxes of the same residue, the choice of minimal braid diagram does not matter, by relation \ref{D: Cherednik Algebra}\ref{D: Cherednik triple crossings}.

\begin{Example}\label{E: Construction B}
	Let everything be as in \autoref{E: Presentation 1}. Then the diagram $B_{\mathfrak{t}_{\mathbf{d}}^\diamond}^\lambda$ of $c_{\mathfrak{t}_{\mathbf{d}}^\diamond}^\lambda$ is:
	
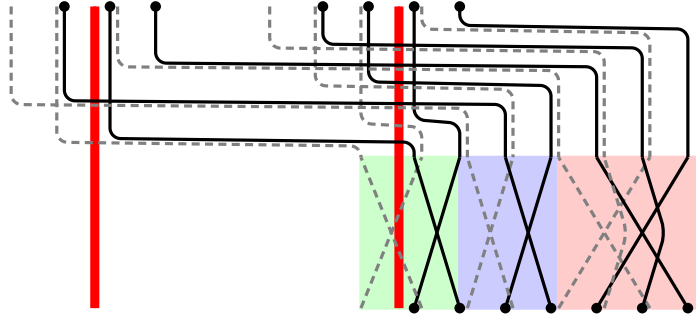
\[
	\begin{tikzpicture}[very thick,baseline,rounded corners]
	\filldraw[color=green!20,sharp corners] (3.5,-2) +(0,0) rectangle +(1.3,2);
	\filldraw[color=blue!20,sharp corners] (4.8,-2) +(0,0) rectangle +(1.3,2);
	\filldraw[color=red!20,sharp corners] (6.1,-2) +(0,0) rectangle +(1.8,2);

	\draw[wei](0,-2) -- (0,2);
	\draw[wei](4,-2) -- (4,2);

	\fill (-0.4,2) circle (2pt);
	\fill (0.2,2) circle (2pt);
	\fill (0.8,2) circle (2pt);

	\fill (4.2,2) circle (2pt);
	\fill (3.6,2) circle (2pt);
	\fill (3,2) circle (2pt);
	\fill (4.8,2) circle (2pt);

	\fill (4.2,-2) circle (2pt);
	\fill (4.8,-2) circle (2pt);
	\fill (5.4,-2) circle (2pt);
	\fill (6,-2) circle (2pt);
	\fill (6.6,-2) circle (2pt);
	\fill (7.2,-2) circle (2pt);
	\fill (7.8,-2) circle (2pt);

	\draw[ghost] (-0.5,2) +(0,0) -- +(0,-1.8) -- +(4,-1.85) -- +(4,-2);
	\draw[solid] (0.2,2) +(0,0) -- +(0,-1.75) -- +(4,-1.8) -- +(4,-2);

	\draw[ghost] (3.5,2) +(0,0) -- +(0,-1.55) -- +(0.8,-1.6) -- +(0.8,-2);
	\draw[solid] (4.2,2) +(0,0) -- +(0,-1.5) -- +(0.6,-1.55) -- +(0.6,-2);

	\draw[ghost] (-1.1,2) +(0,0) -- +(0,-1.3) -- +(6,-1.35) -- +(6,-2);
	\draw[solid] (-0.4,2) +(0,0) -- +(0,-1.25) -- +(5.8,-1.3) -- +(5.8,-2);

	\draw[ghost] (2.9,2) +(0,0) -- +(0,-1.05) -- +(2.6,-1.1) -- +(2.6,-2);
	\draw[solid] (3.6,2) +(0,0) -- +(0,-1) -- +(2.4,-1.05) -- +(2.4,-2);

	\draw[ghost] (0.3,2) +(0,0) -- +(0,-0.8) -- +(5.8,-0.85) -- +(5.8,-2);
	\draw[solid] (0.8,2) +(0,0) -- +(0,-0.75) -- +(5.8,-0.8) -- +(5.8,-2);

	\draw[ghost] (2.3,2) +(0,0) -- +(0,-0.55) -- +(4.4,-0.6) -- +(4.4,-2);
	\draw[solid] (3,2) +(0,0) -- +(0,-0.5) -- +(4.2,-0.55) -- +(4.2,-2);

	\draw[ghost] (4.3,2) +(0,0) -- +(0,-0.3) -- +(3,-0.35) -- +(3,-2);
	\draw[solid] (4.8,2) +(0,0) -- +(0,-0.25) -- +(3,-0.3) -- +(3,-2);

	\draw[ghost] (3.5,-2) +(0,0) -- +(0.8,2);
	\draw[solid] (4.2,-2) +(0,0) -- +(0.6,2);

	\draw[ghost] (4.3,-2) +(0,0) -- +(-0.8,2);
	\draw[solid] (4.8,-2) +(0,0) -- +(-0.6,2);

	\draw[ghost] (4.9,-2) +(0,0) -- +(0.6,2);
	\draw[solid] (5.4,-2) +(0,0) -- +(0.6,2);

	\draw[ghost] (5.5,-2) +(0,0) -- +(-0.6,2);
	\draw[solid] (6,-2) +(0,0) -- +(-0.6,2);

	\draw[ghost] (6.1,-2) +(0,0) -- +(1.2,2);
	\draw[ghost] (7.3,-2) +(0,0) -- +(-1.2,2);
	\draw[solid] (6.6,-2) +(0,0) -- +(1.2,2);
	\draw[solid] (7.8,-2) +(0,0) -- +(-1.2,2);

	\draw[ghost] (6.7,-2) +(0,0) -- +(0.3,1) -- +(0,2);
	\draw[solid] (7.2,-2) +(0,0) -- +(0.3,1) -- +(0,2);
	\end{tikzpicture}
	\]
	The three shadings highlight the boxed areas containing the braid diagrams of the permutations $\sigma_1$, $\sigma_2$, and $\sigma_3$ in $\sigma = (\sigma_1,\sigma_2,\sigma_3)$.
\end{Example}

\begin{Lemma}\label{L: Both presentations same}
	Let $\lambda\in\mathscr{U}_n^\ell(\theta)$ and $\mathbf{d}$ a decomposition path for $\lambda$. As elements of $\mathbb{A}_n^\kappa(\theta)$, $A_{\mathfrak{t}_{\mathbf{d}}\cdot\sigma}^\lambda = B_{\mathfrak{t}_{\mathbf{d}}\cdot\sigma}^\lambda$, for all $\sigma\in\mathfrak{S}_\mathbf{d}$.
\end{Lemma}

\begin{proof}
	We start with Construction A and transform it into Construction B. We show how it works when $\sigma = \omega_0^{\mathbf{d}}$, as for all other permutations $\sigma\in\mathfrak{S}_{\mathbf{d}}$ the argument is similar.

	By the definition of a staggered sequence, all of the non-vertical strings in the horizontal strips of Construction A have the same residue. Hence, we first focus on the solid and ghost strings of a single staggered sequence inside the diagram. This is a subdiagram of the form:
	
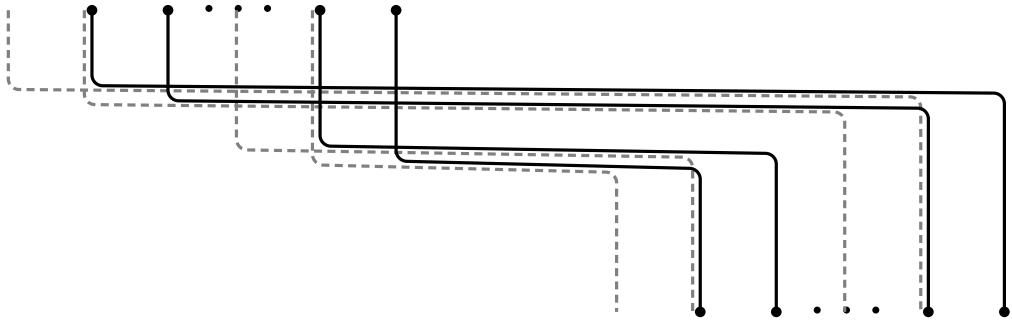
\[
	\begin{tikzpicture}[very thick,baseline,rounded corners]
	\fill (-6,4) circle (2pt);
	\fill (-5,4) circle (2pt);
	\fill (-3,4) circle (2pt);
	\fill (-2,4) circle (2pt);

	\Huge
	\node at (-4,4) {$\cdots$};
	\normalsize

	\fill (2,0) circle (2pt);
	\fill (3,0) circle (2pt);
	\fill (5,0) circle (2pt);
	\fill (6,0) circle (2pt);

	\Huge
	\node at (4,0) {$\cdots$};
	\normalsize

	\draw[ghost] (-7.1,4)  +(0,0) -- +(0,-1.05) -- +(12,-1.15) -- +(12,-4);
	\draw[ghost] (-6.1,4)  +(0,0) -- +(0,-1.25) -- +(10,-1.35) -- +(10,-4);
	\draw[ghost] (-4.1,4)  +(0,0) -- +(0,-1.85) -- +(6,-1.95) -- +(6,-4);
	\draw[ghost] (-3.1,4)  +(0,0) -- +(0,-2.05) -- +(4,-2.15) -- +(4,-4);

	\draw[solid] (-6,4)  +(0,0) -- +(0,-1) -- +(12,-1.1) -- +(12,-4);
	\draw[solid] (-5,4)  +(0,0) -- +(0,-1.2) -- +(10,-1.3) -- +(10,-4);
	\draw[solid] (-3,4)  +(0,0) -- +(0,-1.8) -- +(6,-1.9) -- +(6,-4);
	\draw[solid] (-2,4)  +(0,0) -- +(0,-2) -- +(4,-2.1) -- +(4,-4);
	\end{tikzpicture}
	\]
	By changing the heights of the almost horizontal strings, using relation \ref{D: Cherednik Algebra}\ref{D: Cherednik double crossings}, this subdiagram can be transformed into:
	
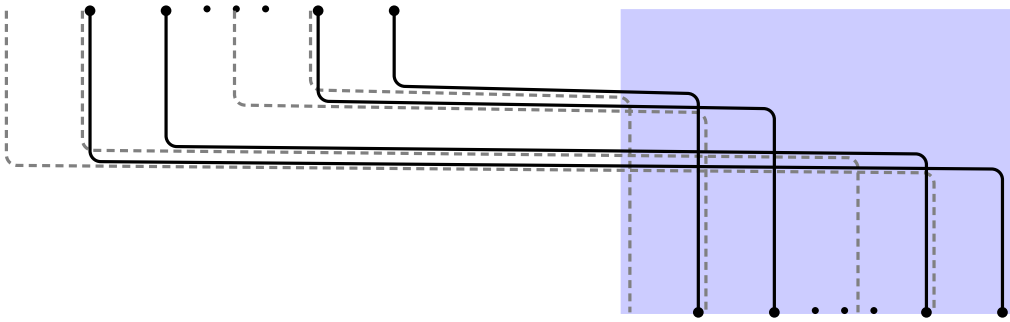
\[
	\begin{tikzpicture}[very thick,baseline,rounded corners]
	\filldraw[color=blue!20,sharp corners] (1,0) +(0,0) rectangle +(5.1,4);

	\fill (-6,4) circle (2pt);
	\fill (-5,4) circle (2pt);
	\fill (-3,4) circle (2pt);
	\fill (-2,4) circle (2pt);

	\Huge
	\node at (-4,4) {$\cdots$};
	\normalsize

	\fill (2,0) circle (2pt);
	\fill (3,0) circle (2pt);
	\fill (5,0) circle (2pt);
	\fill (6,0) circle (2pt);

	\Huge
	\node at (4,0) {$\cdots$};
	\normalsize

	\draw[ghost] (-7.1,4)  +(0,0) -- +(0,-2.05) -- +(12.2,-2.15) -- +(12.2,-4);
	\draw[ghost] (-6.1,4)  +(0,0) -- +(0,-1.85) -- +(10.2,-1.95) -- +(10.2,-4);
	\draw[ghost] (-4.1,4)  +(0,0) -- +(0,-1.25) -- +(6.2,-1.35) -- +(6.2,-4);
	\draw[ghost] (-3.1,4)  +(0,0) -- +(0,-1.05) -- +(4.2,-1.15) -- +(4.2,-4);

	\draw[solid] (-6,4)  +(0,0) -- +(0,-2) -- +(12,-2.1) -- +(12,-4);
	\draw[solid] (-5,4)  +(0,0) -- +(0,-1.8) -- +(10,-1.9) -- +(10,-4);
	\draw[solid] (-3,4)  +(0,0) -- +(0,-1.2) -- +(6,-1.3) -- +(6,-4);
	\draw[solid] (-2,4)  +(0,0) -- +(0,-1) -- +(4,-1.1) -- +(4,-4);
	\end{tikzpicture}
	\]
	However, the same manipulation actually works in the whole diagram $A_{\mathfrak{t}_{\mathbf{d}}\cdot\omega_0^\mathbf{d}}^\lambda$ because by using Construction A all other strings, in particular those of other residues, are vertical in the horizontal strip we are working in. Hence, when we apply relation \ref{D: Cherednik Algebra}\ref{D: Cherednik triple crossings}, the braid relation applies exactly and no correction terms are introduced.

	The shaded area in the diagram above indicates where the permutation $\sigma_j$ acts on the strings of the staggered sequence. The left and right bounds of it are exactly the left and right bound of the boxed area for the braid diagram of $\sigma_j$ in Construction B. Using isotopy to transform the diagram into:
	
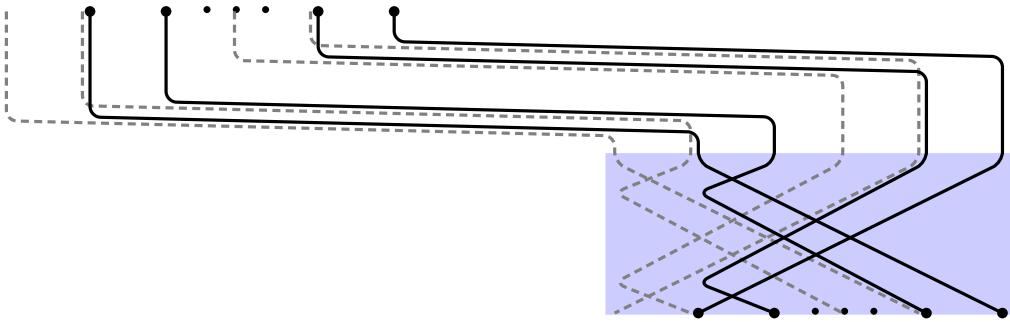
\[
	\begin{tikzpicture}[very thick,baseline,rounded corners]
	\filldraw[color=blue!20,sharp corners] (0.8,0) +(0,0) rectangle +(5.3,2.1);

	\fill (-6,4) circle (2pt);
	\fill (-5,4) circle (2pt);
	\fill (-3,4) circle (2pt);
	\fill (-2,4) circle (2pt);

	\Huge
	\node at (-4,4) {$\cdots$};
	\normalsize

	\fill (2,0) circle (2pt);
	\fill (3,0) circle (2pt);
	\fill (5,0) circle (2pt);
	\fill (6,0) circle (2pt);

	\Huge
	\node at (4,0) {$\cdots$};
	\normalsize

	\draw[ghost] (-7.1,4)  +(0,0) -- +(0,-1.45) -- +(8,-1.65) -- +(8,-2) -- +(12,-4);
	\draw[ghost] (-6.1,4)  +(0,0) -- +(0,-1.25) -- +(8,-1.45) -- +(8,-2) -- +(7,-2.4) -- +(10,-4);
	\draw[ghost] (-4.1,4)  +(0,0) -- +(0,-0.65) -- +(8,-0.85) -- +(8,-2) -- +(5,-3.6) -- +(6,-4);
	\draw[ghost] (-3.1,4)  +(0,0) -- +(0,-0.45) -- +(8,-0.65) -- +(8,-2) -- +(4,-4);

	\draw[solid] (-6,4)  +(0,0) -- +(0,-1.4) -- +(8,-1.6) -- +(8,-2) -- +(12,-4);
	\draw[solid] (-5,4)  +(0,0) -- +(0,-1.2) -- +(8,-1.4) -- +(8,-2) -- +(7,-2.4) -- +(10,-4);
	\draw[solid] (-3,4)  +(0,0) -- +(0,-0.6) -- +(8,-0.8) -- +(8,-2) -- +(5,-3.6) -- +(6,-4);
	\draw[solid] (-2,4)  +(0,0) -- +(0,-0.4) -- +(8,-0.6) -- +(8,-2) -- +(4,-4);
	\end{tikzpicture}
	\]
	Note that this works in the whole diagram $A_{\mathfrak{t}_{\mathbf{d}}\cdot\omega_0^\mathbf{d}}^\lambda$ because, by using Construction A, the area in $A_{\mathfrak{t}_{\mathbf{d}}\cdot\omega_0^\mathbf{d}}^\lambda$, corresponding to the shaded area indicated above, only contains the solid and ghost strings of the staggered sequence. But the diagram above is now in form of Construction B. Hence, applying these manipulations to the strings of all other staggered sequences gives us the diagram $B_{\mathfrak{t}_{\mathbf{d}}\cdot\omega_0^\mathbf{d}}^\lambda$.
\end{proof}

\begin{Proposition}\label{P: Chap 5 Prop}
	Let $\lambda\in\mathscr{U}_n^\ell(\theta)$, $\mathbf{d}$ a decomposition path for $\lambda$, and $\sigma,\tau\in\mathfrak{S}_{\mathbf{d}}$. Then $|\langle B_{\mathfrak{t}_{\mathbf{d}}\cdot\sigma}^\lambda, B_{\mathfrak{t}_{\mathbf{d}}^\diamond\cdot\tau}^\lambda\rangle_\lambda| = \delta_{\sigma^{-1},\tau}$.
\end{Proposition}

\begin{proof}
	Recall that $\mathfrak{t}_{\mathbf{d}}^\diamond = \mathfrak{t}_{\mathbf{d}}\cdot\omega_0^{\mathbf{d}}$ and choose a reduced expression for $\omega_0^{\mathbf{d}}\cdot\tau$. The diagram for the inner product is $B_{\mathfrak{t}_\mathbf{d}\cdot \sigma}^\lambda \circ (B_{\mathfrak{t}_\mathbf{d}^\diamond \cdot \tau}^\lambda)^\ast$:
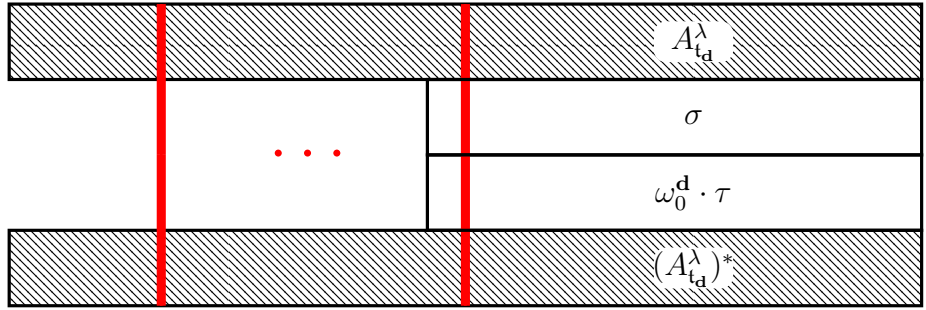
	\[
	\begin{tikzpicture}[very thick,baseline,rounded corners]
	\draw[pattern = north west lines,sharp corners] (-4,1) +(0,0) rectangle +(12,1);
	\draw[fill,color = white] (4.5,1.25) +(0,0) rectangle +(1,0.5);
	\node at (5,1.5) {$A_{\mathfrak{t}_\mathbf{d}}^\lambda$};

	\draw[wei] (-2,0) +(0,0) -- +(0,2);
	\draw[wei] (2,0) +(0,0) -- +(0,2);

	\Huge
	\node[red] at (0,0) {$\cdots$};
	\normalsize

	\draw[sharp corners] (1.5,0) +(0,0) rectangle +(6.5,1);

	\node at (5,0.5) {$\sigma$};

	\draw[pattern = north west lines,sharp corners] (-4,-1) +(0,0) rectangle +(12,-1);
	\draw[fill,color = white] (4.5,-1.25) +(0,0) rectangle +(1,-0.5);
	\node at (5,-1.5) {$(A_{\mathfrak{t}_\mathbf{d}}^\lambda)^\ast$};

	\draw[wei] (-2,0) +(0,0) -- +(0,-2);
	\draw[wei] (2,0) +(0,0) -- +(0,-2);

	\draw[sharp corners] (1.5,0) +(0,0) rectangle +(6.5,-1);

	\node at (5,-0.5) {$\omega_0^{\mathbf{d}}\cdot\tau$};
	\end{tikzpicture}
	\]
	Where the shaded part at the top and bottom are by construction the diagrams $A_{\mathfrak{t}_\mathbf{d}}^\lambda$ and $(A_{\mathfrak{t}_\mathbf{d}}^\lambda)^\ast$ and the boxed areas in the middle correspond to the braid diagrams of the permutations $\sigma$ and $\omega_0^{\mathbf{d}}\cdot\tau$. So combined the middle part of the diagram corresponds to the permutation $\eta := \omega_0^{\mathbf{d}}\cdot\tau\cdot\sigma$. By \autoref{T: Chap 5 main} and \autoref{L: Both presentations same}, if $\eta = \omega_0^{\mathbf{d}}$ then the inner product is $\pm 1$. Now there are three different cases:
	\begin{enumerate}
		\item If $\tau \neq \sigma^{-1}$ and $\ell(\omega_0^{\mathbf{d}}\cdot\tau) + \ell(\sigma) \neq \ell(\omega_0^{\mathbf{d}})$ then, by \autoref{P: Inner product properties},  the inner product is $0$ because the degrees of the elements, see \autoref{L: degree of tableaux}, do not add up to $0$.
		\item If $\tau \neq \sigma^{-1}$ and $\ell(\omega_0^{\mathbf{d}}\cdot\tau) + \ell(\sigma) = \ell(\omega_0^{\mathbf{d}})$ then $\ell(\eta) < \ell(\omega_0^{\mathbf{d}})$ because $\tau \neq \sigma^{-1}$. However, this means that there is a simple reflection $s\in\mathfrak{S}_{\mathbf{d}}$ such that $s^2$ appears in the product $\eta$. In the diagram this means a double $i$-crossing for some $i\in I$, so by relation \ref{D: Cherednik Algebra}\ref{D: Cherednik double crossings}\ref{D: Cherednik double i crossing} the inner product is $0$. Remember that we can find $s$ because we use Construction B and so the braid diagram does not depend on the choice of reduced expressions for $\sigma$ and $\omega_0^{\mathbf{d}}\cdot\tau$, as we can use braid relations here; see above \autoref{E: Construction B}.
		\item If $\tau = \sigma^{-1}$ then $\eta = \omega_0^{\mathbf{d}}$ and the inner product is $\pm 1$, by \autoref{T: Chap 5 main} and \autoref{L: Both presentations same}.
	\end{enumerate}
	This yields the desired result.
\end{proof}

\emph{Note} that we made a special choice of diagram for the elements $c_{\mathfrak{t}_\mathbf{d}\cdot\sigma}^\lambda$ and a priori the inner product could be different for a different choice of diagram. By \autoref{C: R is Cellular Alg}, \autoref{D: Cherednik Algebra}\ref{D: Cherednik triple crossings} and \cite[Proof of Theorem 4.2]{Bowman:ManyCellular}, for an arbitrary diagram of $c_{\mathfrak{t}_\mathbf{d}\cdot\sigma}^\lambda$ we get
\begin{equation*}
c_{\mathfrak{t}_\mathbf{d}\cdot\sigma}^\lambda = B_{\mathfrak{t}_\mathbf{d}\cdot\sigma}^\lambda + \sum_{\substack{\sigma'\in\mathfrak{S}_\mathbf{d} \\ \ell(\sigma') < \ell(\sigma)}} k_{\sigma'} c_{\mathfrak{t}_\mathbf{d}\cdot\sigma'}^\lambda
\end{equation*}
modulo terms of more dominant shape, for scalars $k_{\sigma'}\in K$. So \autoref{P: Chap 5 Prop} implies \autoref{P: Diagonal Gram matrix}.

\begin{proposition*}
	Let $\lambda\in\mathscr{U}_n^\ell(\theta)$, $\mathbf{d}$ a decomposition path for $\lambda$, and $\sigma,\tau\in\mathfrak{S}_{\mathbf{d}}$. Then $|\langle c_{\mathfrak{t}_{\mathbf{d}}\cdot\sigma}^\lambda, c_{\mathfrak{t}_{\mathbf{d}}^\diamond\cdot\tau}^\lambda\rangle_\lambda| = \delta_{\sigma^{-1},\tau}$.
\end{proposition*}

\noindent In particular, this establishes \autoref{T: Inner Product}.

\end{document}